\documentclass[12pt]{amsart}
\usepackage{amsmath, color}
\usepackage{amsfonts}
\usepackage{latexsym}
\usepackage{ifpdf}
\usepackage{graphicx,amssymb,lineno}
\usepackage{amssymb}
\usepackage{mathrsfs}
\usepackage{cite}
\usepackage{extarrows}
\usepackage{ulem}
\allowdisplaybreaks[4]
\numberwithin{equation}{section}

\date{}

\newcommand{\id}{{\rm id}}

\newcommand{\ext}{{\rm ext}}
\renewcommand{\exp}{{\rm exp}}
\newcommand{\End}{{\rm{End}\ts}}

\newcommand{\diag}{{\rm diag}}
\newcommand{\non}{\nonumber}

\newcommand{\wh}{\widehat}
\newcommand{\ot}{\otimes}

\newcommand{\al}{\alpha}
\newcommand{\be}{\beta}

\newcommand{\Ga}{\Gamma}
\newcommand{\ep}{\epsilon}

\newcommand{\de}{\delta}

\newcommand{\hra}{\hookrightarrow}

\newcommand{\ts}{\,}
\newcommand{\pr}{^{\tss\prime}}

\newcommand{\tss}{\hspace{1pt}}

\newcommand{\U}{ {\rm U}}

\newcommand{\CC}{\mathbb{C}\tss}

\newcommand{\ZZ}{\mathbb{Z}\tss}

\newcommand{\Lc}{\mathcal{L}}

\newcommand{\Rc}{\mathcal{R}}

\newcommand{\Ec}{\mathcal{E}}
\newcommand{\Fc}{\mathcal{F}}

\newcommand{\Hc}{\mathcal{H}}

\newcommand{\Xc}{\mathcal{X}}

\newcommand{\gl}{\mathfrak{gl}}

\newcommand{\spa}{\mathfrak{sp}}

\newcommand{\e}{\mathfrak{e}}
\newcommand{\f}{\mathfrak{f}}
\newcommand{\g}{\mathfrak{g}}
\newcommand{\h}{\mathfrak h}

\newcommand{\z}{\mathfrak{z}}

\newcommand{\tra}{ {\rm t}}

\newcommand{\Sym}{\mathfrak S}

\newcommand{\Fand}{\qquad\text{and}\qquad}

\numberwithin{equation}{section}

\newtheorem{thm}{Theorem}[section]
\newtheorem{lem}[thm]{Lemma}
\newtheorem{prop}[thm]{Proposition}
\newtheorem{cor}[thm]{Corollary}
\newtheorem{conj}[thm]{Conjecture}

\theoremstyle{definition}
\newtheorem{defin}[thm]{Definition}

\theoremstyle{remark}
\newtheorem{remark}[thm]{Remark}
\newtheorem{example}[thm]{Example}

\newcommand{\bth}{\begin{thm}}
\renewcommand{\eth}{\end{thm}}
\newcommand{\bpr}{\begin{prop}}
\newcommand{\epr}{\end{prop}}
\newcommand{\ble}{\begin{lem}}
\newcommand{\ele}{\end{lem}}
\newcommand{\bco}{\begin{cor}}
\newcommand{\eco}{\end{cor}}
\newcommand{\bde}{\begin{defin}}
\newcommand{\ede}{\end{defin}}
\newcommand{\bex}{\begin{example}}
\newcommand{\eex}{\end{example}}
\newcommand{\bre}{\begin{remark}}
\newcommand{\ere}{\end{remark}}
\newcommand{\bcj}{\begin{conj}}
\newcommand{\ecj}{\end{conj}}

\newcommand{\bal}{\begin{aligned}}
\newcommand{\eal}{\end{aligned}}
\newcommand{\beq}{\begin{equation}}
\newcommand{\eeq}{\end{equation}}
\newcommand{\ben}{\begin{equation*}}
\newcommand{\een}{\end{equation*}}

\newcommand{\bpf}{\begin{proof}}
\newcommand{\epf}{\end{proof}}

\def\beql#1{\begin{equation}\label{#1}}

\begin{document}

\title{R-Matrix presentation of quantum affine algebra in type $A_{2n-1}^{(2)}$}
\author{Naihuan Jing$^{1, 2}$}
\author{Xia Zhang$^1$}
\author{Ming Liu$^1$}
\address{$^1$ School of Mathematics, South China University of Technology}
\address{$^2$ Department of Mathematics, North Carolina State University} 

\email{jing@math.ncsu.edu}
\email{1134319889@qq.com}
\email{mamliu@scut.edu.cn}
\subjclass[]{}
\keywords{}
\thanks{}
\begin{abstract}
In this paper, we give an RTT presentation of the twisted quantum affine algebra of
type $A_{2n-1}^{(2)}$
and show that it is isomorphic to the Drinfeld new realization via the Gauss decomposition of the
L-operators. This provides the first such presentation for twisted quantum affine algebras with
nontrivial central element.
\end{abstract}
%

\thanks{{\scriptsize
\hskip -0.6 true cm MSC (2010): Primary: 17B37; Secondary:
\newline Keywords: Quantum affine algebras, matrix presentations, Drinfeld realization, symmetric functions.
\newline Corresponding author: mliu@math.scut.edu.cn
}}

\maketitle

\section{Introduction}
\vspace{5 mm}

Quantum groups are usually referred to two important classes: the quantum enveloping algebras
introduced by Drinfeld and Jimbo
as certain deformations of universal enveloping algebras of the Kac-Moody Lie algebras
\cite{dr, j:qd} and the Yangian algebras of Drinfeld \cite{dr} associated with simple Lie algebras. In this first presentation, the quantum algebras are defined by generators and Serre relations
associated with the simple Lie algebras and their generalizations.
The second presentation
was given by Faddeev-Reshetikhin-Tacktajan \cite{frt} in the study of
quantum inverse scattering method using the Yang-Baxter equation and the RTT relation.
An important extension was given by Reshetikhin and Semenov-Tian-Shanski \cite{rs:ce} for the quantum current
algebra with nontrivial central extension.
It was shown by Drinfeld that all finite dimensional
representations of the quantum affine algebras and the Yangian algebras are classified by the third presentation of the quantum algebra---Drinfeld's new realization \cite{dr}. A detailed construction of
the isomorphism between the first and the third presentations for the quantum affine algebras
were given in \cite{b:bg, j:dr} for the untwisted types and
in \cite{da:drr, da:dri} (see also \cite{jz:br}) in the twisted types.

The exact equivalence of the Drinfeld realization and the Drinfeld-Jimbo presentation for the quantum affine algebra was proved by Ding and
Frenkel in 1993 with the help of the
Faddev-Reshetikhin-Tacktajan's RTT formalism. The same identification for the Yangian algebra was done later by Brundan and Kleshchev
in 2005 \cite{BK}. The identification of the Yangian algebras of other classical types were
proved recently by two groups of mathematicians: the identification of the RTT formulation and the Drinfeld realization by Jing-Liu-Molev \cite{jlm:ib} and that of the RTT formulation and the Drinfeld's J-presentation by Guay-Regelski-Wendlandt \cite{grw:eb}.
Similar equivalence between the Drinfeld realization and the FRT formalism in quantum (untwisted) affine cases
have also been completed for all other classical types \cite{HM, jlm:iC, jlm:ibC} via the Gassian generators
of the L-operators.
The R-matrix formulation of quantum super loop algebras and the quantum loop algebra of
$A_{2n-1}^{(2)}$
have been studied
in \cite{Z, LP} as well.
However, the RTT formulation of the twisted quantum affine algebras (with nontrivial central element)
is still unclear.



In this paper, we use the R-matrix method to study the quantum twisted affine algebras by explicitly
describing all Drinfeld currents in the most general situation. The method has the advantage of revealing the inner-relationship between
the untwisted and twisted quantum affine algebras,
and we have completed the grand picture of
characterizing quantum affine algebras: Drinfeld-Jimbo, Drinfeld's current algebra and RTT presentations.
Among the three presentations, the RTT presentation can provide the direct comultiplication formulas for quantum root
vectors, and a triangular decomposition mimicking the classical case.

\section{Twisted quantum affine algebra}
\subsection{Drinfeld-Jimbo and Drinfeld presentations}
 Let $\g$ be the finite dimensional complex simple Lie algebra $\mathfrak{sl}_{2n}$ with the Chevalley generators
 $$
e^{\prime}_{i}, f_{i}^{\prime}, h_{i}^{\prime}, \quad i={1}, \cdots, {2n-1},
$$
and let $\alpha'_i, i\in I_0'=\{1, \dots, 2n-1\}$ be the the simple roots, and let $A_0'=(A'_{ij})$ be the
Cartan matrix defined by $\langle\al'_i|\al_j'\rangle=2A'_{ij}$, $i, j\in I'_0$, where we normalize
the invariant bilinear form such that $\langle\al|\al\rangle=4$. 

Let $\sigma$ be the automorphism of $\g$ defined by $\sigma(e'_i)=e'_{2n-i}, \sigma(f_i')=f'_{2n-i}$,
which induces an diagram automorphism $\sigma$ of the Dynkin diagram $\Gamma_0'$ so that 
$$
\sigma\left(\al'_{i}\right)=\al'_{2n-i}, i=1, \cdots, 2n-1.
$$

It is well-known that the datum $\left(\mathfrak{sl}_{2n}, \sigma\right)$ gives rise to the twisted affine Lie algebra
$\widehat{\g}$ of type
$A_{2n-1}^{(2)}$ with the affine Dynkin diagram ${\Gamma}={\Gamma}_0\cup\{0\}$, where the vertices of ${\Gamma}_0$ consists of the $\sigma$-orbits
of the vertices of ${\Gamma}_0'$. More explicitly, $\al_i=\frac12(\al'_i+\al'_{2n-i})$, $i\in I_0=\{1, \ldots, n\}\simeq I'_0/\sigma$ are the
simple roots for the invariant Lie subalgebra $\g_0$. 
The induced invariant form on the root lattice $\Gamma_0$
is given by
\begin{equation}
(\al_i|\al_j)=\frac12\left(\langle \al_i'|\al_j'\rangle+\langle\al_{\sigma(i)}'|\al_j'\rangle\right)
\end{equation}
where $i, j\in I_0$. Subsequently the finite Cartan matrix $A_0$ of $\Gamma_0$ is induced from the Cartan matrix $A'_0$: for $i, j\in I_0$
\beql{Eq:Cartan matrixA0}
A_{{i} {j}}=\frac{2(\al_i|\al_j)}{(\al_i|\al_i)}=2 \frac{\sum\limits_{u \in \mathbb{Z}_2}{A}'_{\sigma^{u}\left({i}\right) {j}}}
{\sum\limits_{ u \in \mathbb{Z}_2 } {A}'_{\sigma^{u}\left({i}\right) {i}}}.
\eeq
In this case, $A_0$ is of type $C_n$.

Let $\theta_0$ be the highest weight of $\g$ as a $\g_0$-module. Then
$\theta_0$ is the sum of the two highest roots, and thus $\theta_0=\al'_1+2\sum_{i=1}^{2n-2}\al_i'+\al'_{2n-1}=\al_1+2\sum_{i=2}^{n-1}\al_i+\al_n$.
Then $\delta=\al_0+\theta_0$ is the canonical
null root for 
the twisted affine Lie algebra $A_{2n-1}^{(2)}$. It has the simple roots $\al_i$, $i\in I=\{0, \ldots, n\}$, then
\eqref{Eq:Cartan matrixA0} holds for the affine root system.



Let $d_i=\frac12(\al_i|\al_i)$, then $(d_0,d_1,\dots, d_n)=(1,1,\dots,1,2)$ diagonalizes $A$: $d_iA_{ij}=d_jA_{ji}.$

Let $q$ be a nonzero complex number, and set $q_i=q^{d_i}$ for $i=0,\dots,n$. We also need the standard $q$-numbers: 
\beql{kq}
[k]_q=\frac{q^k-q^{-k}}{q-q^{-1}}
\eeq
and the Gausian $q$-numbers:
$[k]_q!=\prod_{s=1}^{k}[s]_q,
{k\brack r}_{q}=\frac{[k]_q!}{[r]_q!\ts[k-r]_q!}$.

The \textit{quantum affine algebra} $U_q(\wh{\g})$ \cite{dr, j:qd} is
a unital associative algebra over $\CC(q)$ with generators $E_{\pm\al_i}$ 
and
$k_i^{\pm 1}$ ($i=0,1,\dots,n$) subject to the defining relations:
\ben
k_ik_i^{-1}=k_i^{-1}k_i=1,\qquad k_ik_j=k_ik_j,
\een
\ben
k_iE_{\pm\al_j}k_i^{-1}=q^{\pm A_{ij}}_iE_{\pm\al_j}, \qquad 
[E_{\al_i},E_{-\al_j}]=\delta_{ij}\frac{k_i-k_i^{-1}}{q_i-q_i^{-1}},
\een
\ben
\bal
\sum_{r=0}^{1-A_{ij}}(-1)^r{{1-A_{ij}}\brack{r}}_{q_i}
(E_{\pm\al_i})^rE_{\pm\al_j}(E_{\pm\al_i})^{1-A_{ij}-r}&=0,\qquad\text{if}\quad i\ne j.
\eal
\een

\bde \cite{d:nr}
The \textit{Drinfeld realization} of the quantum affine algebra $U^{Dr}_q(\wh{\g})$ of type $A^{(2)}_{2n-1}$
is the $\mathbb{C}(q)$ unital associative algebra generated by
$\tilde{x}_{{i},m}^{\pm}$, $\tilde{a}_{{i},l}$, $\tilde{k}_{{i}}^{\pm}$ and $q^{\pm c/2}$ for $i=1,\dots,2n-1$ and
$m,l\in\ZZ$ with $l\ne 0$, subject to the following defining relations:
the elements $q^{\pm c/2}$ are central,
\beq
\tilde{k}_{{2n-i}}=\tilde{k}_{{i}},~~ \tilde{a}_{{2n-i},l}=(-1)^l\tilde{a}_{{i},l},~~\tilde{x}_{{2n-i},m}^{\pm}=(-1)^{m}\tilde{x}_{{i},m}^{\pm},
\eeq
\beq
\tilde{k}_{{i}}\tilde{k}_{{i}}^{-1}=1=\tilde{k}_{{i}}^{-1}\tilde{k}_{{i}},~~ \tilde{k}_{{i}}\tilde{k}_{{j}}=\tilde{k}_{{j}}\tilde{k}_{{i}},
\eeq
\beq
\tilde{k}_{{i}}\tilde{x}_{{j},m}^{\pm}=q^{\pm\sum\limits_{r\in\mathbb{Z}_2}
(\tilde{\al}_{{i}},\tilde{\al}_{\sigma^r ({j})})}
\tilde{x}_{{j},m}^{\pm}\tilde{k}_{{i}}, \qquad
\tilde{k}_{i}\tilde{a}_{{j},l}=\tilde{a}_{{j},l}\tilde{k}_{{i}},
\eeq
\beq
\left[\tilde{a}_{{i}, k}, \tilde{x}_{{j}, l}^{\pm}\right]=\pm \frac{1}{k}\left(\sum_{r \in \mathbb{Z}_2}
\left[k \tilde{A}_{{i},\sigma^{r}({j})} / d_{{i}}\right]_{q_{{i}}} (-1)^{kr}\right)q^{\mp |k|c/2} \tilde{x}_{{j}, k+l}^{ \pm}
\eeq
\beq
\left[\tilde{a}_{{i},k}, \tilde{a}_{{j}, l}\right]=\delta_{k+l, 0} \frac{1}{k}\left(\sum_{r \in \mathbb{Z}_2}\left[k \tilde{A}_{{i},\sigma^{r}({j})} / d_{{i}}\right]_{q_{{i}}} (-1)^{kr}\right)
 \frac{q^{kc}-q^{-kc}}{q_{{j}}-q_{{j}}^{-1}}
\eeq
\beq
\bal
\prod_{s\in \mathbb{Z}_2}\left(u-(-1)^{s}q^{\pm (\tilde{\al}_{{i}},\tilde{\al}_{\sigma^s({j})})}v\right) \tilde{x}_{{i}}^{\pm}(u) \tilde{x}_{{j}}^{\pm}(v)
=\prod_{s\in \mathbb{Z}_2}\left(uq^{\pm (\tilde{\al}_{{i}},\tilde{\al}_{\sigma^s({j})})}-(-1)^{s}v\right) \tilde{x}_{{j}}^{\pm}(v) \tilde{x}_{{i}}^{\pm}(u)
\eal
\eeq
\beq
\left[\tilde{x}_{{i}}^{+}(k), \tilde{x}_{{j}}^{-}(l)\right]=\sum_{s\in \mathbb{Z}_2} \frac{\delta_{\sigma^{s}({i}), {j}} (-1)^{s l}}{q_{{i}}-q_{{i}}^{-1}}\left(q^{\frac{k-l}{2}c} \tilde{\varphi}^+_{{i}}(k+l)-q^{\frac{l-k}{2}c} \tilde{\varphi}^-_{{i}}(k+l)\right)
\eeq
where $\tilde{\varphi}^{\pm}_{{i}}(\pm m)$ 
$\left(m \in \mathbf{Z}_{\geq0}\right)$
are defined by
\ben\quad \sum_{m=0}^{\infty} \tilde{\varphi}^{\pm}_{{i}}(\pm m) z^{\mp m}=\tilde{k}^{\pm 1}_{{i}}\exp \left(\pm \left(q_{\tilde{i}}-q_{{i}}^{-1}\right) \sum_{k=1}^{\infty} \tilde{a}_{{i},\pm k} z^{\mp k}\right);
\een
\beq\bal
\underset{u_1, u_2}{Sym}\{P_{{i}{j}}^{\pm}&(u_{1}, u_{2})(\tilde{x}_{{j}}^{\pm}(v) \tilde{x}_{{i}}^{\pm}(u_{1}) \tilde{x}_{{i}}^{ \pm}(u_{2})
-[2]_{q^{4d_{ij}}}\tilde{x}_{{i}}^{\pm}(u_{1}) \tilde{x}_{{j}}^{\pm}(v) \tilde{x}_{{i}}^{\pm}(u_{2}) \\
&\qquad +\tilde{x}_{{i}}^{\pm}(u_{1}) \tilde{x}_{{i}}^{\pm}(u_{2}) \tilde{x}_{{j}}^{\pm}(v))\}=0, \qquad \mbox{if} \
\tilde{A}_{{i} {j}}=-1, \ \sigma({i}) \neq {j},
\eal\eeq
where $d_{{i}{j}}$ and $P_{{i}{j}}^{\pm}$ are defined as follows: (1)
if $\sigma({i})={i},$ then $d_{{i} {j}}=\frac{1}{2}, P_{{i} {j}}^{ \pm}\left(u_{1}, u_{2}\right)=1$;
(2) if $\tilde{A}_{{i} \sigma({i})}=0$ and $\sigma({j}) \neq {j},$ then
$d_{{i} {j}}=\frac{1}{4}, P_{{i} {j}}^{ \pm}\left(u_{1}, u_{2}\right)=1$;
(3) if ${A}_{{i} \sigma({i})}=0$ and $\sigma({j})={j},$ then
$d_{{i} {j}}=\frac{1}{2}, P_{{i} {j}}^{ \pm}\left(u_{1}, u_{2}\right)
=\frac{u_{1}^{2} q^{ \pm 4}-u_{2}^{2}}{u_{1} q^{ \pm 2}-u_{2}}$;
(4) if $\tilde{A}_{{i} \sigma({i})}=-1,$ then $d_{{i} {j}}=\frac{1}{4}, P_{{i} {j}}^{ \pm}\left(u_{1}, u_{2}\right)=u_{1} q^{ \pm 1}+u_{2}$.
\beq\bal
&\underset{u_1,u_2,u_3}{Sym}\left\{(q^{\frac32} u_{1}^{\mp 1}-[2]_{q^{\frac12}}u_{2}^{\mp 1}+q^{-\frac32} u_{3}^{\mp 1}) \tilde{x}_{{i}}^{ \pm}(u_{1}) \tilde{x}_{{i}}^{ \pm}(u_{2}) \tilde{x}_{{i}}^{\pm}(u_{3})\right\}=0, \tilde{A}_{{i} \sigma({i})}=-1\\
&\underset{u_1,u_2,u_3}{Sym}\left\{(q^{-\frac32} u_{1}^{ \pm 1}-[2]_{q^{\frac12}} u_{2}^{\pm 1}+q^{\frac32} u_{3}^{\pm 1}) \tilde{x}_{{i}}^{\pm}(u_{1}) \bar{x}_{{i}}^{ \pm}(u_{2}) \tilde{x}_{{i}}^{ \pm}(u_{3})\right\}=0, \tilde{A}_{{i} \sigma({i})}=-1.
\eal\eeq
\ede
A simplified Drinfeld presentation of $U^{Dr}_q(\wh{\g})$ in terms of
$A_0$ datum was given in \cite{da:dri}. 
\bde
The \textit{Drinfeld realization} of the quantum affine algebra of type $A^{(2)}_{2n-1}$
is the $\mathbb{C}(q)$ unital associative algebra generated by
$x_{i,m}^{\pm}$, $a_{i,l}$, $k_{i}^{\pm}$ and $q^{\pm c/2}$ for $i=1,\dots,n$ and
$m,l\in\ZZ$ with $l\ne 0$, subject to the following defining relations:
the elements $q^{\pm c/2}$ are central,
\beq
x_{n,2m+1}^{\pm}=0, \tss a_{n,2m+1}=0,
\eeq
\begin{align}
  &k_{i}k_i^{-1}=k_i^{-1}k_{i}=1, \qquad q^{c/2}q^{-c/2}=q^{-c/2}q^{c/2}=1,\\
  &k_ik_j=k_jk_i, \qquad k_i\ts a_{j,k}=a_{j,k}\ts k_i,\qquad
  k_i\ts x_{j,m}^{\pm}\ts k_i^{-1}=q_i^{\pm A_{ij}}x_{j,m}^{\pm}
\end{align}
\begin{align}
  [a_{i,m},a_{j,l}]&=\delta_{m,-l}\ts
  \frac{d_{ij}[mA_{ij}/d_{ij}]_{q_i}}{m}\ts\frac{q^{mc}-q^{-mc}}{q_j-q_j^{-1}},\\[0.4em]
  [a_{i,m}, x_{j,l}^{\pm}]&=\pm \frac{d_{ij}[{m}A_{ij}/d_{ij}]_{q_i}}{{m}}\ts q^{\mp |m|c/2}\ts x^{\pm}_{j,m+l}, \qquad d_{ij}=\max\{d_i,d_j\}\\
  [x_{i,m}^{+},x_{j,l}^{-}]&=\delta_{ij}\frac{q^{(m-l)\ts c/2}\ts\varphi^+_{i,m+l}
  -q^{-(m-l)\ts c/2}\ts\varphi^-_{i,m+l}}{q_i-q_i^{-1}}\delta_{d_j\mid l}\\
  x^{\pm}_{i,m+d_{ij}}x^{\pm}_{j,l}&-q_i^{ A_{ij}}x^{\pm}_{j,l}x^{\pm}_{i,m+d_{ij}}=
  q_i^{A_{ij}}x^{\pm}_{i,m}x^{\pm}_{j,l+d_{ij}}-x^{\pm}_{j,l+d_{ij}}x^{\pm}_{i,m},
\end{align}
\ben
  \sum_{\pi\in \Sym_{r}}\sum_{l=0}^{r}(-1)^l{{r}\brack{l}}_{q_i}
  x^{\pm}_{i,s_{\pi(1)}}\dots x^{\pm}_{i,s_{\pi(l)}}
  x^{\pm}_{j,m}x^{\pm}_{i,s_{\pi(l+1)}}\dots x^{\pm}_{i,s_{\pi(r)}}=0, \quad i\neq j, \quad A_{ij}\in \{0,-1\},
\een
where $r=1-A_{ij}$.
\ben
\bal
\sum_{\pi \in \mathfrak S_{2}}& \left(q(x_{j, s}^{\pm} x_{i, r_{\pi(1)}\pm 1}^{ \pm} x_{i, r_{\pi(2)}}^{\pm}-
[2]_{q^{2}} x_{i, r_{\pi(1)}\pm 1}^{\pm} x_{j, s}^{ \pm} x_{i,r_{\pi(2)}}^{\pm}+x_{i, r_{\pi(1)}\pm 1}^{\pm}
 x_{i, r_{\pi(2)}}^{\pm} x_{j, s}^{\pm})+\right.\\
& \left.q^{-1}(x_{j, s}^{ \pm} x_{i, r_{\pi(1)}}^{ \pm} x_{i, r_{\pi(2)}\pm1}^{\pm}
-[2]_{q^{2}} x_{i,r_{\pi(1)}}^{\pm} x_{j,s}^{\pm} x_{i,r_{\pi(2)}\pm 1}^{\pm}
+x_{i, r_{\pi(1)}}^{\pm} x_{i, r_{\pi(2)} \pm 1}^{ \pm} x_{j, s}^{\pm})\right)=0
\eal
\een
for $A_{ij}=-2$.
The elements
$\varphi^+_{i,m}$ and $\varphi^-_{i,-m}$ with $m\in \ZZ_+$ are defined by
\begin{align}
\label{psiiu}
\varphi^{\pm}_i(u)&:=\sum_{m=0}^{\infty}\varphi^{\pm}_{i,m}u^{\mp m}=k_i^{\pm 1}\ts\exp\big(\pm (q_i-q_i^{-1})
\sum_{s=1}^{\infty}a_{i,\pm s}u^{\mp s}\big) 
\end{align}
whereas $\varphi^+_{i,m}=\varphi^-_{i,-m}=0$ for $m<0$.

\ede
\subsection{Isomorphism between Drinfeld-Jimbo and Drinfeld presentations}\label{subsec:isoDJD}

By using the braid group action, the set of generators of the algebra $U_q(A^{(2)}_{2n-1})$
can be extended to the set of affine root vectors of the form $E_{\al+kd_{\al}\de}$, $F_{\al+kd_{\al}\de}$,
$E_{(kd_i\de,i)}$ and $F_{(kd_i\de,i)}$, where $\al$ runs over the positive roots of $\mathfrak{sp}_{2n}$,
$d_{\al}=(\al,\al)/2$
and $\de$ is the basic imaginary root; see \cite{b:bg, bcp:ac} for details.
The root vectors are used in the explicit isomorphism between the
Drinfeld--Jimbo presentation of the algebra $U_q(A_{2n-1}^{(2)})$
and the ``new realization'' of Drinfeld which goes back to \cite{d:nr}, while detailed arguments
were given by Beck~\cite{b:bg}; see also \cite{jz:br, bcp:ac}.
In particular, for the Drinfeld presentation of the algebra $U_q(\wh{\g}_{N})$
given in the Introduction, we find that
the isomorphism between these presentations is given by
\begin{alignat}{2}
\non
x^{+}_{i,kd_i}&\mapsto o(i)^kE_{\al_i+kd_i\de}, \qquad &x^{-}_{i,-kd_i}
&\mapsto o(i)^kF_{\al_i+kd_i\de},\qquad\qquad  k\geqslant 0,\\
\non
x^{+}_{i,-kd_i}&\mapsto -o(i)^kF_{-\al_i+kd_i\de}\ts k_i^{-1}q^{kd_ic},
\qquad &x^{-}_{i,kd_i}
&\mapsto -o(i)^kq^{-kd_ic}\ts k_i\ts E_{-\al_i+kd_i\de},\quad  k> 0,\\
a_{i,kd_i}&\mapsto o(i)^kq^{-kd_ic/2}E_{(kd_i\de,i)}, \qquad &a_{i,-kd_i}
&\mapsto o(i)^k\tss F_{(kd_i\de,i)}q^{kd_ic/2},\qquad\qquad k> 0,
\non
\end{alignat}
where $o: \{1,2,\dots,n\}\rightarrow \{\pm 1\}$ is a map such that $o(i)=-o(j)$ whenever $A_{ij}<0$
and furthermore $o(i)=1$ when $A_{ij}=-2$, thus $o(n)=-1$. The pattern of $o(i)$ is $((-1)^n, (-1)^{n-1}, \cdots, +, -1)$
or $o(i)=(-1)^{n+1-i}$.

\subsection{Universal R-matrix}
We recall the explicit formulas for the
universal $R$-matrix for the algebra $U_q(\wh{\g})$
obtained by Khoroshkin and Tolstoy~\cite{kt:ur}
and Damiani~\cite{d:rm,da:Un}.

Recall the Cartan matrix $A_0$ for $\spa_{2n}$, the underlying simple Lie algebra in $A^{(2)}_{2n-1}$, defined by \eqref{Eq:Cartan matrixA0}.
The symmetrized Cartan matrix $B_0=DA_0$ is given by
$B_{ij}=(\al_i,\al_j), 1\leq i, j\leq n$.
Let $\tilde{B}=[\tilde{B}_{ij}]$ be the inverse matrix $B_0^{-1}$, then its
entries $\tilde{B}$ are given by
\beql{Bij}
\tilde{B}_{ij}=\begin{cases}
n/4\qquad&\text{for}\quad i=j=n,\\
   j/2\qquad&\text{for}\quad i=n>j,\\
   j\qquad&\text{for}\quad n>i\geqslant j.
  \end{cases}
\eeq
In order to describe the explicit formula for universal R-matrix given in 
\cite{d:rm},
we introduce the matrix $Z^k$ for $k\in \mathbb{Z}_{>0}$ as following:
\beql{Eq:Zij}
z^{k}_{ij}=\begin{cases}
\frac{[n-i]_{q^k}[j]_{q^k}}{[n]_{q^k}}\qquad&\text{for}\quad \text{$k$ is odd},\\
   \frac{[n]_{q^k}}{[2]_{q^{nk}}}\qquad&\text{for}\quad i=j=n, \text{and $k$ is even},\\
   \frac{(-1)^{k/2}2[j]_{q^k}}{[2]_{q^{nk}}}\qquad&\text{for}\quad n=i> j,\text{and $k$ is even},\\
   \frac{[2]_{q^{(n-i)k}}[j]_{q^k}}{[2]_{q^{kn}}}\qquad&\text{for}\quad \text{other cases}.
  \end{cases}
\eeq
Note that $z_{ij}^{k}d_j=z_{ji}^{k}d_i$, so $Z^k$ is symmetrizable.
We adjoin the degree element $d$ to $U_q(A_{2n-1}^{(2)})$ to form the quantum affine algebra
 algebra $\widetilde{U}_q(A_{2n-1}^{(2)})$ with the relations: 
\ben
[d,k_i]=0, \qquad [d,E_{\pm\al_i}]=\pm\delta_{i,0}E_{\pm\al_i}  
\een

For a formal variable $u$, define an automorphism $D_u$ of
the algebra $\widetilde{U}_q(A_{2n-1}^{(2)})\otimes \CC[u,u^{-1}]$
by
\ben
D_u(x^{\pm}_{i,k})=u^{k}x^{\pm}_{i,k},\qquad D_u(a_{i,k})=u^{k}a_{i,k},
\een
for $i=1,\dots, n-1$,
and
\ben
D_u(x^{\pm}_{n,2k})=u^{2k\pm 1}x^{\pm}_{n,2k},\qquad D_u(a_{n,2k})=u^{2k}a_{n,2k}
\qquad D_u(k_{i})=k_i,\qquad D_u(d)=d.
\een

The {\it universal $R$-matrix} of the quantum affine algebra is an element $\mathfrak{R}$ in the completed tensor product
$\widetilde{U}_q(A_{2n-1}^{(2)})\ts\widehat{\otimes}\ts \widetilde{U}_q(A_{2n-1}^{(2)})$ (see Drinfeld~\cite{d:ac})
satisfying certain conditions.
The conditions imply that $\mathfrak{R}$ is a solution of the Yang-Baxter equation
\ben
\mathfrak{R}_{12}\mathfrak{R}_{13}\mathfrak{R}_{23}=\mathfrak{R}_{23}\mathfrak{R}_{13}\mathfrak{R}_{12},
\een
over $\widetilde{U}_q(A_{2n-1}^{(2)})^{\hat{\otimes} 3}$.
Here
$\mathfrak{R}_{ij}$ means that if $A\in V\otimes V$ then $A_{23}=I\otimes A\in V^{\otimes 3}$ etc.
The explicit formula for $\mathfrak{R}$ is given in the $\hbar$-adic settings and we set $q=\exp(\hbar)\in \CC[[\hbar]]$
and regard the quantum affine algebra over $\CC[[\hbar]]$.
Introduce elements $h_1,\dots,h_n$ by
setting $k_i=\exp(\hbar h_i)$.


We will work with the parameter-dependent $R$-matrix defined by
\ben
\Rc(u)=(D_u\otimes \id)\ts \mathfrak{R}\ts q^{c\otimes d+d\otimes c}.
\een
It satisfies the Yang--Baxter equation in the form
\beql{MYBE}
\Rc_{12}(u)\Rc_{13}(uvq^{-c_2})\Rc_{23}(v)
=\Rc_{23}(v)\Rc_{13}(uvq^{c_2})\Rc_{12}(u)
\eeq
where $c_2=1\otimes c\otimes 1$; cf. \cite{fri:qa}.

The R-matrix $\Rc(u)$ has the following triangular decomposition (cf. \cite{da:Un}): 
\beql{rdec}
\Rc(u)=\Rc^{>0}(u)\Rc^{0}(u)\Rc^{<0}(u),
\eeq
where
\ben
\bal
  \Rc^{>0}(u) &= \prod_{\alpha\in \Delta_+}\prod_{k\geqslant 0}
 \exp_{q_{\alpha}}\big((q_{\alpha}^{-1}-q_{\alpha})D_u(E_{\alpha+kd_{\alpha}\delta})\otimes F_{\alpha+kd_{\alpha}\delta}\big), \\
 \Rc^{<0}(u) &= T^{-1}\prod_{\alpha\in \Delta_+}\prod_{k> 0}
 \exp_{q_{\alpha}}\big((q_{\alpha}^{-1}-q_{\alpha})D_u(E_{-\alpha+kd_{\alpha}\delta})\otimes F_{-\alpha+kd_{\alpha}\delta}\big)\ts T\\
  \Rc^{0}(u)&=\exp\Big(
  \sum\limits_{k>0}\sum\limits_{i,j=1}^{n}\frac{1}{d_{i}}\frac{(q_i^{-1}-q_i)(q_j^{-1}-q_j)}{q^{-1}-q}
  \frac{k}{[k]_q}{z}^{k}_{ij}o(i)^{k}o(j)^{k}D_u(E_{kd_i\delta,i})\otimes F_{kd_j\delta,j}\Big)\ts T,
\eal
\een
where $T=\exp(-\hbar \sum_{i, j}^n\tilde{B}_{ij}\tss h_i\otimes h_j)$. 

A straightforward calculation verifies the following formulas for the vector
representation of the quantum affine algebra. As before, we denote by $e_{ij}\in\End\CC^{2n}$
the standard matrix units. Here $V=\sum_{i=1}^n\CC v_i\oplus\CC v_{i'}$ and $e_{ij}v_k=\delta_{jk}e_i$ for $i, j\in\{1, \cdots, n, n', \cdots, 1'\}$.

\bpr\label{prop:FFRep}
The mappings $q^{\pm c/2}\mapsto 1$,
\ben
\bal
    x^{+}_{ik} &\mapsto -q^{-ik}e_{i+1,i}+\xi^{k}q^{ik}e_{i',(i+1)'}, \qquad x^{-}_{ik} \mapsto -q^{-ik}e_{i,i+1}+\xi^{k}q^{ik}e_{(i+1)',i'},\\
    a_{ik} &\mapsto \frac{[k]_{q_i}}{k}\big(q^{-ik}(q^{-k}e_{i+1,i+1}-q^{k}e_{ii})
    +\xi^{k}q^{ik}(q^{-k}e_{i'i'}-q^ke_{(i+1)'(i+1)'})
    \big)\\
    k_i&\mapsto q(e_{i+1,i+1}+e_{i',i'})+q^{-1}(e_{i,i}+e_{(i+1)'(i+1)'})
+\sum\limits_{j\neq i,i+1,i',(i+1)'}e_{jj}
\eal
\een
    for $i=1,\dots,n-1$, and
\ben
\bal    x^{+}_{n,2k} &\mapsto -q_{n}^{-nk}e_{n+1,n}, \qquad x^{-}_{n,2k} \mapsto -q_{n}^{-nk}e_{n,n+1},\\
    a_{n,2k} &\mapsto \frac{[k]_{q_n}}{k}\big(q_n^{-nk}(q_n^{-k}e_{n+1,n+1}-q_n^{k}e_{nn})
    \big)\\
    k_n&\mapsto q^2e_{n+1,n+1}+q^{-2}e_{n,n}+\sum\limits_{j\neq n,n+1}e_{jj},
\eal
\een
define a representation $\pi_V: U_q(A_{2n-1}^{(2)})\to\End V$
of $U_q(A_{2n-1}^{(2)})$ on the vector space $V=\mathbb{C}^{2n}$.
\qed
\epr

It follows from the results of \cite{fri:qa} that
\beql{Rm}
R(u)\triangleq (\pi_{V}\otimes \pi_{V})\ts\Rc(u)=f(u)\overline{R}(u),
\eeq
where
\beql{fu}
f(u)=\xi q^{-2}\prod_{r=0}^{\infty} \frac{\left(1-u \xi^{2 r}\right)\left(1-u q^{-2} \xi^{2 r+1}\right)\left(1-u q^{2} \xi^{2 r+1}\right)\left(1-u \xi^{2 r+2}\right)}{\left(1-u \xi^{2 r+1}\right)^2\left(1-u q^{2} \xi^{2 r+2}\right)\left(1-u q^{-2} \xi^{2 r}\right)},
\eeq
\beql{rbar}
\bal
\overline{R}(u)&=\sum_{i=1}^{2n}e_{ii}\otimes e_{ii}+\frac{u-1}{qu-q^{-1}}\sum_{i\neq j,j'}e_{ii}\otimes e_{jj}
+\frac{q-q^{-1}}{qu-q^{-1}}\sum_{i>j,i\neq j'}e_{ij}\otimes e_{ji}\\
&+\frac{(q-q^{-1})u}{qu-q^{-1}}\sum_{i<j, i\neq j'}e_{ij}\otimes e_{ji}+
\frac{1}{(u-q^{-2})(u-\xi)}\sum_{i,j=1}^{2n} a_{ij}(u)e_{i'j'}\otimes e_{ij},
\eal
\eeq
where $\xi=-q^{-2n}$,
$(\bar{1},\bar{2},...\bar{2n})=(n-1,n-2,...,0,0,...,-n+1)$ and
\beq
a_{ij}(u)=
\left\{
  \begin{array}{ll}
    (q^{-2}u-\xi)(u-1), & \hbox{$i=j$;} \\
    (q^{-2}-1)(q^{\bar{i}-\bar{j}}\xi(u-1)-\de_{ij'}(u-\xi)), & \hbox{$i<j$;} \\
    (q^{-2}-1)u(q^{\bar{i}-\bar{j}}(u-1)-\de_{ij'}(u-\xi)), & \hbox{$i>j$.}
  \end{array}
\right.
\eeq

Introduce the $L$-{\it operators} in $U_q(A^{(2)}_{2n-1})\ot End V$ by the formulas
\beq\label{Eq:Ltilde-operator}
\bal
  \tilde{L}^{+}(u) &= (\id\otimes \pi_{V})\ts \Rc_{21}(uq^{c/2}),\qquad
  \tilde{L}^{-}(u) = (\id\otimes \pi_{V})\ts \Rc_{12}(u^{-1}q^{-c/2})^{-1}
\eal
\eeq
viewed as linear operator over the ring $U_q(A^{(2)}_{2n-1})$ or matrices with entries in $U_q(A^{(2)}_{2n-1})$.

The Yang-Baxter equation \eqref{MYBE} immediately implies the following proposition:

\bpr \label{e:RTT4Ltilde}
The following relations hold
in $U_q(A^{(2)}_{2n-1})\ot EndV^{\ot 2}$ ($u_{\pm}=uq^{\pm c/2}$):
\ben
\bal
R(u/v)\tilde{L}^{\pm}_1(u)\tilde{L}^{\pm}_2(v) &= \tilde{L}^{\pm}_2(v)\tilde{L}^{\pm}_1(u)R(u/v), \\[0.4em]
  R(u_{+}/v_{-})\tilde{L}^{\pm}_1(u)\tilde{L}^{\mp}_2(v) &= \tilde{L}^{\mp}_2(v)\tilde{L}^{\pm}_1(u)R(u_{-}/v_{+}).
\eal
\een
\epr

\section{Extended quantum affine algebra $U^{ext}_q(A^{(2)}_{2n-1})$}

\subsection{Extended twisted quantum affine algebras}
We will introduce the extended twisted quantum affine
algebra $U^{ext}_q(A^{(2)}_{2n-1})$, which contains
$U_q(A^{(2)}_{2n-1})$ as a subalgebra.

\bde\label{def:eqaa}
The {\it extended quantum affine algebra} $U^{\ext}_q(A^{(2)}_{2n-1})$ is an associative algebra
with generators $X^{\pm}_{i,k}$, $h^{+}_{j,m}$, $h^{-}_{j,-m}$ and $q^{c/2}$, where $i=1,\dots,n$, $k\in\ZZ$, $j=1,\dots,n+1$ and $m\in \ZZ_{+}$.
The defining relations are written in terms of generating functions in a formal
variable $u$:
\ben
X^{\pm}_i(u)=\sum_{k\in\ZZ}X^{\pm}_{i,k}\ts u^{-k},\qquad
h^{\pm}_i(u)=\sum_{m=0}^{\infty}h^{\pm}_{i,\pm m}\ts u^{\mp m},
\een
and take the following form: The element $q^{c/2}$ is central and invertible,
\begin{align*}
h_{i,0}^{+}h_{i,0}^{-}=h_{i,0}^{-}h_{i,0}^{+}=h_{n,0}^{+}h_{n+1,0}^{+}=1,\quad
x_{n,2k+1}^{\pm}=h^{\pm}_{n,2k+1}=h^{\pm}_{n+1,2k+1}=0,
\end{align*}
\begin{align}
\label{eq:hihj}
h^{\pm}_i(u)h^{\pm}_j(v)&=h^{\pm}_j(v)h^{\pm}_i(u),\\[0.4em]
\label{eq:hihipm}
f\Big((u_{+}/v_{-})^{\pm 1}\Big)h^{\pm}_i(u)h^{\mp}_i(v)&=f\Big((u_{-}/v_{+})^{\pm 1}\Big)h^{\mp}_i(v)h^{\pm}_i(u),
\\
f\Big((u_{+}/v_{-})^{\pm 1}\Big)\frac{u_{\pm}-v_{\mp}}{qu_{\pm}-q^{-1}v_{\mp}}h^{\pm}_i(u)h^{\mp}_j(v)&=
f\Big((u_{-}/v_{+})^{\pm 1}\Big)\frac{u_{\mp}-v_{\pm}}{qu_{\mp}-q^{-1}v_{\pm}}h^{\mp}_j(v)h^{\pm}_i(u)
\label{hihjpm}
\end{align}
for $i<j$ and $i\neq n$, and
\beql{eq:hnhn'pm}
\bal
f\Big((\frac{u_{+}}{v_{-}})^{\pm 1}\Big)\frac{u_{\pm}^2-v_{\mp}^{2}}{q^2u_{\pm}^2-q^{-2}v_{\mp}^2}h^{\pm}_n(u)h_{n+1}^{\mp}(v)=
f\Big((\frac{u_{-}}{v_{+}})^{\pm 1}\Big)\frac{u_{\mp}^2-v_{\pm}^{2}}{q^2u_{\mp}^2-q^{-2}v_{\pm}^2}h_{n+1}^{\mp}(v)h^{\pm}_{n}(u).
\eal
\eeq
The relations
involving $h^{\pm}_i(u)$ and $X_{j}^{\pm}(v)$ are
\ben
\bal
h_{i}^{\pm}(u)X_{j}^{+}(v)h_{i}^{\pm}(u)^{-1}&
=\frac{u/v_{\pm}-1}{q^{(\ep_i,\alpha_j)}u/v_{\pm}-q^{-(\ep_i,\alpha_j)}}
X_{j}^{+}(v),\\[0.4em]
h_{i}^{\pm}(u)^{-1}X_{j}^{-}(v)h_{i}^{\pm}(u)&
=\frac{u/v_{\mp}-1}{q^{(\ep_i,\alpha_j)}u/v_{\mp}-q^{-(\ep_i,\alpha_j)}}
X_{j}^{-}(v)
\eal
\een
for $i=1,\dots,n$, $j=1,\dots, n-1$
together with
\ben
\bal
h_{i}^{\pm}(u)X_{n}^{+}(v)h_{i}^{\pm}(u)^{-1}&
=\frac{(u/v_{\pm})^2-1}{q^{(\ep_i,\alpha_j)}(u/v_{\pm})^2-q^{-(\ep_i,\alpha_j)}}
X_{n}^{+}(v),\\[0.4em]
h_{i}^{\pm}(u)^{-1}X_{n}^{-}(v)h_{i}^{\pm}(u)&
=\frac{(u/v_{\mp})^2-1}{q^{(\ep_i,\alpha_j)}(u/v_{\mp})^2-q^{-(\ep_i,\alpha_j)}}
X_{n}^{-}(v)
\eal
\een
for $i=1,\dots n$
and
\ben
\bal
h_{n+1}^{\pm}(u)X_{n}^{+}(v)h_{n+1}^{\pm}(u)^{-1}&
=\frac{(u/v_{\pm})^{2}-1}{q^{-2}(u/v_{\pm})^2-q^{2}}
X_{n}^{+}(v),\\[0.4em]
h_{n+1}^{\pm}(u)^{-1}X_{n}^{-}(v)h_{n+1}^{\pm}(u)&
=\frac{(u/v_{\mp})^2-1}{q^{-2}(u/v_{\mp})^2-q^{2}}
X_{n}^{-}(v)
\eal
\een
and
\ben
\bal
h_{n+1}^{\pm}(u)X_{n-1}^{+}(v)h_{n+1}^{\pm}(u)^{-1}&=
\frac{q^{-1}u/v_{\pm}+q}{u/v_{\pm}+1}X^{+}_{n-1}(v),\\
h_{n+1}^{\pm}(u)^{-1}X_{n-1}^{-}(v)h_{n+1}^{\pm}(u)&=
\frac{q^{-1}u/v_{\mp}+q}{u/v_{\mp}+1}X^{-}_{n-1}(v),
\eal
\een
while
\ben
\bal
h_{n+1}^{\pm}(u)X_{i}^{+}(v)=X_{i}^{+}(v)h_{n+1}^{\pm}(u),\qquad
h_{n+1}^{\pm}(u)X_{i}^{-}(v)=X_{i}^{-}(v)h_{n+1}^{\pm}(u),
\eal
\een
for $1\leqslant i\leqslant n-2$. For the relations involving $X^{\pm}_i(u)$ we have
\ben
(u-q^{\pm (\alpha_i,\alpha_j)}v)X_{i}^{\pm}(uq^i)X_{j}^{\pm}(vq^j)
=(q^{\pm (\alpha_i,\alpha_j)}u-v) X_{j}^{\pm}(vq^j)X_{i}^{\pm}(uq^i)
\een
for $i,j=1,\dots,n-1$;
\ben
(u^2-q^{\pm (\alpha_i,\alpha_n)}v^2)X_{i}^{\pm}(uq^i)X_{n}^{\pm}(vq^n)
=(q^{\pm (\alpha_i,\alpha_n)}u^2-v^2) X_{n}^{\pm}(vq^n)X_{i}^{\pm}(uq^i)
\een
for $i=1,\dots,n$
and
\ben
[X_i^{+}(u),X_j^{-}(v)]=
\delta_{ij}(q_i-q_i^{-1})\Big(\delta\big(\frac{uq^{-c}}v\big)h_i^{-}(v_+)^{-1}h_{i+1}^{-}(v_+)
-\delta\big(\frac{uq^{c}}v\big)h_i^{+}(u_+)^{-1}h_{i+1}^{+}(u_+)\Big)
\een
for $j\neq n$
and
\ben
[X_i^{+}(u),X_n^{-}(v)]=
\delta_{in}(q_n-q_n^{-1})\Big(\delta\big((\frac{uq^{-c}}v)^2\big)q^{c}h_i^{-}(v_+)^{-1}h_{i+1}^{-}(v_+)
-\delta\big((\frac{uq^{c}}v)^2\big)q^{-c}h_i^{+}(u_+)^{-1}h_{i+1}^{+}(u_+)\Big)
\een
where the formal delta function $\de(u)=\sum_{r\in\ZZ}\tss u^r$ and the
Serre relations are
\beql{serrex}
\sum_{\pi\in \Sym_{r}}\sum_{l=0}^{r}(-1)^l{{r}\brack{l}}_{q_i}
  X^{\pm}_{i}(u_{\pi(1)})\dots X^{\pm}_{i}(u_{\pi(l)})
  X^{\pm}_{j}(v)\tss X^{\pm}_{i}(u_{\pi(l+1)})\dots X^{\pm}_{i}(u_{\pi(r)})=0,
\eeq
which hold for all $i\neq j$ and $r=1-A_{ij}$. For $A_{ij}=-1$, the relation is
\beq
\bal
&\sum_{\sigma\in \mathfrak{S}_2}\sigma\Big((q^2u_1+u_2)(X^{\pm}_n(v)X_{n-1}^{\pm}(u_1)X_{n-1}^{\pm}(u_2)\\
&-[2]_{q^2}X_{n-1}^{\pm}(u_1)X_{n}^{\pm}(v)X_{n-1}^{\pm}(u_2)
+X_{n-1}^{\pm}(u_1)X_{n-1}^{\pm}(u_2)X_{n}^{\pm}(v))\Big)=0
\eal
\eeq
\ede

Introduce two formal series $z^{\pm}(u)$ 
in $u^{\pm 1}$ 
with coefficients in $U^{\ext}_q(A_{2n-1}^{(2)})$ by ($\xi=-q^{-2n}$)
\beql{zpm}
{z}^{\pm}(u)=\prod_{i=1}^{n-1}{h}^{\pm}_{i}(u\xi q^{2i})^{-1}
\prod_{i=1}^{n}{h}^{\pm}_{i}(u\xi q^{2i-2}){h}^{\pm}_{n+1}(u).
\eeq

\bpr\label{prop:zu}
The coefficients of $z^{\pm}(u)$ are central elements of $U^{\ext}_q(A^{(2)}_{2n-1})$.
\epr

\bpf
We take $z^+(u)$ to show the argument. 
By \eqref{eq:hihj} we have
$z^{+}(u)\ts h_j^{+}(v)=h_j^{+}(v)\ts z^{+}(u)$ for all $j=1,\dots,n+1$.
Also it follows from \ref{def:eqaa}
that $z^{+}(u)\ts X_j^{\pm}(v)=X_j^{\pm}(v)\ts z^{+}(u)$ for $j=1,\dots, n$.
We now check $z^{+}(u)\ts h_{n+1}^{-}(v)=h_{n+1}^{-}(v)\ts z^{+}(u)$ in detail by using
\eqref{fu}. 
By \eqref{eq:hihipm} we have
\ben
\bal
z^{+}(u)h^{-}_{n+1}(v)=f(u_{-}/v_{+})\ts f(u_{+}/v_{-})^{-1}
\prod_{i=1}^{n-1}{h}^{+}_{i}(u\xi q^{2i})^{-1}
\prod_{i=1}^{n}{h}^{+}_{i}(u\xi q^{2i-2})h^{-}_{n+1}(v)h_{n+1}^{+}(u).
\eal
\een
Due to the commutativity \eqref{eq:hihj}, the last product can be rewritten as
\ben
\prod_{i=1}^{n-2}h^{+}_{i}(u\xi q^{2i})^{-1}h^{+}_{i+1}(u\xi q^{2i})h^{+}_{1}(u\xi)
h^{+}_{n-1}(-uq^{-2})^{-1}h_{n}^{+}(-uq^{-2})
h^{-}_{n+1}(v)h_{n+1}^{+}(u).
\een
Now \eqref{hihjpm} and \eqref{eq:hnhn'pm} imply
\begin{multline}
\non
z^{+}(u)h^{-}_{n+1}(v)=f(u_{-}/v_{+})\ts f(u_{+}/v_{-})^{-1}
{}\times \frac{(q^{-2}u_{-}/v_{+}-1)(u_{+}/v_{-}-1)}{(q^{-2}u_{+}/v_{-}-1)(u_{-}/v_{+}-1)}\\
{}\times\prod_{i=1}^{n-2}h^{+}_{i}(u\xi q^{2i})^{-1}h^{+}_{i+1}(u\xi q^{2i})h^{+}_{1}(u\xi)
h^{-}_{n+1}(v)h^{+}_{n-1}(-uq^{-2})^{-1}h_{n}^{+}(-uq^{-2})
h_{n+1}^{+}(u).
\end{multline}

Now, applying \eqref{hihjpm} agian, we come to the relation
\ben
\bal
&z^{+}(u)h^{-}_{n+1}(v)=
\frac{f(u_{-}\xi/v_{+})}{ f(u_{+}\xi/v_{-})}\frac{f(u_{-}/v_{+})}{ f(u_{+}/v_{-})}\\[0.5em]
&\qquad\qquad \times
\frac{(q^{-2}u_{-}/v_{+}-1)(u_{+}/v_{-}-1)}{(q^{-2}u_{+}/v_{-}-1)(u_{-}/v_{+}-1)}
\frac{u_{-}\xi/v_{+}-1}{qu_{-}\xi/v_{+}-q^{-1}}\frac{qu_{+}\xi/v_{-}-q^{-1}}{u_{+}\xi/v_{-}-1}
\ts h^{-}_{n+1}(v)z^{+}(u).
\eal
\een
Then $z^{+}(u)h^{-}_{n+1}(v)=h^{-}_{n+1}(v)z^{+}(u)$ follows due to:
$f(u)f(u\xi)=(\xi q^{-2})^2\frac{(1-uq^{2}\xi)(1-u)}{(1-uq^{-2})(1-u\xi)}$.
\epf

\bpr\label{prop:embed}
The map $q^{c/2}\mapsto q^{c/2}$, $k_i\mapsto h^{+}_{i0}(h_{i+1,0}^{+})^{-1}$ $(1\leq i\leq n-1)$, $k_n\mapsto (h^{+}_{n,0})^2$ and
\ben
\bal
x^{\pm}_{i}(u)&\mapsto (q_i-q_i^{-1})^{-1}X^{\pm}_i(uq^i),\qquad\quad 1\leq i\leq n\\
\varphi^{\pm}_{i}(u)&\mapsto \begin{cases} h^{\mp}_{i+1}(uq^i)\ts h^{\mp}_{i}(uq^i)^{-1}, & 1\leq i\leq n-1\\
q^{\pm c}h^{\mp}_{n+1}(uq^{n})\ts h^{\mp}_{n}(uq^{n})^{-1}, & i=n\end{cases}
\eal
\een
defines an embedding $\varsigma:U_q(A^{(2)}_{2n-1})\hra U^{\ext}_q(A^{(2)}_{2n-1})$.
\epr

\bpf In terms of generating series, it is straightforward to verify that the map $\varsigma$ defines a homomorphism from
$U_q(A^{(2)}_{2n-1})$ to $U^{\ext}_q(A^{(2)}_{2n-1})$.
To show the injectivity, we will construct another
homomorphism $\varrho:U^{ext}_q(A^{(2)}_{2n-1})\to U_q(A^{(2)}_{2n-1})$ such that
the composition $\varrho\circ\varsigma$ is the identity homomorphism on $U_q(A^{(2)}_{2n-1})$.
We extend $U_q(A^{(2)}_{2n-1})$ by adjoining the square roots $k_n^{\pm 1/2}$
to $U_q(A^{(2)}_{2n-1})$ and keep the same notation
for the extended algebras for the rest of the proof.
First we can verify that the mapping $\rho_1: U^{\ext}_q(A^{(2)}_{2n-1})\to U^{\ext}_q(A^{(2)}_{2n-1})$ by
\beq
\bal
&X^{\pm}_j(u)\mapsto X^{\pm}_j(u),\\
&h_i^{\pm}(u)\mapsto \prod_{m=0}^{\infty}z^{+}(u\xi^{-2m-1})z^{+}(u\xi^{-2m-2})^{-1}h_i^{\pm}(u)
\eal
\eeq
is an automorphism of $U^{\ext}_q(A^{(2)}_{2n-1})$.

By setting $\overline{\varphi}^{\pm}_{i}(u)=k_i^{\mp 1}\tss\varphi^{\pm}_{i}(u)$ and
\begin{align*}h^{\pm}_0(u)&=\prod_{m=0}^{\infty} \prod_{j=1}^{n-1}\overline{\varphi}^{\mp}_{j}(u\xi^{-2m}q^j)^{-1}
\overline{\varphi}^{\mp}_{j}(u\xi^{-2m-1}q^j)
\overline{\varphi}^{\mp}_{j}(u\xi^{-2m-1}q^{-j})^{-1}\overline{\varphi}^{\mp}_{j}(u\xi^{-2m-2}q^{-j})\\
& \hskip 3cm \times \prod_{m=0}^{\infty} \overline{\varphi}^{\mp}_{n}(-u\xi^{-2m}q^{n})^{-1}
\overline{\varphi}^{\mp}_{n}(-u\xi^{-2m-1}q^{n})
\end{align*}
we define the map $\rho_2:U^{\ext}_q(A^{(2)}_{2n-1})\to U_q(A^{(2)}_{2n-1})$ by
\ben
X^{\pm}_i(u)\mapsto
\ts (q_i-q_i^{-1})x^{\pm}_{i}(uq^{-i}), ~~i=1,\cdots n,
\een
while for $1\leq i\leq n-1$
\begin{align*}
h^{\pm}_i(u)&\mapsto h^{\pm}_0(u)\prod_{j=1}^{i-1}\overline{\varphi}^{\mp}_{j}(uq^{-j})\times \prod_{j=i}^{n-1}k_i^{\pm}k^{\pm1/2}_n,\\
h^{\pm}_n(u)&\mapsto h^{\pm}_0(u)\prod_{j=1}^{n-1}\overline{\varphi}^{\mp}_{j}(uq^{-j})\times k^{\pm1/2}_nq^{\mp c},\\
h^{\pm}_{n+1}(u)&\mapsto h^{\pm}_0(u) \prod_{j=1}^{n-1}\overline{\varphi}^{\mp}_{j}(uq^{-j})
{ \overline{\varphi}^{\mp}_{n}(uq^{-n})\times k^{\mp 1/2}_n.}
\end{align*}
Direct calculation
shows that the map $\rho_2$ defines a homomorphism from $U^{ext}_q(A^{(2)}_{2n-1})$ to $U_q(A^{(2)}_{2n-1})$.
Furthermore, it is easy to check $(\rho_2\circ \rho_1)\circ \varsigma$ is the identity map on $U_q(A^{(2)}_{2n-1})$ by the formulas of $z^{\pm}(u)$.
\epf
\subsection{ $L$-operators in $U^{ext}_q(A^{(2)}_{2n-1})$}
Proposition \ref{prop:embed} implies the following elements belong to $U^{ext}_q(A^{(2)}_{2n-1})\ot EndV$:
\begin{align}
\label{lplus}
L^{\pm}(u) &:= \tilde{L}^{\pm}(u)\prod_{m=0}^{\infty}
z^{\pm}(u\tss\xi^{-2m-1})\tss z^{\pm}(u\tss\xi^{-2m-2})^{-1}.
\end{align}
Although the entries of $L^{\pm}(u)$ contain a completion of $ZU^{ext}_q(A^{(2)}_{2n-1})$,
it turn outs that the coefficients of the series in $u^{\pm1}$ actually
belong to $U^{ext}_q(A^{(2)}_{2n-1})$ (cf. proof of Thm. \ref{thm:tri-dec}). The following result is
an immediate consequence of Prop. \ref{e:RTT4Ltilde}.
\bpr\label{prop:L-op}
${L}^{\pm}(u)$ satisfy the following relations in $U^{ext}_q(A^{(2)}_{2n-1})\ot EndV^{\ot 2}$ :
\ben
\bal
R(u/v)L^{\pm}_1(u)L^{\pm}_2(v) &= L^{\pm}_2(v)L^{\pm}_1(u)R(u/v), \\[0.4em]
  R(u_{+}/v_{-})L^{+}_1(u)L^{-}_2(v) &= L^{-}_2(v)L^{+}_1(u)R(u_{-}/v_{+}),
\eal
\een
\epr

The triangular decomposition \eqref{rdec} of the universal R-matrix $\mathcal{R}(u)$ and the construction of $L^{\pm}(u)$
imply the following corresponding decomposition of $L^{\pm}(u)$:
\beq
L^{\pm}(u)=F^{\pm}(u)H^{\pm}(u)E^{\pm}(u),
\eeq
where
\beq
\bal
F^{+}(u)&=(\id\otimes \pi_{V})\ts \Rc^{>0}_{21}(u_{+}),\quad
E^{+}(u)=M^{-1}(\id\otimes \pi_{V})\ts \Rc^{<0}_{21}(u_{+}),\\
H^{+}(u)
&=(\id\otimes \pi_{V})\mathcal{R}^{0}_{21}(u_{+})M
{}\times \prod_{m=0}^{\infty} z^{+}(u\xi^{-2m-1})z^{+}(u\xi^{-2m-2})^{-1},
\eal
\eeq
\beq
\bal
F^{-}(u)&=(\id\otimes \pi_{V})\ts \Rc^{>0}(1/u_{+})^{-1}M,\quad
E^{-}(u)=(\id\otimes \pi_{V})\ts \Rc^{<0}(1/u_{+})^{-1},\\
H^{-}(u)
&=M^{-1}(\id\otimes \pi_{V})\mathcal{R}^{0}(1/u_{+})^{-1}
{}\times \prod_{m=0}^{\infty} z^{-}(u\xi^{-2m-1})z^{-}(u\xi^{-2m-2})^{-1},
\eal
\eeq
where $M=\sum\limits_{i\neq n}1\ot e_{ii}+q^{-c}\ot e_{nn}$.

We define the following generating series for the Drinfeld generators $x^{\pm}_{i,k}$ of $U_q(A^{(2)}_{2n-1})$:
\ben
\bal
x_i^{-}(u)^{\geqslant 0}&=\sum\limits _{k\geqslant 0}x^{-}_{i,-k}u^{k},
\qquad\quad x^{+}_{i}(u)^{> 0}=\sum\limits _{k> 0}x^{+}_{i,-k}u^{k},\\
x_i^{-}(u)^{< 0}&=\sum\limits _{k> 0}x^{-}_{i,k}u^{-k},
\qquad\quad x^{+}_{i}(u)^{\leqslant 0}=\sum\limits _{k\geqslant 0}x^{+}_{i,k}u^{-k}.
\eal
\een
For $i=1,\dots, n-1$, set
\ben
\bal
f_i^{+}(u)&=(q_i -q_i^{-1})x_i^{-}(u_+q^{-i})^{\geqslant 0},\qquad\quad
e_i^{+}(u)=(q_i -q_i^{-1})x^{+}_{i}(u_{-}q^{-i})^{> 0},\\
f_{i}^{-}(u)&=(q_i^{-1}-q_i)x^{-}_{i}(u_{-}q^{-i})^{<0},\qquad\quad
e_i^{-}(u)=(q_i^{-1}-q_i)x^{+}_{i}(u_{+}q^{-i})^{\leqslant 0}
\eal
\een
and
\ben
\bal
f_n^{+}(u)&=(q_n -q_n^{-1})u_{+}x_n^{-}(u_+q^{-n})^{\geqslant 0},\qquad\quad
e_n^{+}(u)=(q_{n}-q_{n}^{-1})(u_{-})^{-1}x^{+}_{n}(u_{-}q^{-n})^{>0},\\
f_n^{-}(u)&=(q_n^{-1}-q_n)u_{-}x^{-}_{n}(u_{-}q^{-n})^{<0},\qquad\quad
e_n^{-}(u)=(q_n^{-1}-q_n)(u_{+})^{-1}x_n^{+}(u_{+}q^{-n})^{\leqslant 0}.
\eal
\een

We will need the following result, which is similarly proved as in \cite[Lem. 5.3]{jlm:iC}.
\ble\label{lem:K} The image $(\id\ot \pi_V)\big(T_{21}\big)$ is the diagonal matrix
\begin{equation*}
\diag\Big[\prod\limits_{a=1}^{n-1}k_ak_n^{1/2}, \prod\limits_{a=2}^{n-1}k_ak_n^{1/2}, \dots, k_n^{1/2}, k_n^{-1/2}, k_{n-1}^{-1}k_{n}^{-1/2},
\dots, \prod\limits_{a=1}^{n-1}k_a^{-1}k_n^{-1/2}\Big]
\end{equation*}
\ele

\bth\label{thm:tri-dec}
In $U^{ext}_q(A^{(2)}_{2n-1})$, we have
\ben
F^{\pm}(u)=
\begin{bmatrix}
  1 & &  & &  &  &  &  \\
  f_{1}^{\pm}(u) & \qquad 1 &  &  & & &\bigcirc  &  \\
   &  \qquad \ddots& & \ddots&  &  &  & \\
   &  &  & f^{\pm}_n(u)& 1&  &  &  \\[0.5em]
   &  &  &  & -f^{\pm}_{n-1}(u\xi q^{2(n-1)}) & 1 & &  \\
   &  \qquad \bigstar & &  &  &  \ddots& \ddots &  \\[0.2em]
   &  &  &  &  &  &  -f_{1}^{\pm}(u \xi q^{2})& 1
\end{bmatrix},
\een
\ben
E^{\pm}(u)=
\begin{bmatrix}
  1 & e_1^{\pm}(u) &  &  &  &  &  &  \\
    & 1 & e_2^{\pm}(u) &  &  &  &  \bigstar& \\
   &  & \ddots & \ddots &  &  &  &  \\
   &  &  & 1 & e_n^{\pm}(u) &  &  &  \\
   & &  &  & \ddots & \ddots &  &  \\
   &  & &  &  & 1 &  -e_2^{\pm}(u\xi q^4)&  \\[0.5em]
   & \bigcirc &  &  &  &  &  1& -e_1^{\pm}(u\xi q^2) \\
  &  &  &  &  &  &  & 1
\end{bmatrix},
\een
and
\ben
H^{\pm}(u)=\diag\tss\big[h_1^{\pm}(u),\dots,h^{\pm}_{n}(u),
z^{\pm\tss [1]}(u)\tss h_{n}^{\pm}(u\xi^{[1]})^{-1},\dots,
z^{\pm\tss [n]}(u)\tss h_1^{\pm}(u\xi^{[n]})^{-1}\big].
\een
\eth
\bpf
In order to study the $(i+1,i)-$ entry in $F^{+}(u)$, we only need to
evaluate the image of
the product
\ben
\prod_{k\geqslant 0}
 \exp_{q_i}\big((q_i^{-1}-q_i) F_{\alpha_i+kd_i\delta}\otimes D_{u_{+}}( E_{\alpha_i+kd_i\delta})\big)
\een
for simple roots $\al_i$ ($1\leq i\leq n$). In the following let $t_i=q_i-q_i^{-1}$.
By the isomorphism of Sec.~\ref{subsec:isoDJD}, the above is rewritten as
\ben
\prod_{k\geqslant 0}
 \exp_{q_i}\left(-t_ix^{-}_{i,-kd_{i}}\otimes D_{u_{+}}(x^{+}_{i,kd_{i}})\right).
\een
Suppose $i\leqslant n-1$, then $d_i=1$. By Proposition~\ref{prop:FFRep} it follows that
\begin{align}
\label{itpiv}
&(\id\otimes \pi_V)\prod_{k\geqslant 0}
 \exp_{q_i}\big(-t_i(u_{+})^{kd_i}x^{-}_{i,-kd_i}\otimes x^{+}_{i,kd_{i}}\big)\\
& =\prod_{k\geqslant 0}\exp_{q_i}\left(t_i(u_{+}q^{-i})^{k}x^{-}_{i,-k}\otimes e_{i+1,i}
 -t_i(u_{+}\xi q^{i})^{k}x^{-}_{i,-k}\otimes e_{i',(i+1)'}\right).
 \non
\end{align}
Expanding the $q$-exponent, we can write this expression in the form
\ben
\bal
 1&+t_i\sum\limits_{k\geqslant 0}x^{-}_{i,-k}(u_{+}q^{-i})^{k}\otimes e_{i+1,i}
 -t_i\sum\limits_{k\geqslant 0}x^{-}_{i,-k}(u_{+}\xi q^{i})^{k}\otimes e_{i',(i+1)'}\\
 {}=1&+t_ix^{-}_{i}(u_{+}q^{-i})^{\geqslant 0}\otimes e_{i+1,i}
 -t_ix^{-}_{i}(u_{+}\xi q^{i})^{\geqslant 0}\otimes e_{i',(i+1)'}
\eal
\een
which coincides with $1+f^{+}_{i}(u)\otimes e_{i+1,i}-f^{+}_{i}(u\xi q^{2i})\otimes e_{i',(i+1)'}$,
as required.
Now we consider the case $i=n$, so $d_n=2$ and $D_u(x^{+}_{n,2k})=u^{2k+1}x^{+}_{n,2k}$, it follows from
Prop.~\ref{prop:FFRep} that
\begin{align}
&(\id\otimes \pi_V)\prod_{k\geqslant 0}
 \exp_{q_n}\left(-t_nx^{-}_{n,-kd_{n}}\otimes D_{u_{+}}(x^{+}_{n,kd_{n}})\right)\\
&=\prod_{k\geqslant 0}\exp_{q_n}\left(t_nu_{+}x^{-}_{n,-2k}(u_{+}q^{-n})^{2k}\otimes e_{n+1,n}\right).
 \non
\end{align}
Expanding the $q$-exponent, we can write this expression in the form
\ben
\bal
 1+t_n\sum\limits_{k\geqslant 0}u_{+}x^{-}_{n,-2k}(u_{+}q^{-n})^{2k}\otimes e_{n+1,n}
 =1+t_nu_{+}x^{-}_{n}(u_{+}q^{-n})^{\geqslant 0}\otimes e_{n+1,n}
\eal
\een
which coincides with $1+f^{+}_{n}(u)\otimes e_{n+1,n}$
as required.

For the $(i,i+1)-$th entry in $E^{+}(u)$, we
first evaluate the image of the product
\begin{align*}
&T_{21}^{-1}\prod_{k> 0}
 \exp_{q_i}\left(-t_iF_{-\alpha_i+kd_i\delta}\otimes D_{u_{+}}(E_{-\alpha_i+kd_{i}\delta})\right) T_{21}\\
=&T_{21}^{-1}\prod_{k> 0}
 \exp_{q_i}\left(-t_iq^{-kd_{i}c}x^{+}_{i,-kd_{i}}k_i\otimes D_{u_{+}}(q^{kd_{i}c}k_i^{-1}x^{-}_{i,kd_{i}})\right)T_{21}.
\end{align*}
For $i=1,...,n-1$, $d_i=1$
 and $D_{u}(x^{-}_{i,k})=u^{k}x^{-}_{i,k}$, by the action of $\pi_V$ (Prop. \ref{prop:FFRep}) we have
\ben
\bal
&(id\otimes \pi_V)\prod_{k> 0}
 \exp_{q_i}\left(-t_iq^{-kd_{i}c}x^{+}_{i,-kd_{i}}k_i\otimes D_{u_{+}}(q^{kd_{i}c}k_i^{-1}x^{-}_{i,kd_{i}})\right)\\
& =
 \prod_{k> 0}exp_{q_i}\left(q_it_i u_{-}^{k}x^{+}_{i,-k}k_i\otimes(q^{-ik}e_{i,i+1}
 -\xi^{k}q^{ik} e_{(i+1)',i'})\right)\\
&=1+q_it_i\sum\limits_{k> 0}x^{+}_{i,-k}(u_{-}q^{-i})^{k} k_i\otimes e_{i,i+1}
 -q_it_i\sum\limits_{k\geq 0}x^{+}_{i,-k}(u_{-}\xi q^{i})^{k}k_i\otimes e_{(i+1)',i'}.
\eal
\een
Then by Lemma \ref{lem:K}, we have for $i<n-1$,
\ben
\bal
(id\otimes \pi_V)&
\big(T_{21}^{-1}
\exp_{q_i}\left(-t_iq^{-kd_{i}c}x^{+}_{i,-kd_{i}}k_i\otimes D_{u_{+}}(q^{kd_{i}c}k_i^{-1}x^{-}_{i,kd_{i}})\right)
T_{21}\big)\\
 &=1+q_it_i
\sum\limits_{k> 0}\left(\prod_{j=i}^{n-1}k_j^{-1} k_n^{-1/2}x^{+}_{i,-k}\prod_{j=i}^{n-1}k_j k_n^{1/2}\right)
(u_{-}q^{-i})^{k} \otimes e_{i,i+1}\\
 &-q_it_i
\sum\limits_{k\geq 0}\left(\prod_{j=i+1}^{n-1}k_j k_n^{1/2}x^{+}_{i,-k}\prod_{j=i+1}^{n-1}k_j^{-1} k_n^{-1/2}
\right)(u_{-}\xi q^{i})^{k}k_i\otimes e_{(i+1)',i'}.\\
\eal
\een
By using the fact $k_ix_{j,k}^{\pm}k_i^{-1}=q_i^{\pm A_{ij}}x_{jk}^{\pm}$, the last equation is equal to
\ben
\bal
&=1+t_ix^{+}_{i}(u_{-}q^{-i})^{> 0}\otimes e_{i,i+1}
 -t_ix^{+}_{i}(u_{-} \xi q^{i})^{> 0}\otimes e_{(i+1)',i'}\\
&=1+e^{+}_{i}(u)\otimes e_{i,i+1}
 -e^{+}_{i}(u\xi q^{2i})\otimes e_{(i+1)',i'}.
 \eal
\een

Thus, we get the $(i,i+1)$-th entry of $E^{+}(z)$ is $e^{+}_{i}(z)$,
and the corresponding $((i+1)',i')$-th entry is $-e^{+}_{i}(z\xi q^{2i})$.

For $i= n-1$, $q_{n-1}=q$, so we have
\ben
\bal
(id&\ot \pi_V)\bigg(T_{21}^{-1}
\exp_{q}\big(-tq^{-kc}x^{+}_{n-1,-k}k_{n-1}\otimes D_{u_{+}}(q^{kc}k_{n-1}^{-1}x^{-}_{n-1,k})\big)
T_{21}\bigg)\\
 &=1+qt
\sum\limits_{k> 0}\left(k_{n-1}^{-1} k_n^{-1/2}x^{+}_{n-1,-k}k_{n-1} k_n^{1/2}\right)
(u_{-}q^{-(n-1)})^{k} \otimes e_{n-1,n}\\
 &-qt
\sum\limits_{k\geq 0}\left(k_n^{1/2}x^{+}_{n-1,-k}k_n^{-1/2}
\right)(u_{-}\xi q^{n-1})^{k}\otimes e_{n',(n-1)'}.
\eal
\een
By using the fact $k_ix_{j,k}^{\pm}k_i^{-1}=q_i^{\pm A_{ij}}x_{jk}^{\pm}$,
 we have the last equation is simplified to
\ben
\bal
&=1+tx^{+}_{n-1}(u_{-}q^{-(n-1)})^{> 0}\otimes e_{n-1,n}
 -tx^{+}_{n-1}(u_{-} \xi q^{n-1})^{> 0}\otimes e_{n',(n-1)'}\\
&=1+e^{+}_{n-1}(z)\otimes e_{n-1,n}
 -e^{+}_{n-1}(z\xi q^{2(n-1)})\otimes e_{n',(n-1)'}.
 \eal
\een
Thus the $(n-1,n)$-th and $(n',(n-1)')$-th entries of $E^{+}(z)$ are $e^{+}_{n-1}(z)$ and $-e^{+}_{n-1}(z\xi q^{2(n-1)})$.

For $i=n$, by the action of $\pi_V$ (see Proposition \ref{prop:FFRep}) we have
\ben
\bal
&(id\otimes \pi_V)
\prod_{k> 0}
 \exp_{q_n}\left(-t_nq^{-kd_{n}c}x^{+}_{n,-kd_{n}}k_n\otimes D_{u_{+}}(q^{kd_{n}c}k_n^{-1}x^{-}_{n,kd_{n}})\right)\\
&=(id\otimes \pi_V)\prod_{k> 0}
 exp_{q_n}\left(-t_n(u_{-})^{2k}(u_{+})^{-1}x^{+}_{n,-2k}k_n\otimes q^{2kc}k_n^{-1}x^{-}_{n,2k}\right)\\
& =
 \prod_{k> 0}exp_{q_n}\left(q_nt_n(u_{-}q^{-n})^{2k}(u_{+})^{-1}x^{+}_{n,-2k}k_n\otimes e_{n,n+1}\right)\\
&=1+q_nt_n\sum\limits_{k> 0}(u_{-}q^{-n})^{2k}(u_{+})^{-1}x^{+}_{n,-2k} k_n\otimes e_{n,n+1}\\
\eal
\een
Then by Lemma \ref{lem:K} and also $k_ix_{j,k}^{\pm}k_i^{-1}=q_i^{\pm A_{ij}}x_{jk}^{\pm}$, we have that
\ben
\bal
&(id\otimes \pi_V)
\big(T_{21}^{-1}
\prod_{k> 0}\exp_{q_n}\left(-t_nq^{-kd_{n}c}x^{+}_{n,-kd_{n}}k_n\otimes D_{u_{+}}(q^{kd_{n}c}k_n^{-1}x^{-}_{n,kd_{n}})\right)
T_{21}\big)\\
 &=1+q_{n}t_n
\sum\limits_{k> 0}\left(k_n^{-1/2}x^{+}_{n,-2k}k_{n}^{1/2}\right)(u_{-}q^{-n})^{2k}(u_{+})^{-1}
 \otimes e_{n,n+1}\\
&=1+t_n(u_{+})^{-1}
x^{+}_{n}(u_{-}q^{-n})^{>0}
 \otimes e_{n,n+1}=1+q^{-c}e^{+}_{n}(z)\otimes e_{n,n+1}.
 \eal
\een
Thus, we get the $(n,n+1)$-th entry of $E^{+}(z)$ is $e^{+}_{n}(z)$.

Now we consider $H^{+}(u)$, in which we only treat the $(n,n)$- and $(n',n')-$ entries.  

By definition and the action of $\pi_V$, 
the $(n,n)$-entry of $(id\ot \pi_V)\mathcal{R}^{0}_{21}(u_{+})$ is
\begin{align*} 
&\exp\left[\sum\limits_{k>0}\sum\limits_{j=1}^{n-1}-t_j
q^{-nk}z_{n-1,j}^{k}a_{j,-k}u^k \right]
\exp\left[\sum\limits_{k>0}\sum\limits _{j=1}^{n-1}t_j
(-1)^{k}q^{-(n-1)2k}z_{nj}^{2k} a_{j,-2k}u^{2k}\right]\\
&\cdot\exp\left[\sum\limits_{k>0}-t_n\left(q^{-2nk}(-1)^k
z_{n-1,n}^{2k}-q^{-(n-1)2k}z_{n,n}^{2k}\right) a_{n,-2k}u^{2k}\right]
\end{align*}
Due to \eqref{Eq:Zij} for $z^k_{ij}$, the above can be rewritten as (note $\xi=-q^{-2n}$)   
\beq
\exp\left[\sum\limits_{k>0}\sum\limits_{j=1}^{n-1}t_j
\frac{q^{jk}-q^{-jk}}{1+\xi^{-k}}a_{j,-k}u^k\right]
\exp\left[\sum\limits_{k>0}t_n
\frac{(-q^{n})^{2k}}{1+\xi^{-2k}}a_{n,-2k}u^{2k}\right]
\eeq
Expanding the fractions into power series, we can write the expression as
\begin{multline}
\non
\exp\Big(\sum_{k>0}\sum_{j=1}^{n-1}\sum_{m=0}^{\infty}(-1)^mt_j
\xi^{-mk}(q^{jk}-q^{-jk})u^{k}a_{j,-k}\Big)\\
{}\times \exp\Big(\sum_{k>0}
\sum_{m=0}^{\infty}(-1)^mt_n\xi^{-m2k}(-q^{n})^{2k}u^{2k}a_{n,-2k}\Big),
\end{multline}
which becomes the following by using $\overline{\varphi}^-_{i}(u)$ (see \eqref{psiiu}):
\begin{multline}
\non
\prod_{m=0}^{\infty} \prod_{j=1}^{n-1}\overline{\varphi}^-_j(u\xi^{-2m}q^{j})^{-1}
\overline{\varphi}^-_j(u\xi^{-2m}q^{-j})
\overline{\varphi}^-_j(u\xi^{-2m-1}q^{j})\overline{\varphi}^-_j(\xi^{-2m-1}q^{-j})^{-1}\\
{}\times \prod_{m=0}^{\infty} \overline{\varphi}^-_{n}(-u\xi^{-2m}q^{n})^{-1}
\overline{\varphi}^-_{n}(-u\xi^{-2m+1}q^{n})\times \overline{\varphi}^-_{n}(-u\xi q^{n})^{-1}.
\end{multline}
Setting $\overline{h}^{+}_{i}(u)=t_i^{-1}\tss {h}^{+}_{i}(u)$ with
$t_i=h_{i,0}^{+}$, it follows from Proposition~\ref{prop:embed} that the above is
\begin{multline}
\non
\prod_{m=0}^{\infty} \prod_{j=1}^{n-1}\overline{h}^{+}_{j}(u\xi^{-2m}q^{2j})
\overline{h}^{+}_{j}(u\xi^{-2m-1}q^{2j})^{-1}
\times\prod_{m=0}^{\infty} \prod_{j=1}^{n}\overline{h}^{+}_{j}(u\xi^{-2m}q^{2j-2})^{-1}
\overline{h}^{+}_{j}(u\xi^{-2m-1}q^{2j-2})\\
\times\prod_{m=0}^{\infty}\overline{h}^{+}_{n+1}(u\xi^{-2m-1})^{-1}
\overline{h}^{+}_{n+1}(u\xi^{-2m-2})\times \overline{h}^{+}_n(u)q^c.
\end{multline}
One then invokes \eqref{zpm} for $z^{\pm}(u)$ to conclude that
\ben
\exp\Big(\sum_{k>0}\sum_{j=1}^{n}t_j\tilde{B}_{1,j}(q^{k})u^{k}a_{j,-k}\Big)
=\prod_{m=0}^{\infty}z^{+}(u\xi^{-2m-1})^{-1}z^{+}(u\xi^{-2m-2})\times \overline{h}^{+}_n(u)q^c.
\een
In particular, the coefficients of $u^n$ in the infinite product (occurring in \eqref{lplus}) belong to the algebra $U^{\ext}_q(\wh{\spa}_{2n})$.

Furthermore, Lemma~\ref{lem:K} implies that
the $(n,n)$-entry of the matrix
$(\id\ot \pi_V)(T_{21})$ equals
$k_n^{1/2}=t_n$, so the $(n,n)$-entry of the matrix $H^{+}(u)$ is $h_n^{+}(u)$.

By definition 
and recalling the action of $\pi_V$,
we have $(n',n')-$ entry of $(id\ot \pi_V)\mathcal{R}^{0}_{21}(u_{+})$ is
\beql{eq:hn'}
\bal
&\exp\left[\sum\limits_{k>0}\sum\limits_{j=1}^{n-1}t_j
\xi^kq^{nk}z_{n-1,j}^{k}a_{j,-k}u^k\right]
\times\exp\left[\sum\limits_{k>0}\sum\limits _{j=1}^{n-1}t_j
(-1)^{k}q^{-(n+1)2k}z_{nj}^{2k} a_{j,-2k}u^{2k} \right]\\
&\times\exp\left[\sum\limits_{k>0}-t_n
\xi^{2k}\left(
-(-1)^k q^{-2nk}z_{n-1,n}^{2k}+q^{-(n+1)2k}z_{nn}^{2k}\right)a_{n,-2k}u^{2k}\right]
\eal
\eeq

Recalling the meaning of $z^k_{ij}$, \eqref{eq:hn'} is rewritten as
\beq
\exp\left[\sum\limits_{k>0}\sum\limits_{j=1}^{n-1}t_j
\frac{q^{jk}-q^{-jk}}{1+\xi^{-k}}a_{j,-k}u^k\right]
\exp\left[\sum\limits_{k>0}t_n
\frac{(-q^{-n})^{2k}}{1+\xi^{-2k}}a_{n,-2k}u^{2k}\right]
\eeq
Expanding the fractions to power series, the above can be written as
\begin{multline}
\non
\exp\Big(\sum_{k>0}\sum_{j=1}^{n-1}\sum_{m=0}^{\infty}(-1)^mt_j
\xi^{-mk}(q^{jk}-q^{-jk})u^{k}a_{j,-k}\Big)\\
{}\times \exp\Big(\sum_{k>0}
\sum_{m=0}^{\infty}(-1)^mt_n\xi^{-m2k}(q^{-n})^{2k}u^{2k}a_{n,-2k}\Big),
\end{multline}
which is converted into the following by using the series $\overline{\varphi}^-_{i}(u)$:
\begin{multline}
\non
\prod_{m=0}^{\infty} \prod_{j=1}^{n-1}\overline{\varphi}^-_j(u\xi^{-2m}q^{j})^{-1}
\overline{\varphi}^-_j(u\xi^{-2m}q^{-j})
\overline{\varphi}^-_j(u\xi^{-2m-1}q^{j})\overline{\varphi}^-_j(\xi^{-2m-1}q^{-j})^{-1}\\
{}\times \prod_{m=0}^{\infty} \overline{\varphi}^-_{n}(u\xi^{-2m}q^{-n})
\overline{\varphi}^-_{n}(u\xi^{-2m-1}q^{-n})^{-1},
\end{multline}
and this is exactly (by Prop.~\ref{prop:embed}):
\begin{multline}
\non
\prod_{m=0}^{\infty} \prod_{j=1}^{n-1}\overline{h}^{+}_{j}(u\xi^{-2m}q^{2j})
\overline{h}^{+}_{j}(u\xi^{-2m-1}q^{2j})^{-1}
\times\prod_{m=0}^{\infty} \prod_{j=1}^{n}\overline{h}^{+}_{j}(u\xi^{-2m}q^{2j-2})^{-1}
\overline{h}^{+}_{j}(u\xi^{-2m-1}q^{2j-2})\\
\times\prod_{m=0}^{\infty}\overline{h}^{+}_{n+1}(u\xi^{-2m-1})^{-1}
\overline{h}^{+}_{n+1}(u\xi^{-2m-2})\times \overline{h}^{+}_{n+1}(u).
\end{multline}
Finally using the formula \eqref{zpm} of $z^{\pm}(u)$, we get that 
\ben
\exp\Big(\sum_{k>0}\sum_{j=1}^{n}t_j\tilde{B}_{1,j}(q^{k})u^{k}a_{j,-k}\Big)
=\prod_{m=0}^{\infty}z^{+}(u\xi^{-2m-1})^{-1}z^{+}(u\xi^{-2m-2})\times \overline{h}^{+}_{n+1}(u).
\een
As a consequence, the coefficients of $u^n$ in the infinite product (also \eqref{lplus}) is an element in the algebra $U^{\ext}_q(\wh{\spa}_{2n})$.

Moreover, Lemma~\ref{lem:K} implies that
the $(n,n)$-entry of the matrix
$(\id\ot \pi_V)(T_{21})$ is
$k_n^{-1/2}=t_n^{-1}$, hence
the $(n,n)$-entry of the matrix $H^{+}(u)$ is $$h_{n+1}^{+}(u)=z^{\pm\tss [1]}(u)\tss h_{n}^{\pm}(u\xi^{[1]})^{-1}.$$
\epf
\section{R-matrix realization of $U_q(A^{(2)}_{2n-1})$}
Recall that two related R-matrices have been defined in \eqref{Rm}. Accordingly two R-matrix algebras $U(R)$ and $U(\overline{R})$ are defined in this section.

\subsection{The algebras $U(R)$ and $U(\overline{R})$ }

\bde
The associative algebra $U(R)$ over $\CC(q)$
is generated by an invertible central element $q^{c/2}$ and
elements  ${l}^{\pm}_{ij}[\mp m]$, $1\leqslant i,j\leqslant 2n, m\in \ZZ_{+}$,
subject to the following defining relations:
\ben
{l}_{ij}^{+}[0]={l}^{-}_{ji}[0]=0\quad\text{for}\quad i<j,~~ {l}^{+}_{nn}[0]\ts {l}_{n+1,n+1}^{+}[0]=1
\Fand {l}^{+}_{ii}[0]\ts {l}_{ii}^{-}[0]={l}_{ii}^{-}[0]\ts{l}^{+}_{ii}[0]=1,
\een
and the remaining relations are written in matrix forms as follows. Let
\begin{equation}\label{liju}
{L}^{\pm}(u)=\sum\limits_{i,j=1}^{2n}{l}_{ij}^{\pm}(u)\ot e_{ij}\in
U(R)[[u,u^{-1}]]\ot\End\CC^{2n}.
\end{equation}
where the entries ${l}^{\pm}_{ij}(u)=\sum_{m=0}^{\infty}{l}^{\pm}_{ij}[\mp m]\ts u^{\pm m}$ are $U(R)$-valued formal power series.
Consider the tensor product algebra $U(R)\ot\End\CC^{2n}\ot\End\CC^{2n}$
and set 
\beql{lonetwo}
{L}^{\pm}_1(u)=\sum\limits_{i,j=1}^{2n}e_{ij}\ot 1\ot {l}_{ij}^{\pm}(u)
\Fand
{L}^{\pm}_2(u)=\sum\limits_{i,j=1}^{2n}1\ot e_{ij}\ot {l}_{ij}^{\pm}(u).
\eeq
Then the remaining defining relations take the form
\begin{align}
\label{rllss}
{R}(u/v)L^{\pm}_{1}(u)L^{\pm}_2(v)&=L^{\pm}_2(v)L^{\pm}_{1}(u){R}(u/v),\\[0.4em]
{R}(uq^c/v)L^{+}_1(u)L^{-}_2(v)&=L^{-}_2(v)L^{+}_1(u){R}(uq^{-c}/v),
\label{rllmp}
\end{align}
\ede

Similarly,
the {\it algebra $U(\overline{R})$ over $\CC(q)$}
is generated by an invertible central element $q^{c/2}$ and
elements  ${\ell}^{\tss\pm}_{ij}[\mp m]$ with $1\leqslant i,j\leqslant 2n$ and $m\in \ZZ_{+}$
such that
\ben
{\ell}_{ij}^{+}[0]={\ell}^{-}_{ji}[0]=0\quad\text{for}\quad i<j,\quad  {\ell}^{+}_{nn}[0]\ts {\ell}_{n+1,n+1}^{+}[0]=1
\Fand {\ell}^{+}_{ii}[0]\ts {\ell}_{ii}^{-}[0]={\ell}_{ii}^{-}[0]\ts{\ell}^{+}_{ii}[0]=1.
\een
The remaining defining relations of $U(\overline{R})$ take the form
\begin{align}
\label{gen rel1}
\overline{R}(u/v)\ts\Lc^{\pm}_{1}(u)\ts\Lc^{\pm}_2(v)&
=\Lc^{\pm}_2(v)\ts\Lc^{\pm}_{1}(u)\ts\overline{R}(u/v),\\[0.4em]
\overline{R}(u\tss q^c/v)\ts\Lc^{+}_1(u)\ts\Lc^{-}_2(v)&
=\Lc^{-}_2(v)\ts\Lc^{+}_1(u)\ts\overline{R}(u\tss q^{-c}/v),
\label{gen rel2}
\end{align}
where $\Lc^{\pm}(u)$ are defined similarly as in \eqref{liju}-\eqref{lonetwo}.

\bre
The R-matrix $\overline{R}(u)$ obeys the unitarity property: $\overline{R}_{12}(u)\overline{R}_{21}(1/u)=1$,
which implies that
relation \eqref{gen rel2} can be written in the equivalent form
\beql{gen rel3}
\overline{R}(u\tss q^{-c}/v)\ts\Lc^{-}_1(u)\ts\Lc^{+}_2(v)
=\Lc^{+}_2(v)\ts\Lc^{-}_1(u)\ts\overline{R}(u\tss q^{c}/v).
\eeq
\ere
\bre\label{rem:gln}
The defining relations satisfied by the series
$\ell^{\tss\pm}_{ij}(u)$
with $1\leqslant i,j\leqslant n$ coincide with those for the quantum affine algebra
$U_q(\wh{\gl}_n)$ in \cite{df:it}.
\ere

The algebras
$U(R)$ and $U(\overline{R})$ are related via a Heisenberg algebra as follows.
Introduce the Heisenberg algebra $\Hc_q(n)$
with generators $q^c$ and $\be_r$, $r\in\ZZ\setminus \{0\}$ subject to
the defining relations 
\ben
\big[\be_r,\be_s\big]=\de_{r, -s}\ts \al_r,\qquad r\geqslant 1,
\een
and $q^c$ is central and invertible, where the scalars $\al_r$ are defined by the expansion
\ben
\exp\ts\sum_{r=1}^{\infty}\al_r u^r=\frac{g(u\tss q^{-c})}{g(u\tss q^{c})}.
\een
So we have the identity
\ben
g(u\ts q^{c}/v)\ts \exp\ts\sum_{r=1}^{\infty}\be_r u^r\cdot
\ts \exp\ts\sum_{s=1}^{\infty}\be_{-s} v^{-s}
=g(u\ts q^{-c}/v)\ts \exp\ts\sum_{s=1}^{\infty}\be_{-s} v^{-s}
\cdot \exp\ts\sum_{r=1}^{\infty}\be_r u^r.
\een

\bpr\label{prop:homheis}
The mappings
\beql{heisext}
L^+(u)\mapsto \exp\ts\sum_{r=1}^{\infty}\be_{-r} u^{-r}\cdot
\Lc^{\ts +}(u),\qquad
L^-(u)\mapsto \exp\ts\sum_{r=1}^{\infty}\be_{r} u^{r}\cdot
\Lc^{\ts -}(u),
\eeq
define a homomorphism $U(\overline{R})\to \Hc_q(n)\ot_{\CC[q^c,\ts q^{-c}]}U(R)$.
\qed
\epr
Let $t$ be the matrix transposition of $\End\CC^{2n}$ given by $e_{ij}^{t}=e_{j',i'}$, and can be lifted to
a partial transpose on tensor products of $\End\CC^{2n}$: the partial transpose $\tra_a$
refers to applying $t$ to the $a$-th copy in the tensor product.
The following {\it crossing symmetry} relations are
satisfied by the $R$-matrices:
\begin{align}\label{RDRD}
 \overline{R}(u)D_1\overline{R}(u\tss\xi)^{\tra_1}D_1^{-1}
 &=\frac{(u-q^2)(u\tss\xi-1)}{(1-u)(1-u\tss\xi q^2)},\\[0.4em]
 {R}(u)D_1{R}(u\tss\xi)^{\tra_1}D_1^{-1}&=\xi^2 q^{-2},
 \label{crsymr}
\end{align}
where the diagonal matrix $D=\diag\{q^{\overline{1}},\dots ,q^{\overline{2n}}\}$
and the meaning of the subscripts
is the same as in \eqref{lonetwo}.

\bpr\label{prop:central}
In the algebras $U(R)$ and $U(\overline{R})$ we have
the relations
\beql{DLDL}
D{L}^{\pm}(u\tss\xi)^{\tra}D^{-1}{L}^{\pm}(u)
={L}^{\pm}(u)D{L}^{\pm}(u\tss\xi)^{\tra}D^{-1}={z}^{\pm}(u)\ts1,
\eeq
and
\beql{DLbarDLbar}
D\Lc^{\pm}(u\tss\xi)^{\tra}D^{-1}\Lc^{\pm}(u)
=\Lc^{\pm}(u)D\Lc^{\pm}(u\tss\xi)^{\tra}D^{-1}=\z^{\pm}(u)\ts1,
\eeq
for certain series ${z}^{\pm}(u)$ and $\z^{\pm}(u)$ with coefficients
in the respective algebra.
\epr

\bpf The proof is similar for both cases, so we consider the algebra $U(\overline{R})$.
Multiply both sides of \eqref{gen rel1} by $u/v-\xi$ and set
$u/v=\xi$ to get
\beql{QL1L2}
Q\tss\Lc^{\pm}_1(u\tss\xi)\Lc^{\pm}_2(u)=\Lc^{\pm}_2(u)\Lc^{\pm}_1(u\tss\xi)\tss Q,
\eeq
where $Q=\sum\limits_{i,j=1}^{2n}q^{\overline{i}-\overline{j}}e_{i',j'}\ot e_{ij}=D_1^{-1}P^{\ts \tra_1}D_1$
and $P=\sum\limits_{i,j=1}^{2n}e_{ij}\ot e_{ji}$.
Therefore, \eqref{QL1L2} takes the form
\beql{pdld}
P^{\ts\tra_1}D_1\Lc_1^{\pm}(u\tss\xi)D_1^{-1}\Lc_2^{\pm}(u)
=\Lc_2^{\pm}(u)D_1\Lc_1^{\pm}(u\tss\xi)D_1^{-1}P^{\ts\tra_1}.
\eeq
The image of the operator $P^{\ts\tra_1}$ in $\End(\CC^{N})^{\ot 2}$
is one-dimensional,
so that each side of this equality
must be equal to $P^{\ts\tra_1}$ times a certain series $\z^{\pm}(u)$
with coefficients in $U(\overline{R})$. Observe that
$P^{\ts\tra_1}D_1=P^{\ts\tra_1}D^{-1}_2$
and $P^{\ts\tra_1}\Lc_1^{\pm}(u\tss\xi)=P^{\ts\tra_1}\Lc_2^{\pm}(u\tss\xi)^{\tra}$
and so we get
\ben
P^{\ts\tra_1}D_2\Lc_2^{\pm}(u\tss\xi)^{\tra}D_2^{-1}\Lc_2^{\pm}(u)
=\Lc_2^{\pm}(u)D_2\Lc_2^{\pm}(u\tss\xi)^{\tra}D_2^{-1}P^{\ts\tra_1}=\z^{\pm}(u)\tss P^{\ts\tra_1}.
\een
The required relations now follow by taking trace of the first copy of $\End\CC^{N}$.
\epf

\bpr\label{prop:centr}
All coefficients of the series $z^{+}(u)$ and $z^{-}(u)$
belong to the center of the algebra $U(R)$.
\epr

\bpf
We will verify that $z^{+}(u)$ commutes with all series $l_{ij}^-(v)$; the remaining
cases follow by similar or simpler arguments.
By the defining relations \eqref{rllmp} we can write
\ben
\bal
D_1L^{+}_1(u\tss\xi)^{\tra}D_1^{-1}L^{+}_1(u)L^{-}_2(v)
=D_1L^{+}_1(u\tss\xi)^{\tra}D_1^{-1}{R}(uq^c/v)^{-1}L^{-}_2(v)L^{+}_1(u){R}(uq^{-c}/v).
\eal
\een
By \eqref{crsymr} the right hand side equals
\ben
\xi^{-2}q^2D_1L^{+}_1(u\tss\xi)^{\tra}{R}(u\tss\xi q^c/v)^{\tra_1}L^{-}_2(v)
D_1^{-1}L^{+}_1(u){R}(uq^{-c}/v).
\een
Applying the patrial transposition $\tra_1$ to both sides in
\eqref{rllmp} we get the relation
\ben
L^{+}_1(u\tss\xi)^{\tra}{R}(u\tss\xi q^c/v)^{\tra_1}L^{-}_2(v)
=L^{-}_2(v){R}(u\tss\xi q^{-c}/v)^{\tra_1}L^{+}_1(u\tss\xi)^{\tra}.
\een
Hence, using \eqref{crsymr} and \eqref{DLDL} we obtain
\ben
\bal
z^{+}(u)L^{-}_2(v)&=D_1L^{+}_1(u\tss\xi)^{\tra}D_1^{-1}L^{+}_1(u)L^{-}_2(v)\\[0.4em]
&=\xi^{-2}q^2 L^{-}_2(v)D_1 {R}(u\tss\xi q^{-c}/v)^{\tra_1}
D_1^{-1}D_1L^{+}_1(u\tss\xi)^{\tra}D_1^{-1}L^{+}_1(u){R}(uq^{-c}/v)\\[0.4em]
&=\xi^{-2}q^2 L^{-}_2(v)D_1 {R}(u\tss\xi q^{-c}/v)^{\tra_1}
D_1^{-1}z^{+}(u){R}(uq^{-c}/v)=L^{-}_2(v)z^{+}(u),
\eal
\een
as required.
\epf

\bre\label{rem:noncent}
The crossing symmetry \eqref{crsymr} of the $R$-matrix $R(u)$ is essential
for Proposition~\ref{prop:centr}.
Although the coefficients of the series $\z^{+}(u)$ and $\z^{-}(u)$ are central
in the respective subalgebras of $U(\overline{R})$ generated by the coefficients of the
series $\ell^+_{ij}(u)$ and $\ell^-_{ij}(u)$, they are not central in the entire
algebra $U(\overline{R})$. We will give the explicit formulas of $\z^{\pm}(u)$ and
${z}^{\pm}(u)$ in Section \ref{Sec:guass-decomposition}.
\qed
\ere

\subsection{Quasi-determinant and quantum minors}
Let $A=[a_{ij}]$ be an $N\times N$ matrix over a ring with $1$.
Denote by $A^{ij}$ the matrix obtained from $A$
by deleting the $i$-th row
and $j$-th column. Suppose that the matrix
$A^{ij}$ is invertible.
The $ij$-th {\it quasi-determinant} of $A$
is defined by the formula
\ben
|A|_{ij}=a_{ij}-r^{\tss j}_i(A^{ij})^{-1}\ts c^{\tss i}_j,
\een
where $r^{\tss j}_i$ is the row matrix obtained from the $i$-th
row of $A$ by deleting the element $a_{ij}$, and $c^{\tss i}_j$
is the column matrix obtained from the $j$-th
column of $A$ by deleting the element $a_{ij}$; see
\cite{gr:dm}. In particular,
the four quasi-determinants of a $2\times 2$ matrix $A$ are
\ben
\bal
&|A|_{11}=a^{}_{11}-a^{}_{12}\ts a_{22}^{-1}\ts a^{}_{21},\qquad
|A|_{12}=a^{}_{12}-a^{}_{11}\ts a_{21}^{-1}\ts a^{}_{22},\\
&|A|_{21}=a^{}_{21}-a^{}_{22}\ts a_{12}^{-1}\ts a^{}_{11},\qquad
|A|_{22}=a^{}_{22}-a^{}_{21}\ts a_{11}^{-1}\ts a^{}_{12}.
\eal
\een
The quasi-determinant $|A|_{ij}$ is also denoted
by boxing the entry $a_{ij}$,
\ben
|A|_{ij}=\left|\begin{matrix}a_{11}&\dots&a_{1j}&\dots&a_{1N}\\
                                   &\dots&      &\dots&      \\
                             a_{i1}&\dots&\boxed{a_{ij}}&\dots&a_{iN}\\
                                   &\dots&      &\dots&      \\
                             a_{N1}&\dots&a_{Nj}&\dots&a_{NN}
                \end{matrix}\right|.
\een

Note that 
elements of the
tensor product algebra $U(\overline{R})\ot \End(\CC^{N})^{\ot m}$ can be viewed
as operators on the space $(\CC^{N})^{\ot m}$ with coefficients in $U(\overline{R})$.
Accordingly, for such an element
\ben
X=\sum_{a_i,b_i}\ts X^{a_1\dots a_m}_{\ts b_1\dots b_m}
\ot e_{a_1b_1}\ot\dots\ot e_{a_mb_m}
\een
where $e_{ij}=|i\rangle\langle j|$, so the minor 
\beql{matrelem}
X^{a_1\dots a_m}_{\ts b_1\dots b_m}=\langle a_1,\dots, a_m\ts |\ts X\ts |\ts b_1,\dots, b_m\ts\rangle
\eeq
and the tensors $X\ts |\ts b_1,\dots, b_m\ts\rangle$ and $\langle a_1,\dots, a_m\ts |\ts X$ are defined as usual.

Consider the algebra $U(\overline{R})$ and for any $2\leqslant i,j\leqslant 2\pr$
introduce the quasideterminant
\ben
s^{\pm}_{ij}(u)=\left|\begin{matrix}
\ell^{\tss\pm}_{11}(u)&\ell^{\tss\pm}_{1j}(u)\\[0.2em]
\ell^{\tss\pm}_{i1}(u)&\boxed{\ell^{\tss\pm}_{ij}(u)}
\end{matrix}\right|
=\ell^{\tss\pm}_{ij}(u)-\ell^{\tss\pm}_{i1}(u)\tss
\ell^{\tss\pm}_{11}(u)^{-1}\tss \ell^{\tss\pm}_{1j}(u).
\een

Let the power series
$\ell^{\tss\pm\tss a_1 a_2}_{~~\tss b_1 b_2}(u)$ ({\it quantum minors})
with coefficients in ${U}(\overline{R})$ be defined by
\beql{quamintau}
\ell^{\tss\pm\tss a_1 a_2}_{~~\tss b_1 b_2}(u)=
\langle a_1,a_2\ts |\ts \wh{R}(q^{-2})\ts \Lc^{\pm}_1(u)\ts \Lc^{\pm}_2(u\tss q^{2})\ts |
\ts b_1,b_2\rangle,
\eeq
where $a_i,b_i\in\{1,\dots,N\}$ and we set
\beql{rhat}
\wh{R}(u)=\frac{uq-q^{-1}}{u-1}\ts\overline{R}(u).
\eeq
The following symmetry properties are straightforward to verify.

\ble\label{lem:skewsymm}\quad
(i)\quad
If $a_1\ne a'_2$ and $a_1<a_2$ then $\ell^{\tss\pm\tss a_1 a_2}_{~~\tss b_1 b_2}(u)
=-q^{-1}\ts\ell^{\tss\pm\tss a_2 a_1}_{~~\tss b_1 b_2}(u)$.
\medskip

\noindent
(ii)\quad
If $b_1\ne b'_2$ and $b_1<b_2$ then $\ell^{\tss\pm\tss a_1 a_2}_{~~\tss b_1 b_2}(u)
=-q\ts \ell^{\tss\pm\tss a_1 a_2}_{~~\tss b_2 b_1}(u)$.
\qed
\ele

\ble\label{lem:toneone}
For any $2\leqslant i,j\leqslant 2\pr$ we have
\beql{quasi-minor}
s^{\pm}_{ij}(u)=\ell^{\tss\pm}_{11}(uq^{-2})^{-1}\tss \ell^{\tss\pm\tss 1 i}_{~~\tss 1 j}(uq^{-2}).
\eeq
Moreover,
\beql{commtoo}
\big[\ell^{\tss\pm}_{11}(u), \ell^{\tss\pm\tss 1 i}_{~~\tss 1 j}(v)\big]=0
\eeq
and
\beql{pm}
\frac{q^{-1}u_{\pm}-q\tss v_{\mp}}{u_{\pm}-v_{\mp}}\ts
\ell_{11}^{\tss\pm}(u)\ts\ell^{\tss\mp\tss 1i}_{~~1j}(v)
=\frac{q^{-1}u_{\mp}-q\tss v_{\pm}}{u_{\mp}-v_{\pm}}\ts
\ell^{\tss\mp\tss 1i}_{~~1j}(v)\ts\ell_{11}^{\tss\pm}(u).
\eeq
\ele

The following consequences of Lemma \ref{lem:toneone} are easily seen:
for $2\leqslant i,j\leqslant 2^{\tss\prime}$ we have
\beql{LS}
[\ell^{\tss\pm}_{11}(u),s^{\pm}_{ij}(v)]=0
\eeq
and
\beql{LPMSMP}
\frac{u_{\pm}-v_{\mp}}{qu_{\pm}-q^{-1}v_{\mp}}\ell_{11}^{\tss\pm}(u)s^{\mp}_{ij}(v)
=\frac{u_{\mp}-v_{\pm}}{qu_{\mp}-q^{-1}v_{\pm}}s^{\mp}_{ij}(v)\ts\ell_{11}^{\tss\pm}(u).
\eeq

\subsection{Homomorphism from $U({R}^{[n-1]})$ to $U({R}^{[n]})$}
Now we establish a relation between the algebras $U(\overline{R})$
associated with the Lie algebras $\g_{N-2}$ and $\g_{N}$. Let us denote by $\overline{R}^{\tss[n]}$ the
R-matrix associated with rank $n$.
Consider the algebra $U(\overline{R}^{\tss[n-1]})$ and let
the indices
of the generators $\ell^{\tss\pm}_{ij}[\mp m]$ range over the sets
$2\leqslant i,j\leqslant 2\pr$ and $m=0,1,\dots$, where
$i\pr=N-i+1$, as before.

\bth\label{thm:embed}
The maps $q^{\pm c/2}\mapsto q^{\pm c/2}$ and
\beql{embedgen}
\ell^{\tss\pm}_{ij}(u)\mapsto \left|\begin{matrix}
\ell^{\tss\pm}_{11}(u)&\ell^{\tss\pm}_{1j}(u)\\[0.2em]
\ell^{\tss\pm}_{i1}(u)&\boxed{\ell^{\tss\pm}_{ij}(u)}
\end{matrix}\right|,\qquad 2\leqslant i,j\leqslant 2\pr,
\eeq
define a homomorphism ${U}(\overline{R}^{\tss[n-1]})\to {U}(\overline{R}^{\tss[n]})$.
\eth

The following result can be checked directly. 
\ble
For $2\leq i,j\leq 2'$, we have
\begin{multline}\label{R23ran}
\wh{R}^{\tss[n]}_{12}(q^{-2})\wh{R}^{\tss[n]}_{34}(q^{-2})\overline{R}^{\tss[n]}_{14}(aq^{-2})
\overline{R}^{\tss[n]}_{24}(a)\overline{R}^{\tss[n]}_{13}(a)
\overline{R}^{\tss[n]}_{23}(aq^2)\ts|\ts1,i,1,j\rangle\\[0.4em]
{}=\frac{aq^{-1}-q}{aq-q^{-1}}\ts\wh{R}^{\tss[n]}_{12}(q^{-2})
\wh{R}^{\tss[n]}_{34}(q^{-2})\overline{R}^{\tss[n-1]}_{24}(a)\ts|\ts1,i,1,j\rangle
\end{multline}
and
\begin{multline}\label{lanR23}
\langle 1,i,1,j\ts|\ts\overline{R}^{\tss[n]}_{23}(aq^2)\overline{R}^{\tss[n]}_{13}(a)
\overline{R}^{\tss[n]}_{24}(a)\overline{R}^{\tss[n]}_{14}(aq^{-2})
\wh{R}^{\tss[n]}_{12}(q^{-2})\wh{R}^{\tss[n]}_{34}(q^{-2})\\[0.4em]
{}=\frac{aq^{-1}-q}{aq-q^{-1}}\ts\langle 1,i,1,j\ts|\ts\overline{R}^{\tss[n-1]}_{24}(a)
\wh{R}^{\tss[n]}_{12}(q^{-2})\wh{R}^{\tss[n]}_{34}(q^{-2}).
\end{multline}
\ele
\bpf
To prove the theorem, introduce the matrices
\ben
\Gamma^{\pm}(u)=\sum\limits_{i,j=2}^{2^{\tss\prime}} e_{ij}\otimes
\ell^{\tss\pm\tss 1 i}_{~~ 1 j}(u)\in \End \CC^{N}\ot U(\overline{R}^{\tss[n]}).
\een
Our next step is to verify that the following relations hold in
the algebra $U(\overline{R}^{\tss[n]})$:
\begin{align}\non
\overline{R}^{\tss[n-1]}(u/v)\Gamma^{\pm}_1(u)\Gamma^{\pm}_2(v)
&=\Gamma^{\pm}_2(v)\Gamma^{\pm}_1(u)\overline{R}^{\tss[n-1]}(u/v),\\[0.4em]
\frac{q^{-1}{u_+}-qv_{-}}{q{u_+}-q^{-1}v_{-}}\ts
\overline{R}^{\tss[n-1]}(uq^c/v)\Gamma^{+}_1(u)\Gamma^{-}_2(v)&
=\frac{q^{-1}{u_-}-qv_{+}}{q{u_-}-q^{-1}v_{+}}
\ts\Gamma^{-}_2(v)\Gamma^{+}_1(u)\overline{R}^{\tss[n-1]}(uq^{-c}/v).
\non
\end{align}
The calculations are quite similar in both cases so we only give details for the first relation.
The Yang--Baxter equation and the defining relations for the algebra $U(\overline{R}^{\tss[n]})$ give
\begin{multline}\non
\overline{R}^{\tss[n]}_{23}\big(\frac{u q^2}{v}\big)\overline{R}^{\tss[n]}_{13}\big(\frac uv\big)
\overline{R}^{\tss[n]}_{24}\big(\frac{u}{v}\big)\overline{R}^{\tss[n]}_{14}\big(\frac{u}{v q^{2}}\big)
\wh{R}^{\tss[n]}_{12}(q^{-2})\Lc^{\pm}_1(u)\ts \Lc^{\pm}_2(uq^{2})
\wh{R}^{\tss[n]}_{34}(q^{-2})\Lc^{\pm}_3(v)\ts \Lc^{\pm}_4(vq^{2})\\[0.4em]
=\wh{R}^{\tss[n]}_{34}(q^{-2})\Lc^{\pm}_3(v)\ts \Lc^{\pm}_4(vq^{2})
\wh{R}^{\tss[n]}_{12}(q^{-2})\Lc^{\pm}_1(u)\ts \Lc^{\pm}_2(uq^{2})
\overline{R}^{\tss[n]}_{14}\big(\frac{u}{v q^{2}}\big)\overline{R}^{\tss[n]}_{24}\big(\frac{u}{v}\big)
\overline{R}^{\tss[n]}_{13}\big(\frac uv\big)\overline{R}^{\tss[n]}_{23}\big(\frac{u q^2}{v}\big).
\end{multline}
Hence, assuming that $2\leqslant i,j,k,l\leqslant 2^{\tss\prime}$
and applying \eqref{R23ran} and \eqref{lanR23}
we get
\begin{multline}\non
\langle 1,k,1,l\ts|\ts\overline{R}^{\tss[n-1]}_{24}\big(\frac{u}{v}\big)
\wh{R}^{\tss[n]}_{12}(q^{-2})\Lc^{\pm}_1(u)\ts \Lc^{\pm}_2\big(uq^{2}\big)
\wh{R}^{\tss[n]}_{34}(q^{-2})\Lc^{\pm}_3(v)\ts \Lc^{\pm}_4(vq^{2})\ts|\ts1,i,1,j\rangle\\[0.4em]
=\langle 1,k,1,l\ts|\ts\wh{R}^{\tss[n]}_{34}(q^{-2})\Lc^{\pm}_3(v)\ts \Lc^{\pm}_4(vq^{2})
\wh{R}^{\tss[n]}_{12}(q^{-2})\Lc^{\pm}_1(u)\ts \Lc^{\pm}_2(uq^{2})
\overline{R}^{\tss[n]}_{24}\big(\frac{u}{v}\big)\ts|\ts1,i,1,j\rangle,
\end{multline}
which is equivalent to
\ben
\overline{R}^{\tss[n-1]}_{24}(u/v)\Gamma^{\pm}_2(u)\Gamma^{\pm}_4(v)
=\Gamma^{\pm}_4(v)\Gamma^{\pm}_2(u)\overline{R}^{\tss[n-1]}_{24}(u/v),
\een
as required.
Finally, set
\ben
S^{\pm}(u)=\sum\limits_{2\leqslant i,j\leqslant 2\pr}e_{ij}\ot s^{\pm}_{ij}(u).
\een
By Lemma \ref{lem:toneone},
\ben
S^{\pm}(u)=\ell^{\tss\pm}_{11}(uq^{-2})^{-1}\Ga^{\pm}(uq^{-2})
\een
and
\ben
\frac{q^{-1}u_{\pm}-qv_{\mp}}{u_{\pm}-v_{\mp}}\ts\ell_{11}^{\tss\pm}(u)\ts\Ga^{\mp}(v)
=\frac{q^{-1}u_{\mp}-qv_{\tss\pm}}{u_{\mp}-v_{\pm}}\ts\Ga^{\mp}(v)\ts\ell_{11}^{\tss\pm}(u).
\een
The above relations for the matrices $\Ga^{\pm}(u)$ imply
\ben
\bal
\overline{R}^{\tss[n-1]}(u/v)S^{\pm}_1(u)S^{\pm}_2(v)
&=S^{\pm}_2(v)S^{\pm}_1(u)\overline{R}^{\tss[n-1]}(u/v),\\[0.4em]
\overline{R}^{\tss[n-1]}(uq^{\pm c}/v)S^{\pm}_1(u)S^{\mp}_2(v)
&=S^{\mp}_2(v)S^{\pm}_1(u)\overline{R}^{\tss[n-1]}(uq^{\mp c}/v),
\eal
\een
which completes the proof.
\epf
The following generalizes
Theorem~\ref{thm:embed} by using the Sylvester theorem for quasideterminants
\cite{gr:dm}, \cite{kl:mi}; cf. the proof of the Yangian algebra \cite[Thm~3.7]{jlm:ib}.
Fix a positive integer $m$ such that
$m< n$.
Suppose that the generators $\ell_{ij}^{\tss\pm}(u)$ of the algebra $U(\overline{R}^{\tss[n-m]})$ are
labelled by the indices
$m+1\leqslant i,j\leqslant (m+1)\pr$ with $i\pr=N-i+1$ as before.

\bth\label{thm:red}
For $m\leq n-1$,
the map
\beql{redu}
\ell^{\tss\pm}_{ij}(u)\mapsto \left|\begin{matrix}
\ell^{\tss\pm}_{11}(u)&\dots&\ell^{\tss\pm}_{1m}(u)&\ell^{\tss\pm}_{1j}(u)\\
\dots&\dots&\dots&\dots\\
\ell^{\tss\pm}_{m1}(u)&\dots&\ell^{\tss\pm}_{mm}(u)&\ell^{\tss\pm}_{mj}(u)\\[0.2em]
\ell^{\tss\pm}_{i1}(u)&\dots&\ell^{\tss\pm}_{im}(u)&\boxed{\ell^{\tss\pm}_{ij}(u)}
\end{matrix}\right|,\qquad m+1\leqslant i,j\leqslant (m+1)\pr,
\eeq
defines a homomorphism $\varphi^+_m:U(\overline{R}^{\tss[n-m]})\to U(\overline{R}^{\tss[n]})$.
\qed
\eth

All the maps \eqref{redu} are compatible in the following sense, as seen by using the Sylvester theorem for quasideterminants.
Write the maps $\varphi^+_m=\psi^{(n)}_m$ for $U(\overline{R}^{\tss[n]})$, and for each admissible
$l$ we have the corresponding homomorphism
\ben
\psi^{(n-l)}_m:U(\overline{R}^{\tss[n-l-m]})\to U(\overline{R}^{\tss[n-l]})
\een
given by \eqref{redu}. Then we have the equation:
\beql{consist}
\psi^{(n)}_l\circ\psi^{(n-l)}_m=\psi^{(n)}_{l+m}.
\eeq

Let $\{a_1,\dots,a_k\}$ and $\{b_1,\dots,b_k\}$ be subsets of $\{1,\dots,N\}$, where
$a_1<a_2<\dots<a_k$ and $b_1<b_2<\dots<b_k$ such that
$a_i\neq a_j'$ and $b_i\neq b_j'$ for all $i,j$. Introduce the
{\it type $A$ quantum minors} as the matrix elements \eqref{matrelem}:
\begin{multline}\non
\ell^{\tss\pm~a_1,\dots,a_k}_{\ts\ts~~b_1,\dots,b_k}(u)
=\langle a_1,\dots,a_k\ts|\ts \wh{R}_{k-1,k}(\wh{R}_{k-2,k}\wh{R}_{k-2,k-1})
\dots(\wh{R}_{1,k}\dots\wh{R}_{1,2})\\[0.3em]
{}\times\Lc_1^{\pm}(u)\Lc_2^{\pm}(uq^2)\dots\Lc_k^{\pm}(uq^{2k-2})\ts|\ts b_1,\dots,b_k\rangle,
\end{multline}
where $\wh{R}_{ij}=\wh{R}_{ij}(q^{2(i-j)})$.  Explicitly we have:
\ben
\bal
\ell^{\tss\pm~a_1,\dots,a_k}_{\ts\ts~~b_1,\dots,b_k}(u)
&=\sum_{\sigma\in \Sym_k}(-q)^{-l(\sigma)}\ell^{\tss\pm}_{a_{\sigma(1)}b_1}(u)
\dots\ell^{\tss\pm}_{a_{\sigma(k)}b_k}(uq^{2k-2})\\
&=\sum_{\sigma\in \Sym_k}(-q)^{l(\sigma)}
\ell^{\tss\pm}_{a_{k}b_{\sigma(k)}}(uq^{2k-2})\dots\ell^{\tss\pm}_{a_{1}b_{\sigma(1)}}(u),
\eal
\een
where $l(\sigma)$ is the number of inversions of the permutation $\sigma\in \Sym_k$.
The assumption on the indices $a_i$ and $b_i$ implies that
the quantum minors satisfy the same form
as those for the quantum affine algebra $U_q(\wh\gl_n)$. By applying
$R$-matrix calculations (cf.~\cite{hm:qa}, \cite{kl:mi}, \cite[Ch.~1]{m:yc}), we have that
for $1\leqslant i,j\leqslant k$
\ben
\bal[]
[\ell^{\tss\pm}_{a_ib_j}(u),\ell^{\tss\pm~a_1,a_2,\dots,a_k}_{\ts\ts~~b_1,b_2,\dots,b_k}(v)]&=0,\\[0.4em]
 \displaystyle\prod_{a=1}^{k-1}\frac{u_{\pm}q^{-k}-v_{\mp}q^{k}}{u_{\pm}q^{1-k}-v_{\mp}q^{k-1}}
\ell^{\tss\pm}_{a_ib_j}(u)\ell^{\tss\mp~a_1,a_2,\dots,a_k}_{\ts\ts~~b_1,b_2,\dots,b_k}(v)&=
\displaystyle\prod_{a=1}^{k-1}\frac{u_{\mp}q^{-k}-v_{\pm}q^{k}}{u_{\mp}q^{1-k}-v_{\pm}q^{k-1}}
\ell^{\tss\mp~a_1,a_2,\dots,a_k}_{\ts\ts~~b_1,b_2,\dots,b_k}(v)\ell^{\tss\pm}_{a_ib_j}(u).
\eal
\een

\bco\label{cor:commu}
Under the assumptions of Theorem~\ref{thm:red} we have
\ben
\bal
\big[\ell^{\tss\pm}_{ab}(u),\varphi^+_m(\ell^{\tss\pm}_{ij}(v))\big]&=0,\\[0.4em]
\frac{u_{\pm}-v_{\mp}}{qu_{\pm}-q^{-1}v_{\mp}}\ell^{\tss\pm}_{ab}(u)\varphi^+_m(\ell^{\mp}_{ij}(v))
{}&=\frac{u_{\mp}-v_{\pm}}{qu_{\mp}-q^{-1}v_{\pm}}\varphi^+_m(\ell^{\mp}_{ij}(v))\ell^{\tss\pm}_{ab}(u),
\eal
\een
for all $1\leqslant a,b\leqslant m$ and $m+1\leqslant i,j\leqslant (m+1)\pr$.
\eco

\bpf
Both formulas are verified by the fact that
the quasideterminants and quantum minors satisfy the following type A relations:
(cf. similar ones for the Yangian \cite[Sec.~3]{jlm:ib}.
\ben
\left|\begin{matrix}
\ell^{\tss\pm}_{11}(u)&\dots&\ell^{\tss\pm}_{1m}(u)&\ell^{\tss\pm}_{1j}(u)\\
\dots&\dots&\dots&\dots\\
\ell^{\tss\pm}_{m1}(u)&\dots&\ell^{\tss\pm}_{mm}(u)&\ell^{\tss\pm}_{mj}(u)\\[0.2em]
\ell^{\tss\pm}_{i1}(u)&\dots&\ell^{\tss\pm}_{im}(u)&\boxed{\ell^{\tss\pm}_{ij}(u)}
\end{matrix}\right|={\ell\ts}^{\pm\tss1\dots\ts m}_{~~1\dots\ts m}(uq^{-2m})^{-1}\cdot
{\ell\ts}^{\pm\tss1\dots\ts m\ts i}_{~~1\dots\ts m\ts j}(uq^{-2m}).
\een\epf


\section{Gauss decomposition}\label{sec:gauss decom}
\subsection{Gaussian generators}

It is known that the matrix $L^{\pm}(u)$ of generators for $U(R^{[n]})$ admits
a unique
{\it Gauss decomposition}. Namely, $L^{\pm}(u)$ can be uniquely factored as
\beql{gaussdec}
L^{\pm}(u)=F^{\pm}(u)H^{\pm}(u)E^{\pm}(u).
\eeq
where
\ben
F^{\pm}(u)=\begin{bmatrix}
1&0&\dots&0\ts\\
f^{\pm}_{21}(u)&1&\dots&0\\
\vdots&\vdots&\ddots&\vdots\\
f^{\pm}_{2n\ts1}(u)&f^{\pm}_{2n\ts2}(u)&\dots&1
\end{bmatrix},
\qquad
E^{\pm}(u)=\begin{bmatrix}
\ts1&e^{\pm}_{12}(u)&\dots&e^{\pm}_{1\ts 2n}(u)\ts\\
\ts0&1&\dots&e^{\pm}_{2\ts 2n}(u)\\
\vdots&\vdots&\ddots&\vdots\\
0&0&\dots&1
\end{bmatrix},
\een
and $H^{\pm}(u)=\diag\ts\big[h^{\pm}_1(u),\dots,h^{\pm}_{2n}(u)\big]$. Moreover, the entries can be expressed in terms of
quasideterminants.
as follows \cite{gr:dm}, \cite{gr:tn}:
\beql{hmqua}
h^{\pm}_i(u)=\begin{vmatrix} l^{\pm}_{1\tss 1}(u)&\dots&l^{\pm}_{1\ts i-1}(u)&l^{\pm}_{1\tss i}(u)\\
                          \vdots&\ddots&\vdots&\vdots\\
                         l^{\pm}_{i-1\ts 1}(u)&\dots&l^{\pm}_{i-1\ts i-1}(u)&l^{\pm}_{i-1\ts i}(u)\\[0.2em]
                         l^{\pm}_{i\tss 1}(u)&\dots&l^{\pm}_{i\ts i-1}(u)&\boxed{l^{\pm}_{i\tss i}(u)}\\
           \end{vmatrix},\qquad i=1,\dots,2n,
\eeq
\beql{eijmlqua}
e^{\pm}_{ij}(u)=h^{\pm}_i(u)^{-1}\ts\begin{vmatrix} l^{\pm}_{1\tss 1}(u)&\dots&
l^{\pm}_{1\ts i-1}(u)&l^{\pm}_{1\ts j}(u)\\
                          \vdots&\ddots&\vdots&\vdots\\
                         l^{\pm}_{i-1\ts 1}(u)&\dots&l^{\pm}_{i-1\ts i-1}(u)&l^{\pm}_{i-1\ts j}(u)\\[0.2em]
                         l^{\pm}_{i\tss 1}(u)&\dots&l^{\pm}_{i\ts i-1}(u)&\boxed{l^{\pm}_{i\tss j}(u)}\\
           \end{vmatrix}
\eeq
and
\beql{fijlmqua}
f^{\pm}_{ji}(u)=\begin{vmatrix} l^{\pm}_{1\tss 1}(u)&\dots&l^{\pm}_{1\ts i-1}(u)&l^{\pm}_{1\tss i}(u)\\
                          \vdots&\ddots&\vdots&\vdots\\
                         l^{\pm}_{i-1\ts 1}(u)&\dots&l^{\pm}_{i-1\ts i-1}(u)&l^{\pm}_{i-1\ts i}(u)\\[0.2em]
                         l^{\pm}_{j\ts 1}(u)&\dots&l^{\pm}_{j\ts i-1}(u)&\boxed{l^{\pm}_{j\tss i}(u)}\\
           \end{vmatrix}\ts h^{\pm}_i(u)^{-1}
\eeq
for $1\leqslant i<j\leqslant 2n$. The same formulas hold for the entries of the respective triangular matrices $\Fc^{\pm}(u)$ and
$\Ec^{\pm}(u)$ and the diagonal matrices $$\Hc^{\pm}(u)=\diag\ts[\h^{\pm}_{1}(u),\dots,\h^{\pm}_{2n}(u)]$$
in terms of the formal series $\ell^{\tss\pm}_{ij}(u)$ via the Gauss decomposition
\ben
\Lc^{\pm}(u)=\Fc^{\pm}(u)\ts\Hc^{\pm}(u)\ts\Ec^{\pm}(u)
\een
for the algebra $U(\overline{R}^{\tss[n]})$. We will denote by
$\e_{ij}(u)$ and $\f_{ji}(u)$ the entries of the respective matrices
$\Ec^{\pm}(u)$ and $\Fc^{\pm}(u)$ for $i<j$.

The following is immediate from the Gaussian generators.
\bpr\label{prop:corrgauss}
Under the homomorphism $U(\overline{R})\to \Hc_q(n)\ot_{\CC[q^c,\ts q^{-c}]}U(R)$
provided by Proposition~\ref{prop:homheis} we have
\ben
\bal
  \e^{\pm}_{ij}(u) & \mapsto e^{\pm}_{ij}(u), \\
  \f^{\pm}_{ij}(u) & \mapsto f^{\pm}_{ij}(u), \\
  \h^{\pm}_{i}(u) &\mapsto \exp\sum\limits _{k=1}^{\infty}\be_{\mp k}u^{\mp k}\cdot h^{\pm}_{i}(u).
\eal
\een
\epr
%

\subsection{Images of the generators under the homomorphism $\psi_m$}
For $0\leqslant m< n$, we use the superscript
$[n-m]$ to indicate the submatrices with rows and columns labelled by
$m+1,m+2,\dots,(m+1)'$. Thus we set
\ben
\Fc^{\pm[n-m]}(u)=\begin{bmatrix}
1&0&\dots&0\ts\\
\f^{\pm}_{m+2\ts m+1}(u)&1&\dots&0\\
\vdots&\ddots&\ddots&\vdots\\
\f^{\pm}_{(m+1)'\tss m+1}(u)&\dots&\f^{\pm}_{(m+1)'\ts (m+2)'}(u)&1
\end{bmatrix},
\een
\ben
\Ec^{\pm[n-m]}(u)=\begin{bmatrix} 1&\e^{\pm}_{m+1\tss m+2}(u)&\ldots&\e^{\pm}_{m+1\tss (m+1)'}(u)\\
                        0&1&\ddots &\vdots\\
                         \vdots&\vdots&\ddots&\e^{\pm}_{(m+2)'\tss(m+1)'}(u)\\
                         0&0&\ldots&1\\
           \end{bmatrix}
\een
and $\Hc^{\pm[n-m]}(u)=\diag\ts\big[\h^{\pm}_{m+1}(u),\dots,\h^{\pm}_{(m+1)'}(u)\big]$. Also we set that
\ben
\Lc^{\pm[n-m]}(u)=\Fc^{\pm[n-m]}(u)\ts \Hc^{\pm[n-m]}(u)\ts \Ec^{\pm[n-m]}(u)=(\ell^{\tss\pm[n-m]}_{ij}(u)).
\een
The following result can be easily checked by the same argument as in the Yangian case (cf. \cite[Prop.~4.1]{jlm:ib}):
\bpr\label{prop:gauss-consist}
The series $\ell^{\tss\pm[n-m]}_{ij}(u)$ coincides with the image of
the generator series $\ell^{\tss\pm}_{ij}(u)$
of the extended quantum affine algebra $U(\overline{R}^{\tss[n-m]})$
under the homomorphism \eqref{redu}:
\ben
\ell^{\pm[n-m]}_{ij}(u)=\psi_m\big(\ell^{\pm}_{ij}(u)\big),\qquad m+1\leqslant i,j\leqslant (m+1)'.
\een
\epr
%

The following can be shown by using Proposition~\ref{prop:gauss-consist}.
\bco\label{cor:guass-embed}
The following relations hold in $U(\overline{R}^{\tss[n]})${\rm :}
\beql{subRTT}
\overline{R}^{\tss[n-m]}_{12}(u/v)\ts \Lc^{\pm[n-m]}_1(u)\ts \Lc^{\pm[n-m]}_2(v)
=\Lc^{\pm[n-m]}_2(v)\ts \Lc^{\pm[n-m]}_1(u)\ts \overline{R}^{\tss[n-m]}_{12}(u/v),
\eeq
\beql{subRTT2}
\overline{R}^{\tss[n-m]}_{12}(u_{+}/v_{-})\ts \Lc^{+[n-m]}_1(u)\ts \Lc^{-[n-m]}_2(v)
=\Lc^{-[n-m]}_2(v)\ts \Lc^{+[n-m]}_1(u)\ts \overline{R}^{\tss[n-m]}_{12}(u_{-}/v_{+}).
\eeq
\eco
%

The following proposition is the same as in the untwisted classical cases (BCD) \cite{jlm:ib,jlm:iC,jlm:ibC}.
\bpr\label{prop:relmone}
Suppose that
$m+1\leqslant j,k,l\leqslant (m+1)'$ and $j\neq l'$. Then the following relations hold
in $U(\overline{R}^{\tss[n]})${\rm :} if $j=l$ then
\begin{align}\label{ELMPj=l}
\e_{mj}^{\pm}(u)\ell^{\mp [n-m]}_{kl}(v)
&=\frac{qu_{\mp}-q^{-1}v_{\pm}}{u_{\mp}-v_{\pm}}\ell^{\mp [n-m]}_{kj}(v)\e_{ml}^{\pm}(u)
-\frac{(q-q^{-1})u_{\mp}}{u_{\mp}-v_{\pm}}\ell^{\mp [n-m]}_{kj}(v)\e_{mj}^{\mp}(v),\\
\e_{mj}^{\pm}(u)\ell^{\tss\pm [n-m]}_{kl}(v)
&=\frac{qu-q^{-1}v}{u-v}\ell^{\tss\pm [n-m]}_{kj}(v)\e_{ml}^{\pm}(u)
-\frac{(q-q^{-1})u}{u-v}\ell^{\tss\pm [n-m]}_{kj}(v)\e_{mj}^{\pm}(v);
\non
\end{align}
if $j<l$ then
\begin{align}\label{ELMPj<l}
[\e^{\pm}_{mj}(u),\ell^{\mp [n-m]}_{kl}(v)]
&=\frac{(q-q^{-1})v_{\pm}}{u_{\mp}-v_{\pm}}\ell^{\mp [n-m]}_{kj}(v)\e^{\pm}_{ml}(u)-
\frac{(q-q^{-1})u_{\mp}}{u_{\mp}-v_{\pm}}\ell^{\mp [n-m]}_{kj}(v)\e^{\mp}_{ml}(v),\\
[\e^{\pm}_{mj}(u),\ell^{\tss\pm [n-m]}_{kl}(v)]
&=\frac{(q-q^{-1})v}{u-v}\ell^{\tss\pm [n-m]}_{kj}(v)\e^{\pm}_{ml}(u)-
\frac{(q-q^{-1})u}{u-v}\ell^{\tss\pm [n-m]}_{kj}(v)\e^{\pm}_{ml}(v);
\non
\end{align}
if $j>l$ then
\ben
\bal[]
[\e^{\pm}_{mj}(u),\ell^{\mp [n-m]}_{kl}(v)]
&=\frac{(q-q^{-1})u_{\mp}}{u_{\mp}-v_{\pm}}\ell^{\mp [n-m]}_{kj}(v)(\e^{\pm}_{ml}(u)-\e_{ml}^{\mp}(v)),\\
[\e^{\pm}_{mj}(u),\ell^{\tss\pm [n-m]}_{kl}(v)]
&=\frac{(q-q^{-1})u}{u-v}\ell^{\tss\pm [n-m]}_{kj}(v)(\e^{\pm}_{ml}(u)-\e_{ml}^{\pm}(v)).
\eal
\een
\epr

Similar arguments prove the following counterpart of Proposition~\ref{prop:relmone}
involving the generator series $\f^{\pm}_{ji}(u)$.

\bpr\label{prop:relmonf}
Suppose that
$m+1\leqslant j,k,l\leqslant (m+1)'$ and $j\neq k'$. Then the following relations hold
in $U(\overline{R}^{\tss[n]})${\rm :} if $j=k$ then
\ben
\bal
\f_{jm}^{\pm}(u)\ell_{jl}^{\mp [n-m]}(v)
&=\frac{u_{\pm}-v_{\mp}}{qu_{\pm}-q^{-1}v_{\mp}}\ell_{jl}^{\mp [n-m]}(v)\f_{jm}^{\pm}(u)
+\frac{(q-q^{-1})v_{\mp}}{qu_{\pm}-q^{-1}v_{\mp}}\ts\f_{jm}^{\mp}(v)\ell_{jl}^{\mp [n-m]}(v),\\
\f_{jm}^{\pm}(u)\ell_{jl}^{\pm [n-m]}(v)
&=\frac{u-v}{qu-q^{-1}v}\ell_{jl}^{\pm [n-m]}(uv)\ts\f_{jm}^{\pm}(u)
+\frac{(q-q^{-1})v}{qu-q^{-1}v}\f_{jm}^{\pm}(v)\ell_{jl}^{\pm [n-m]}(v);
\eal
\een
if $j<k$ then
\ben
\bal[]
[\f_{jm}^{\pm}(u),\ell_{kl}^{\mp [n-m]}(v)]
&=\frac{(q-q^{-1})v_{\mp}}{u_{\pm}-v_{\mp}}\f_{km}^{\mp}(v)\ell_{jl}^{\mp [n-m]}(v)
-\frac{(q-q^{-1})u_{\pm}}{u_{\pm}-v_{\mp}}\f_{km}^{\pm}(u)\ell_{jl}^{\mp [n-m]}(v),\\
[\f_{jm}^{\pm}(u),\ell_{kl}^{\pm [n-m]}(v)]
&=\frac{(q-q^{-1})v}{u-v}\f_{km}^{\pm}(v)\ell_{jl}^{\pm [n-m]}(v)
-\frac{(q-q^{-1})u}{u-v}\f_{km}^{\pm}(u)\ell_{jl}^{\pm [n-m]}(v);
\eal
\een
if $j>k$ then
\ben
\bal[]
[\f_{jm}^{\pm}(u),\ell_{kl}^{\mp [n-m]}(v)]
&=\frac{(q-q^{-1})v_{\mp}}{u_{\pm}-v_{\mp}}(\f_{km}^{\mp}(v)-\f_{km}^{\pm}(u))
\ell_{jl}^{\mp [n-m]}(v),\\
[\f_{jm}^{\pm}(u),\ell_{kl}^{\pm [n-m]}(v)]
&=\frac{(q-q^{-1})v}{u-v}(\f_{km}^{\pm}(v)-\f_{km}^{\pm}(u))\ell_{jl}^{\pm [n-m]}(v).
\eal
\vspace{-0.7cm}
\een
\qed
\epr

\section{Isomorphism between $U^{ext}_q(A^{(2)}_{2n-1})$ and $U({R}^{[n]})$}\label{Sec:guass-decomposition}
In this section, we derive the commuting relations among Gaussian generators in $U({\overline{R}}^{[n]})$
and $U({R}^{[n]})$. By Proposition \ref{prop:corrgauss} only the case of $U({\overline{R}}^{[n]})$ needs to be considered.

To do so, we will frequently employ the following Laurent series with coefficients in
the respective algebras $U(R^{[n]})$:
\begin{align}
\label{Xi}
X^{+}_i(u)=e^{+}_{ii+1}(u_{+})-e_{ii+1}^{-}(u_{-}),
\qquad X^{-}_i(u)=f^{+}_{i+1,i}(u_{-})-f^{-}_{i+1,i}(u_{+})
\end{align}
for $i=1,\dots,n-1$
 and
\begin{align}
\label{Xnp}
X^{+}_n(u)&=u\left(e^{+}_{n,n+1}(u_{+})-e^{-}_{n,n+1}(u_{-})\right)\\
\label{Xnm}
X^{-}_n(u)&=\frac{1}{u}\left(f^{+}_{n+1,n}(u_{-})-f^{-}_{n+1,n}(u_{+})\right)
\end{align}
and then $\Xc^{\pm}_i(u)$ are similarly defined for $U(\overline{R}^{\tss[n]})$.

\subsection{Type A relations}

Due to the observation made in Remark~\ref{rem:gln}
and the quasideterminant formulas \eqref{hmqua}, \eqref{eijmlqua}
and \eqref{fijlmqua}, some of the relations between the Gaussian generators follow
from those for the quantum affine algebra
$U_q(\wh{\mathfrak{gl}}_n)$; see \cite{df:it}. To reproduce them, set
\ben
\Lc^{A\pm}(u)=\sum_{i,j=1}^n e_{ij}\otimes \ell^{\tss\pm}_{ij}(u)
\een
and consider the $R$-matrix used in \cite{df:it}: 
\begin{multline}
R_A(u)=\sum_{i=1}^n e_{ii}\otimes e_{ii}+\frac{u-1}{qu-q^{-1}}
\sum_{i\neq j} e_{ii}\otimes e_{jj}\\
+\frac{q-q^{-1}}{qu-q^{-1}}\sum_{i>j}e_{ij}\otimes e_{ji}
+\frac{(q-q^{-1})u}{qu-q^{-1}}\sum_{i<j}e_{ij}\otimes e_{ji}.
\label{rtypea}
\end{multline}
Comparing this with the $R$-matrix \eqref{rbar}, we immediately see the
relations in the algebra
$U(\overline{R}^{\tss[n]})${\rm:}
\ben
\bal
R_A(u/w)\Lc^{A\pm}_{1}(u)\Lc^{A\pm}_2(v)&=\Lc^{A\pm}_2(v)\Lc^{A\pm}_{1}(u){R}_A(u/v),\\[0.3em]
R_A(uq^c/v)\Lc^{A+}_1(u)\Lc^{A-}_2(v)&=\Lc^{A-}_2(v)\Lc^{A+}_1(u){R}_A(uq^{-c}/v).
\eal
\een
Therefore we get the following relations for the Gaussian generators (cf.
\cite{df:it}), where we use the notation similar to \eqref{Xi}.

\bpr\label{TypeArelation}
In
the algebra $U(\overline{R}^{\tss[n]})$ we have
\ben
\h^{\pm}_i(u)\h^{\pm}_j(v)=\h^{\pm}_j(v)\h^{\pm}_i(u),\qquad
\h^{\pm}_i(u)\h^{\mp}_i(v)=\h^{\mp}_i(v)\h^{\pm}_i(u), \quad\text{for}\quad 1\leq i,j\leq n.
\een
\ben
\frac{u_{\pm}-v_{\mp}}{qu_{\pm}-q^{-1}v_{\mp}}\h^{\pm}_i(u)\h^{\mp}_j(v)=
\frac{u_{\mp}-v_{\pm}}{qu_{\mp}-q^{-1}v_{\pm}}\h^{\mp}_j(v)\h^{\pm}_i(u), \quad\text{for}\quad 1\leq i<j\leq n.
\een
Moreover,
\ben
\bal
\h_{i}^{\pm}(u)\Xc_{j}^{+}(v)
&=\frac{u_{\mp}-v}{q^{(\ep_i,\alpha_j)}u_{\mp}-q^{-(\ep_i,\alpha_j)}v} \Xc_{j}^{+}(v)\h_{i}^{\pm}(u),\\
\h_{i}^{\pm}(u)\Xc_{j}^{-}(v)
&=\frac{q^{(\ep_i,\alpha_j)}u_{\pm}-q^{-(\ep_i,\alpha_j)}v}{u_{\pm}-v} \Xc_{j}^{-}(v)\h_{i}^{\pm}(u), \quad\text{for}\quad 1\leq i\leq n, 1\leq j< n.
\eal
\een
and
\beq\label{eq:xixj}
(u-q^{\pm (\alpha_i,\alpha_j)}v)\Xc_{i}^{\pm}(uq^i)\Xc_{j}^{\pm}(vq^j)
=(q^{\pm (\alpha_i,\alpha_j)}u-v) \Xc_{j}^{\pm}(vq^j)\Xc_{i}^{\pm}(uq^i),
\eeq
\ben
[\Xc_i^{+}(u),\Xc_j^{-}(v)]=\delta_{ij}(q-q^{-1})
\Big(\delta\big(u\ts q^{-c}/v\big)\ts\h_i^{-}(v_+)^{-1}\h_{i+1}^{-}(v_+)
-\delta\big(u\ts q^{c}/v\big)\ts\h_i^{+}(u_+)^{-1}\h_{i+1}^{+}(u_+)\Big)
\een
for $1\leq i,j<n$,
together with the Serre relations
\beql{rel:serrextypeB}
\sum_{\pi\in \Sym_{r}}\sum_{l=0}^{r}(-1)^l{{r}\brack{l}}_{q_i}
  \Xc^{\pm}_{i}(u_{\pi(1)})\dots \Xc^{\pm}_{i}(u_{\pi(l)})
  \Xc^{\pm}_{j}(v)\tss \Xc^{\pm}_{i}(u_{\pi(l+1)})\dots \Xc^{\pm}_{i}(u_{\pi(r)})=0,
\eeq
for all $1\leq i\neq j<n$, here $r=1-A_{ij}$.
\qed
\epr

\bre\label{re:anotherA}
Consider the inverse matrices $\Lc^{\pm}(u)^{-1}=[\ell^{\tss\pm}_{ij}(u)']_{i,j=1,\dots,2n}$.
By the defining relations
\eqref{gen rel1} and \eqref{gen rel2}, we have
\ben
\bal
\Lc^{\pm}_{1}(u)^{-1}\Lc^{\pm}_2(v)^{-1}\overline{R}^{\tss[n]}(u/v)
&=\overline{R}^{\tss[n]}(u/v)\Lc^{\pm}_2(v)^{-1}\Lc^{\pm}_{1}(u)^{-1},\\
\Lc^{-}_2(v)^{-1}\Lc^{+}_1(u)^{-1}\overline{R}^{\tss[n]}({u}q^{c}/v)
&=\overline{R}^{\tss[n]}({u}q^{-c}/v)\Lc^{-}_2(v)^{-1}\Lc^{+}_1(u)^{-1}.
\eal
\een
So we can obtain another family of generators of the algebra $U(\overline{R}^{\tss[n]})$
satisfying the defining relations of $U_q(\wh{\gl}_n)$. Namely, the same relations
are also satisfied by
the coefficients of the series
$\ell^{\tss\pm}_{ij}(u)'$ with $i,j=n',\dots,1'$. In particular,
by taking the inverse matrices, we get a Gauss decomposition
for the matrix $[\ell^{\tss\pm}_{ij}(u)']_{i,j=n',\dots,1'}$ from
the Gauss decomposition of the matrix $\Lc^{\pm}(u)$.
\qed
\ere

\subsection{Relations for the last root generators}
In this subsection, we will study the relations among the last
root generators in $U(\overline{R})$: $\Xc^{\pm}_n(u)$, $\h^{\pm}_{n}(u)$ and $\h^{\pm}_{n+1}(u)$.
By Corollary \ref{cor:guass-embed}, we have
\beql{subRTTl}
\overline{R}^{[1]}_{12}(u/v)\ts \Lc^{\pm[1]}_1(u)\ts {\Lc}^{\pm[1]}_2(v)
={\Lc}^{\pm[1]}_2(v)\ts {\Lc}^{\pm[1]}_1(u)\ts \overline{R}^{[1]}_{12}(u/v),
\eeq
\beql{subRTT2}
\overline{R}^{[1]}_{12}(u_{\pm}/v_{\mp})\ts {\Lc}^{\pm[1]}_1(u)\ts {\Lc}^{\mp[1]}_2(v)
={\Lc}^{\mp[1]}_2(v)\ts {\Lc}^{\pm[1]}_1(u)\ts \overline{R}^{[1]}_{12}(u_{\mp}/v_{\pm}),
\eeq
where
\beq
\overline{R}^{[1]}(u)=\sum_{i=n}^{n+1} e_{ii}\otimes e_{ii}+\frac{u^2-1}{(u^2-q^{-4})q^2}\sum_{i\neq j}e_{ii}\otimes e_{jj}
+\frac{1-q^{-4}}{u^2-q^{-4}}u\sum_{i\neq j}e_{ij}\otimes e_{ji}.
\eeq

\bpr
In $U(\overline{R})$, the following relations hold for $i,j=n,n+1$
\beq
\h^{\pm}_i(u)\h_j^{\pm}(v)=\h_j^{\pm}(v)\h^{\pm}_i(u),
\eeq
\beq
\h^{\pm}_i(u)\h_i^{\mp}(v)=\h_i^{\mp}(v)\h^{\pm}_i(u),
\eeq

\beql{eq:hnhn+1pm}
\frac{(u_{\pm}/v_{\mp})^{2}-1}{q^2(u_{\pm}/ v_{\mp})^2-q^{-2}}\h^{\pm}_n(u)\h_{n+1}^{\mp}(v)=
\frac{(u_{\mp}/v_{\pm})^{2}-1}{q^2(u_{\mp}/ v_{\pm})^2-q^{-2}}\h_{n+1}^{\mp}(v)\h^{\pm}_{n}(u).
\eeq

\beql{eq:hnxnp}
\h_{n}^{\pm}(u)\Xc^{+}_n(v)=\frac{(u/v_{\pm})^2-1}{q^2(u/v_{\pm})^2-q^{-2}}\Xc^{+}_n(v)\h_{n}^{\pm}(u),
\eeq

\beq
\Xc^{-}_{n}(v)\h_{n}^{\pm}(u)=\frac{(u/v_{\mp})^2-1}{q^2(u/v_{\mp})^2-q^{-2}}\h_{n}^{\pm}(u)\Xc^{-}_{n}(v).
\eeq

\beql{eq:hn+1xnp}
\h_{n+1}^{\pm}(u)\Xc^{+}_n(v)=\frac{(u/v_{\pm})^2-1}{q^{-2}(u/v_{\pm})^2-q^{2}}\Xc^{+}_n(v)\h_{n+1}^{\pm}(u),
\eeq

\beq
\Xc^{-}_{n}(v)\h_{n+1}^{\pm}(u)=\frac{(u/v_{\mp})^2-1}{q^{-2}(u/v_{\mp})^2-q^{2}}\h_{n+1}^{\pm}(u)\Xc^{-}_{n}(v).
\eeq

\beql{Eq:XnXnp}
\Xc_{n}^{+}(u)\Xc_{n}^{+}(v)=\frac{q^{2}(u/v)^2-q^{-2}}{q^{-2}(u/v)^2-q^2}\Xc_{n}^{+}(v)\Xc_{n}^{+}(u),
\eeq
\beql{Eq:XnXnm}
\Xc_{n}^{-}(u)\Xc_{n}^{-}(v)=\frac{q^{-2}(u/v)^2-q^{2}}{q^{2}(u/v)^2-q^{-2}}\Xc_{n}^{-}(v)\Xc_{n}^{-}(u).
\eeq
\epr

\bpf
We only prove \eqref{eq:hnhn+1pm}, \eqref{eq:hnxnp}, \eqref{eq:hn+1xnp} and \eqref{Eq:XnXnp}, as the other relations can be treated similarly.

Relation \eqref{subRTT2} can be rewritten as:
\beq
{\Lc}^{\mp[1]}_2(v)^{-1}\overline{R}^{[1]}_{12}(u_{\pm}/v_{\mp}){\Lc}^{\pm[1]}_1(u)
={\Lc}^{\pm[1]}_1(u)\overline{R}^{[1]}_{12}(u_{\mp}/v_{\pm}){\Lc}^{\mp[1]}_2(v)^{-1}.
\eeq
Thus, we get
\beq
\frac{(u_{\pm}/v_{\mp})^{2}-1}{q^2(u_{\pm}/ v_{\mp})^2-q^{-2}}\ell^{\pm \ts [1]}_{n,n}(u)\ell_{n+1,n+1}^{\mp\ts [1]}(v)'=
\frac{(u_{\mp}/v_{\pm})^{2}-1}{q^2(u_{\mp}/ v_{\pm})^2-q^{-2}}\ell_{n+1,n+1}^{\mp\ts [1]}(v)'\ell^{\pm \ts [1]}_{n,n}(u),
\eeq
where $\ell^{\pm}_{ij}(u)'$ is the $(i,j)-$ entry of $\Lc^{\pm}(u)^{-1}$.
Then we can get \eqref{eq:hnhn+1pm} by using the Guass decomposition of  $\Lc^{\pm}(u)$.

By \eqref{subRTT2}, we have
\ben
\ell^{\pm \ts [1]}_{n,n}(u)\ell_{n,n+1}^{\mp\ts [1]}(v)=
\frac{(u_{\mp}/v_{\pm})^2-1}{q^2(u_{\mp}/v_{\pm})^2-q^{-2}}\ell_{n,n+1}^{\mp\ts [1]}(v)\ell^{\pm \ts [1]}_{n,n}(u)
+\frac{(q^2-q^{-2})u_{\mp}/v_{\pm}}{q^2(u_{\mp}/v_{\pm})^2-q^{-2}}\ell_{n,n}^{\mp\ts [1]}(v)\ell^{\pm \ts [1]}_{n,n+1}(u),
\een
which is equivalent to the following form written in Gaussian generators:
\ben
\h^{\pm}_{n}(u)\h^{\mp}_n(v)\e_{n,n+1}^{\mp}(v)=
\frac{(u_{\mp}/v_{\pm})^2-1}{q^2(u_{\mp}/v_{\pm})^2-q^{-2}}\h^{\mp}_n(v)\e_{n,n+1}^{\mp}(v)\h^{\pm}_{n}(u)
+\frac{(q^2-q^{-2})u_{\mp}/v_{\pm}}{q^2(u_{\mp}/v_{\pm})^2-q^{-2}}\h_{n}^{\mp}(v)\h^{\pm}_{n}(u)\e^{\pm}_{n,n+1}(u).
\een
Furthermore, by \eqref{subRTT2} we have
\beq
\ell^{\pm \ts [1]}_{n,n}(u)\ell_{n,n}^{\mp\ts [1]}(v)=
\ell_{n,n}^{\mp\ts [1]}(v)\ell^{\pm \ts [1]}_{n,n}(u),
\eeq
which is equivalent to $\h^{\pm}_{n}(u)\h^{\mp}_n(v)=\h^{\mp}_n(v)\h^{\pm}_{n}(u)$.
Thus, by the invertibility of $\h^{\mp}_n(v)$, we have
\beql{eq:hnenpm}
\h^{\pm}_{n}(u)\e_{n,n+1}^{\mp}(v)=
\frac{(u_{\mp}/v_{\pm})^2-1}{q^2(u_{\mp}/v_{\pm})^2-q^{-2}}\e_{n,n+1}^{\mp}(v)\h^{\pm}_{n}(u)
+\frac{(q^2-q^{-2})u_{\mp}/v_{\pm}}{q^2(u_{\mp}/v_{\pm})^2-q^{-2}}\h^{\pm}_{n}(u)\e^{\pm}_{n,n+1}(u).
\eeq
Similarly, we also have
\ben
\h^{\pm}_{n}(u)\e_{n,n+1}^{\pm}(v)=
\frac{(u/v)^2-1}{q^2(u/v)^2-q^{-2}}\e_{n,n+1}^{\pm}(v)\h^{\pm}_{n}(u)
+\frac{(q^2-q^{-2})u/v}{q^2(u/v)^2-q^{-2}}\h^{\pm}_{n}(u)\e^{\pm}_{n,n+1}(u).
\een
Thus, \eqref{eq:hnxnp} follows by using the definition of $\Xc_n^{+}(u)$.

Note that \eqref{subRTT2} has the following equivalent form:
\beq
{\Lc}^{\pm[1]}_1(u)^{-1}{\Lc}^{\mp[1]}_2(v)^{-1}\overline{R}^{[1]}_{12}(u_{\pm}/v_{\mp})
=\overline{R}^{[1]}_{12}(u_{\mp}/v_{\pm}){\Lc}^{\mp[1]}_2(v)^{-1}{\Lc}^{\pm[1]}_1(u)^{-1},
\eeq
so we have that
\ben
\ell^{\pm \ts [1]}_{n+1,n+1}(u)'\ell_{n+1,n+1}^{\mp\ts [1]}(v)'=
\ell_{n+1,n+1}^{\mp\ts [1]}(v)'\ell^{\pm \ts [1]}_{n+1,n+1}(u)',
\een
and
\ben
\bal
\ell^{\pm \ts [1]}_{n,n+1}(u)'\ell_{n+1,n+1}^{\mp\ts [1]}(v)'&=
\frac{(u_{\mp}/v_{\pm})^2-1}{q^2(u_{\mp}/v_{\pm})^2-q^{-2}}\ell_{n+1,n+1}^{\mp\ts [1]}(v)'\ell^{\pm \ts [1]}_{n,n+1}(u)'\\
&+\frac{(q^2-q^{-2})u_{\mp}/v_{\pm}}{q^2(u_{\mp}/v_{\pm})^2-q^{-2}}\ell_{n,n+1}^{\mp\ts [1]}(v)'\ell^{\pm \ts [1]}_{n+1,n+1}(u)',
\eal
\een
which are equivalent to the following Gaussian generator form:
\ben
\h^{\pm}_{n+1}(u)\h_{n+1}^{\mp}(v)=
\h_{n+1}^{\mp}(v)\h^{\pm}_{n+1}(u),
\een
\ben
\bal
\e^{\pm}_{n,n+1}(u)\h^{\pm}_{n+1}(u)^{-1}\h_{n+1}^{\mp}(v)^{-1}&=
\frac{(u_{\mp}/v_{\pm})^2-1}{q^2(u_{\mp}/v_{\pm})^2-q^{-2}}\h_{n+1}^{\mp}(v)^{-1}\e^{\pm}_{n,n+1}(u)\h^{\pm}_{n+1}(u)^{-1}\\
&+\frac{(q^2-q^{-2})u_{\mp}/v_{\pm}}{q^2(u_{\mp}/v_{\pm})^2-q^{-2}}
\e_{n,n+1}^{\mp}(v)\h_{n+1}^{\mp}(v)^{-1}\h^{\pm}_{n+1}(u)^{-1}.
\eal
\een
Then we get that
\ben
\h_{n+1}^{\pm}(u)\e_{n,n+1}^{\mp}(v)=\frac{(u_{\mp}/v_{\pm})^2-1}{q^{-2}(u_{\mp}/v_{\pm})^2-q^{2}}\e_{n,n+1}^{\mp}(v)\h_{n+1}^{\pm}(u)
+\frac{(q^{-2}-q^{2})u_{\mp}/v_{\pm}}{q^{-2}(u_{\mp}/v_{\pm})^2-q^{2}}\h_{n+1}^{\pm}(u)\e_{n,n+1}^{\pm}(u).
\een
Similarly, we also have
\ben
\h_{n+1}^{\pm}(u)\e_{n,n+1}^{\pm}(v)=\frac{(u/v)^2-1}{q^{-2}(u/v)^2-q^{2}}\e_{n,n+1}^{\pm}(v)\h_{n+1}^{\pm}(u)
+\frac{(q^{-2}-q^{2})u/v}{q^{-2}(u/v)^2-q^{2}}\h_{n+1}^{\pm}(u)\e_{n,n+1}^{\pm}(u).
\een
Therefore, by the definition of $\Xc_n^{+}(u)$ we have proved \eqref{eq:hn+1xnp}.

From \eqref{subRTT2}, we have
$\ell^{\pm \ts [1]}_{n,n+1}(u)\ell_{n,n+1}^{\mp\ts [1]}(v)=\ell_{n,n+1}^{\mp \ts [1]}(v)\ell^{\pm\ts [1]}_{n,n+1}(u).$
Then using \eqref{eq:hnenpm}, we have
\ben
\bal
\e_{n,n+1}^{\pm}(u)\e_{n,n+1}^{\mp}(v)&-
\frac{q^2(u_{\mp}/v_{\pm})^2-q^{-2}}{q^{-2}(u_{\mp}/v_{\pm})^2-q^{2}}\e_{n,n+1}^{\mp}(v)\e_{n,n+1}^{\pm}(u)\\
&=-\frac{(q^2-q^{-2})u_{\mp}/v_{\pm}}{q^{-2}(u_{\mp}/v_{\pm})^2-q^{2}}\e_{n,n+1}^{\pm}(u)^2
-\frac{(q^2-q^{-2})u_{\mp}/v_{\pm}}{q^{-2}(u_{\mp}/v_{\pm})^2-q^{2}}\e_{n,n+1}^{\mp}(v)^2.
\eal
\een
Similarly, we have
\ben
\bal
\e_{n,n+1}^{\pm}(u)\e_{n,n+1}^{\pm}(v)&-
\frac{q^2(u/v)^2-q^{-2}}{q^{-2}(u/v)^2-q^{2}}\e_{n,n+1}^{\pm}(v)\e_{n,n+1}^{\pm}(u)\\
&=-\frac{(q^2-q^{-2})u/v}{q^{-2}(u/v)^2-q^{2}}\e_{n,n+1}^{\pm}(u)^2
-\frac{(q^2-q^{-2})u/v}{q^{-2}(u/v)^2-q^{2}}\e_{n,n+1}^{\mp}(v)^2.
\eal
\een
Thus, \eqref{Eq:XnXnp} follows by the definition of $\Xc_n^{+}(u)$.
\epf

\bpr
In $U(\overline{R})$, we have
\beql{eq:xnxnpm}
\left[\Xc_{n}^{+}(u), \Xc_{n}^{-}(v)\right]
=(q_n-q_{n}^{-1})\Bigg\{\delta\big((u_{-}/v_{+})^2\big)q^{c}\h_{n+1}^{-}(v)\h_{n}^{-}(v)^{-1}-
\delta\big((u_{+}/v_{-})^2\big)q^{-c}\h_{n+1}^{+}(u)\h_{n}^{+}(u)^{-1}\Bigg\}
\eeq
\epr

\bpf
It follows from \eqref{subRTT2} that
\beql{eq:lnn+1ln+1n}
\bal
&\frac{(u_{+}/v_{-})^2-1}{q^2(u_{+}/v_{-})^2-q^{-2}}\ell_{n,n+1}^{+\ts [1]}(u)\ell_{n+1,n}^{-\ts [1]}(v)
+\frac{(q^2-q^{-2})u_{+}/v_{-}}{q^2(u_{+}/v_{-})^2-q^{-2}}\ell_{n,n}^{+\ts [1]}(u)\ell_{n+1,n+1}^{-\ts [1]}(v)\\
&=\frac{(u_{-}/v_{+})^2-1}{q^2(u_{-}/v_{+})^2-q^{-2}}\ell_{n+1,n}^{-\ts [1]}(v)\ell_{n,n+1}^{+\ts [1]}(u)
+\frac{(q^2-q^{-2})u_{-}/v_{+}}{q^2(u_{-}/v_{+})^2-q^{-2}}\ell_{n+1,n+1}^{-\ts [1]}(v)\ell_{n,n}^{+\ts [1]}(u)
\eal
\eeq
The Gauss decomposition of $\Lc^{\pm [1]}(u)$ implies that the right hand side of \eqref{eq:lnn+1ln+1n} can be written as
\ben
\bal
\f^{-}_{n+1,n}(v)&\left(\frac{(u_{-}/v_{+})^2-1}{q^2(u_{-}/v_{+})^2-q^{-2}}\ell_{n,n}^{-\ts [1]}(v)\ell_{n,n+1}^{+\ts [1]}(u)
+\frac{(q^2-q^{-2})u_{-}/v_{+}}{q^2(u_{-}/v_{+})^2-q^{-2}}\ell_{n,n+1}^{-\ts [1]}(v)\ell_{n,n}^{+\ts [1]}(u)\right)\\
&+\frac{(q^2-q^{-2})u_{-}/v_{+}}{q^2(u_{-}/v_{+})^2-q^{-2}}\h_{n+1}^{-}(v)\h_{n}^{+}(u).
\eal
\een
Note that
\ben
\ell^{+\ts [1]}_{n,n+1}(u)\ell^{-\ts [1]}_{n,n}(v)=\left(\frac{(u_{-}/v_{+})^2-1}{q^2(u_{-}/v_{+})^2-q^{-2}}\ell_{n,n}^{-\ts [1]}(v)\ell_{n,n+1}^{+\ts [1]}(u)
+\frac{(q^2-q^{-2})u_{-}/v_{+}}{q^2(u_{-}/v_{+})^2-q^{-2}}\ell_{n,n+1}^{-\ts [1]}(v)\ell_{n,n}^{+\ts [1]}(u)\right),
\een
so the right hand side of \eqref{eq:lnn+1ln+1n} equals to
\ben
\f^{-}_{n+1,n}(v)\h^{+}_{n}(u)\e_{n,n+1}^{+}(u)\h^{-}_{n}(v)
+\frac{(q^2-q^{-2})u_{-}/v_{+}}{q^2(u_{-}/v_{+})^2-q^{-2}}\h_{n+1}^{-}(v)\h_{n}^{+}(u),
\een
and the relations between $f^{-}_{n+1,n}(v)$ and $h_{n}^{+}(u)$
\ben
\f^{-}_{n+1,n}(v)\h^{+}_{n}(u)
=\frac{(u_{+}/v_{-})^2-1}{q^2(u_{+}/v_{-})^2-q^{-2}}\h^{+}_{n}(u)\f^{-}_{n+1,n}(v)
+\frac{(q^2-q^{-2})u_{+}/v_{-}}{q^2(u_{+}/v_{-})^2-q^{-2}}\f_{n+1,n}^{+}(u)\h_{n}^{+}(u),
\een
brings the right hand side of \eqref{eq:lnn+1ln+1n} to
\ben
\bal
&\frac{(u_{+}/v_{-})^2-1}{q^2(u_{+}/v_{-})^2-q^{-2}}\h^{+}_{n}(u)\f^{-}_{n+1,n}(v)\e_{n,n+1}^{+}(u)\h_{n}^{-}(v)
+\frac{(q^2-q^{-2})u_{+}/v_{-}}{q^2(u_{+}/v_{-})^2-q^{-2}}\f_{n+1,n}^{+}(u)\h_{n}^{+}(u)\e_{n,n+1}^{+}(u)\h_{n}^{-}(v)\\
&+\frac{(q^2-q^{-2})u_{-}/v_{+}}{q^2(u_{-}/v_{+})^2-q^{-2}}\h_{n+1}^{-}(v)\h_{n}^{+}(u).
\eal
\een
Now comparing both sides of \eqref{eq:lnn+1ln+1n}, we get
\ben
\bal
\frac{(u_{+}/v_{-})^2-1}{q^2(u_{+}/v_{-})^2-q^{-2}}&
\h^{+}_{n}(u)\left[\e_{n,n+1}^{+}(u),\f^{-}_{n+1,n}(v)\right]\h_{n}^{-}(v)\\
&=
\frac{(q^2-q^{-2})u_{-}/v_{+}}{q^2(u_{-}/v_{+})^2-q^{-2}}\h_{n+1}^{-}(v)\h_{n}^{+}(u)
-\frac{(q^2-q^{-2})u_{+}/v_{-}}{q^2(u_{+}/v_{-})^2-q^{-2}}\h_{n+1}^{+}(u)\h_{n}^{-}(v)
\eal
\een
Using \eqref{eq:hnhn+1pm} and the invertibility of $h_{n}^{\pm}(u)$, we have
\ben
\left[\e_{n,n+1}^{+}(u),\f^{-}_{n+1,n}(v)\right]=
\frac{(q^2-q^{-2})u_{-}/v_{+}}{(u_{-}/v_{+})^2-1}\h_{n+1}^{-}(v)\h_{n}^{-}(v)^{-1}
-\frac{(q^2-q^{-2})u_{+}/v_{-}}{(u_{+}/v_{-})^2-1}\h_{n+1}^{+}(u)\h_{n}^{+}(u)^{-1}.
\een
Similarly, we can have
\ben
\left[\e_{n,n+1}^{\pm}(u),\f^{\mp}_{n+1,n}(v)\right]=
\frac{(q^2-q^{-2})u_{\mp}/v_{\pm}}{(u_{\mp}/v_{\pm})^2-1}\h_{n+1}^{\mp}(v)\h_{n}^{\mp}(v)^{-1}
-\frac{(q^2-q^{-2})u_{\pm}/v_{\mp}}{(u_{\pm}/v_{\mp})^2-1}\h_{n+1}^{\pm}(u)\h_{n}^{\pm}(u)^{-1},
\een
\ben
\left[\e_{n,n+1}^{\pm}(u),\f_{n+1,n}^{\pm}(v)\right]=\frac{(q^2-q^{-2})u/v}{(u/v)^2-1}\left(\h_{n+1}^{\pm}(v)\h_{n}^{\pm}(v)^{-1}
-\h_{n+1}^{\pm}(u)\h_{n}^{\pm}(u)^{-1}\right).
\een
Then using the definition of $\Xc_n^{\pm}(u)$, we have proved \eqref{eq:xnxnpm}.
\epf

\subsection{Formulas of $\z^{\pm}(z)$}

Write the relation \eqref{DLDL} in the form
\beql{DLbarD}
D{\Lc}^{\pm}(u\xi)^{t}D^{-1}={\Lc}^{\pm}(u)^{-1}{\z}^{[n] \pm}(u).
\eeq
Using the Gauss decomposition of ${\Lc}^{\pm}(z)$ and taking the $(2n,2n)$-th
entry on both sides of \eqref{DLbarD}, we have
\beql{h1h1pr}
{\h}_{1}^{\pm}(u\xi)={\h}_{1'}^{\pm}(u)^{-1}{\z}^{[n] \pm}(u).
\eeq

\bpr\label{prop:eiprei}
In $U(\overline{R}^{[n]})$, we have the following relations hold for $i=1,\dots ,n$:
\beql{ei'ei}
{\e}^{\pm}_{(i+1)'\ts i'}(u)=-{\e}_{i,i+1}^{\pm}(u\xi q^{2i}),
\eeq
\beql{fi'fi}
{\f}^{\pm}_{i'\ts (i+1)'}(u)=-{\f}_{i+1,i}^{\pm}(u\xi q^{2i}).
\eeq

\epr

\bpf
For $1\leq i\leq n-1$, by Proposition \ref{prop:gauss-consist} and \ref{prop:central}
we have
\beql{Lbar[n-i+1] z}
{\Lc}^{[n-i+1]\pm}(u)^{-1}\z^{[n-i+1]\pm}(u)=D^{[n-i+1]}{\Lc}^{[n-i+1]\pm}(u\xi q^{2i-2})'(D^{[n-i+1]})^{-1},
\eeq
where
\ben
D^{[n-i+1]}=\diag (q^{\bar{i}},...,q^{\bar{i'}}).
\een
Taking the $(i',i')$-th element on both sides of \eqref{Lbar[n-i+1] z}, we have
\beql{hihipr}
{\h}_i^{\pm}(u\xi q^{2i-2})={\h}_{i'}^{\pm}(u)^{-1}{\z}^{[n-i+1] \pm}(u).
\eeq
By taking the $((i+1)',i')$-th entry of \eqref{Lbar[n-i+1] z}, we obtain
\ben
-\e^{\pm}_{(i+1)',i'}(u)\h^{\pm}_{i'}(u)^{-1}{z}^{[n-i+1] \pm}(u)
=q{\h}_i^{\pm}(u\xi q^{2i-2}){\e}_{i,i+1}^{\pm}(u\xi q^{2i-2})
\een
Using \eqref{hihipr}, we have
\beql{ei'hi}
-\e^{\pm}_{(i+1)',i'}(u){\h}_i^{\pm}(u\xi q^{2i-2})
=q{\h}_i^{\pm}(u\xi q^{2i-2}){e}_{i,i+1}^{\pm}(u\xi q^{2i-2}).
\eeq
Furthermore, we have
\ben
q{\h}_i^{\pm}(u){\e}_{i,i+1}^{\pm}(u)={\e}_{i,i+1}^{\pm}(uq^2){\h}_i^{\pm}(u),
\een
which is the result of \cite{df:it}.

Thus, the equation \eqref{ei'hi} is equivalent to
\beql
-\e^{\pm}_{(i+1)',i'}(u){\h}_i^{\pm}(u\xi q^{2i-2})
={\e}_{i,i+1}^{\pm}(u\xi q^{2i}){\h}_i^{\pm}(u\xi q^{2i-2}).
\eeq
The invertibility of ${\h}_i^{\pm}(u\xi q^{2i-2})$ proves \eqref{ei'ei}.
Other relations are proved similarly.
\epf

\bpr
In $U(\overline{R}^{[n]})$, we have that
\beql{en'en}
{\e}^{\pm}_{(n+1)'\ts n'}(u)=-{\e}_{n,n+1}^{\pm}(-u),
\eeq
\beql{fn'fn}
{\f}^{\pm}_{n'\ts (n+1)'}(u)=-{\f}_{n+1,n}^{\pm}(-u),
\eeq
Thus, we have
\beq
\Xc_{n,2m+1}^{\pm}=0.
\eeq
\epr
Now we are in a position to give an explicit formula of $\bar{z}^{[n]\pm}(z)$ in terms of
$\bar{h}^{\pm}_{i}(z)$.

\bpr
In $U(\overline{R}^{[n]})$, we have
\beq
{\z}^{[n]\pm}(u)=\prod_{i=1}^{n-1}{\h}^{\pm}_{i}(u\xi q^{2i})^{-1}
\prod_{i=1}^{n}{\h}^{\pm}_{i}(u\xi q^{2i-2}){\h}^{\pm}_{n+1}(u).
\eeq
\epr

\bpf
Taking the $(2',2')$-th entry on both sides of \eqref{Lbar[n-i+1] z} and expressing the
entries of
${\Lc}^{[n]\pm}(u)^{-1}$ and ${\Lc}^{[n]\pm}(u\xi)^{'}$ in terms of the Gauss generators,
we have
\beq
{\h}^{\pm}_2(u\xi)+{\f}^{\pm}_{21}(u\xi){\h}^{\pm}_{1}(u\xi){\e}^{\pm}_{12}(u\xi)
=\big({\h}^{\pm}_{2'}(u)^{-1}
+{\e}^{\pm}_{2',1'}(u){\h}^{\pm}_{1'}(u)^{-1}{\f}^{\pm}_{1',2'}(u)\big){\z}^{[n]\pm}(u).
\eeq
Since ${\z}^{[n]\pm}(u)$ are central in the subalgebras generated by ${\Lc}^{[n]\pm}(z)$, using
\eqref{h1h1pr}, we can rewrite the above equation as
\ben
{\h}^{\pm}_{2'}(u)^{-1}{\z}^{[n]\pm}(u)=
{\h}^{\pm}_2(u\xi)+{\f}^{\pm}_{21}(u\xi){\h}^{\pm}_{1}(u\xi){\e}^{\pm}_{12}(u\xi)
-{\e}^{\pm}_{2',1'}(u){\h}^{\pm}_{1}(u\xi){\f}^{\pm}_{1',2'}(u)
\een
Now apply Proposition \ref{prop:eiprei} to obtain
\ben
{\h}^{\pm}_{2'}(u)^{-1}{\z}^{[n]\pm}(u)=
{\h}^{\pm}_2(u\xi)+{\f}^{\pm}_{21}(u\xi){\h}^{\pm}_{1}(u\xi){\e}^{\pm}_{12}(u\xi)
-{\e}^{\pm}_{12}(u\xi q^2){\h}^{\pm}_{1}(u\xi){\f}^{\pm}_{21}(u\xi q^2).
\een
Since
\ben
{\h}^{\pm}_{1}(u){\e}^{\pm}_{12}(u)=q^{-1}{\e}^{\pm}_{12}(uq^2){\h}^{\pm}_{1}(u),
~~~
{\h}^{\pm}_{1}(u){\f}^{\pm}_{21}(uq^2)=q^{-1}{\f}^{\pm}_{21}(u){\h}^{\pm}_{1}(u)
\een
by the results of \cite{df:it}, we have
\ben
{\h}^{\pm}_{2'}(u)^{-1}{\z}^{[n]\pm}(u)=
{\h}^{\pm}_2(u\xi)-q^{-1}[{\e}^{\pm}_{12}(u\xi q^2),{\f}^{\pm}_{21}(u\xi)]\bar{\h}^{\pm}_{1}(u\xi).
\een
Again using the results of \cite{df:it}
\ben
[{\e}^{\pm}_{12}(u),{\f}^{\pm}_{21}(v)]=
\frac{u(q-q^{-1})}{u-v}({\h}^{\pm}_{2}(v){\h}^{\pm}_{1}(v)^{-1}
-{\h}^{\pm}_{2}(u){\h}^{\pm}_{1}(u)^{-1}),
\een
we get
\ben
{\h}^{\pm}_{2'}(u)^{-1}{\z}^{[n]\pm}(u)=
{\h}^{\pm}_2(u\xi q^2){\h}^{\pm}_1(u\xi q^2)^{-1}{\h}^{\pm}_{1}(u\xi).
\een
Since ${\z}^{[n-1]\pm}(u)={\h}^{\pm}_{2'}(u){\h}^{\pm}_2(u\xi q^2)$,
we get a recurrence relation
\ben
{\z}^{[n]\pm}(u)=
{\h}^{\pm}_1(u\xi q^2)^{-1}{\h}^{\pm}_{1}(u\xi){\z}^{[n-1]\pm}(u).
\een
Then inductively we prove the proposition.
\epf
Similarly, we can get the same formula for $\z^{\pm}(u)\in U(R^{[n]})$.
\bpr
In $U({R}^{[n]})$, we have
\beq
{z}^{[n]\pm}(u)=\prod_{i=1}^{n-1}{h}^{\pm}_{i}(u\xi q^{2i})^{-1}
\prod_{i=1}^{n}{h}^{\pm}_{i}(u\xi q^{2i-2}){h}^{\pm}_{n+1}(u).
\eeq
\epr

\subsection{Other relations}

\bpr\label{Prophn+1X}
In $U(\overline{R}^{[n]})$, we have
\beql{hihn+1D}
{\h}^{\pm}_i(u){\h}^{\pm}_{n+1}(v)={\h}^{\pm}_{n+1}(v){\h}^{\pm}_i(u), ~~~~i=1,2,...,n,
\eeq

\beql{eq:hihn+1pm}
\frac{u_{\pm}/v_{\mp}-1}{qu_{\pm}/v_{\mp}-q^{-1}}{\h}^{\pm}_i(u){\h}^{\mp}_{n+1}(v)
=\frac{u_{\mp}/v_{\pm}-1}{qu_{\mp}/v_{\pm}-q^{-1}}{\h}^{\pm}_{n+1}(v){\h}^{\pm}_i(u), ~~~~i=1,2,...,n-1,
\eeq

\beql{eq:hn+1xjp}
{\h}_{n+1}^{\pm}(u){\Xc}_j^{+}(v)={\Xc}_j^{+}(v){\h}_{n+1}^{\pm}(u),~~~j\leq n-2,
\eeq

\beq
{\h}_{n+1}^{\pm}(u){\Xc}_j^{-}(v)={\Xc}_j^{-}(v){\h}_{n+1}^{\pm}(u),~~~j\leq n-2.
\eeq

\beql{eq:hn+1xn-1p}
{\h}_{n+1}^{\pm}(u)^{-1}{\Xc}_{n-1}^{+}(v){\h}_{n+1}^{\pm}(u)
=\frac{u/v_{\pm}+1}{q^{-1}u/v_{\pm}+q}{\Xc}_{n-1}^{+}(v)
\eeq

\beql{eq:hn+1xn-1m}
{\h}_{n+1}^{\pm}(u){\Xc}_{n-1}^{-}(v){\h}_{n+1}^{\pm}(u)^{-1}
=\frac{u/v_{\mp}+1}{q^{-1}u/v_{\mp}+q}{\Xc}_{n-1}^{-}(v)
\eeq

\beql{eq:hixnp}
{\h}_i^{\pm}(u){\Xc}_n^{+}(v)={\Xc}_n^{+}(v){\h}_{i}^{\pm}(u),~~~i\leq n-1
\eeq

\beql{eq:hixnm}
{\h}_{i}^{\pm}(u){\Xc}_n^{-}(v)={\Xc}_n^{+}(v){\h}_{i}^{\pm}(u),~~~j\leq n-1,
\eeq

\beql{eq:xixn}
{\Xc}_i^{\pm}(u){\Xc}^{\pm}_n(v)={\Xc}^{\pm}_n(v){\Xc}_i^{\pm}(u),~~~i< n-1,
\eeq
\epr

\bpf We only check \eqref{eq:hihn+1pm}-\eqref{eq:hixnp} and \eqref{eq:xixn}, as the others are trivial or can be treated similarly.

It follows from Corollary \ref{cor:commu} that
\beql{hiln+1C}
\begin{aligned}
&\frac{u_{\pm}/v_{\mp}-1}{qu_{\pm}/v_{\mp}-q^{-1}}{\h}^{\pm}_i(u)\big({\h}^{\mp}_{n+1}(v)
+{\f}_{n+1,n}^{\mp}(v){\h}_{n}^{\mp}(v){\e}_{n,n+1}^{\mp}(v)\big)=\\
&\frac{u_{\mp}/v_{\pm}-1}{qu_{\mp}/v_{\pm}-q^{-1}}\big({\h}^{\pm}_{n+1}(v)
+{\f}_{n+1,n}^{\mp}(v){\h}_{n}^{\mp}(v){e}_{n,n+1}^{\mp}(v)\big){\h}^{\pm}_i(u),
\end{aligned}
\eeq
and
\ben
\begin{aligned}
\frac{u_{\pm}/v_{\mp}-1}{qu_{\pm}/v_{\mp}-q^{-1}}{\h}^{\pm}_i(u){\f}_{n+1,n}^{\mp}(v){\h}_{n}^{\mp}(v)=
\frac{u_{\mp}/v_{\pm}-1}{qu_{\mp}/v_{\pm}-q^{-1}}{\f}_{n+1,n}^{\mp}(v){\h}_{n}^{\mp}(v){\h}^{\pm}_i(u),
\end{aligned}
\een
for $i=1,\dots n-1$.
Thus, the left hand side of \eqref{hiln+1C} equals to
\ben
 \frac{u_{\pm}/v_{\mp}-1}{qu_{\pm}/v_{\mp}-q^{-1}}{\h}^{\pm}_i(u){\h}^{\mp}_{n+1}(v)
+\frac{u_{\mp}/v_{\pm}-1}{qu_{\mp}/v_{\pm}-q^{-1}}{\f}_{n+1,n}^{\mp}(v){\h}_{n}^{\mp}(v){\h}^{\pm}_i(u){\e}_{n,n+1}^{\mp}(v).
\een
Again by Corollary \ref{cor:commu} it follows that
\ben
\begin{aligned}
\frac{u_{\pm}/v_{\mp}-1}{qu_{\pm}/v_{\mp}-q^{-1}}{\h}^{\pm}_i(u){\h}_{n}^{\mp}(v)=
\frac{u_{\mp}/v_{\pm}-1}{qu_{\mp}/v_{\pm}-q^{-1}}{\h}_{n}^{\mp}(v){\h}^{\pm}_i(u), ~~~~i=1,2,...,n-1.
\end{aligned}
\een
Then the left hand side of \eqref{hiln+1C} can be written as
\ben
 \frac{u_{\pm}/v_{\mp}-1}{qu_{\pm}/v_{\mp}-q^{-1}}{\h}^{\pm}_i(u){\h}^{\mp}_{n+1}(v)
+\frac{u_{\pm}/v_{\mp}-1}{qu_{\pm}/v_{\mp}-q^{-1}}{\f}_{n+1,n}^{\mp}(v){\h}^{\pm}_i(u)
{\h}_{n}^{\mp}(v){\e}_{n,n+1}^{\mp}(v))
\een
Also by Corollary \ref{cor:commu}, we have that
\ben
\begin{aligned}
\frac{u_{\pm}/v_{\mp}-1}{qu_{\pm}/v_{\mp}-q^{-1}}{\h}^{\pm}_i(u){\h}_{n}^{\mp}(v){\e}_{n,n+1}^{\mp}(v)=
\frac{u_{\mp}/v_{\pm}-1}{qu_{\mp}/v_{\pm}-q^{-1}}{\h}_{n}^{\mp}(v){\e}_{n,n+1}^{\mp}(v){\h}^{\pm}_i(u), ~~~~i=1,2,...,n-1,
\end{aligned}
\een
Finally, we get the left hand side of \eqref{hiln+1C} as the following form:
\ben
 \frac{u_{\pm}/v_{\mp}-1}{qu_{\pm}/v_{\mp}-q^{-1}}{\h}^{\pm}_i(u){\h}^{\mp}_{n+1}(v)+
\frac{u_{\mp}/v_{\pm}-1}{qu_{\mp}/v_{\pm}-q^{-1}}{\f}_{n+1,n}^{\mp}(v){\h}_{n}^{\mp}(v){\e}_{n,n+1}^{\mp}(v)){\h}^{\pm}_i(u)
\een
which then proves \eqref{eq:hihn+1pm}.

It follows from Remark \ref{re:anotherA} and \cite{df:it} that
\beq
{\h}_{n'}^{\pm}(u){\e}_{i',(i-1)'}^{\pm}(v)
={\e}_{i',(i-1)'}^{\pm}(v){\h}_{n'}^{\pm}(u),~~~
{\h}_{n'}^{\pm}(u){\e}_{i',(i-1)'}^{\mp}(v)
={\e}_{i',(i-1)'}^{\mp}(v){\h}_{n'}^{\pm}(u),
\eeq

\beq
{\h}_{n'}^{\pm}(u)^{-1}{\e}_{n',(n-1)'}^{\pm}(v){\h}_{n'}^{\pm}(u)
=\frac{qu/v-q^{-1}}{u/v-1}{\e}_{n',(n-1)'}^{\pm}(v)
-\frac{q-q^{-1}}{u/v-1}{\e}_{n',(n-1)'}^{\pm}(u),
\eeq

\beq
{\h}_{n'}^{\pm}(u)^{-1}{\e}_{n',(n-1)'}^{\mp}(v){\h}_{n'}^{\pm}(u)
=\frac{qu_{\mp}/v_{\pm}-q^{-1}}{u_{\mp}/v_{\pm}-1}{\e}_{n',(n-1)'}^{\mp}(v)
-\frac{q-q^{-1}}{u_{\mp}/v_{\pm}-1}{\e}_{n',(n-1)'}^{\pm}(u).
\eeq
They can be rewritten as follows (using Prop. \ref{prop:eiprei}): 
\beq
{\h}_{n'}^{\pm}(u){\e}_{i-1,i}^{\pm}(v)
={\e}_{i-1,i}^{\pm}(v){\h}_{n'}^{\pm}(u),~~~
{\h}_{n'}^{\pm}(u){\e}_{i-1,i}^{\mp}(v)
={\e}_{i-1,i}^{\mp}(v){\h}_{n'}^{\pm}(u),
\eeq
\beq
\bal
{\h}_{n'}^{\pm}(u)^{-1}{\e}_{n-1,n}^{\pm}(-vq^{-2}){\h}_{n'}^{\pm}(u)
&=\frac{qu/v-q^{-1}}{u/v-1}{\e}_{n-1,n}^{\pm}(-vq^{-2})\\
&-\frac{q-q^{-1}}{u/v-1}{\e}_{n-1,n}^{\pm}(-uq^{-2}),
\eal
\eeq

\beq
\bal
{\h}_{n'}^{\pm}(u)^{-1}{\e}_{n-1,n}^{\mp}(-vq^{-2}){\h}_{n'}^{\pm}(u)
&=\frac{qu_{\mp}/v_{\pm}-q^{-1}}{u_{\mp}/v_{\pm}-1}{\e}_{n-1,n}^{\mp}(-vq^{-2})\\
&-\frac{q-q^{-1}}{u_{\mp}/v_{\pm}-1}{\e}_{n-1,n}^{\pm}(-uq^{-2}),
\eal
\eeq
then \eqref{eq:hn+1xjp}-\eqref{eq:hn+1xn-1m} follow by using the definitions of $\Xc^{\pm}(u)$.

It follows from Corollary \ref{cor:commu} that for $i<n$
\beq
\frac{u_{\pm}/v_{\mp}-1}{qu_{\pm}/v_{\mp}-q^{-1}}{\h}_i^{\pm}(u){\h}_n^{\mp}(v){\e}_{n,n+1}^{\mp}(v)
=\frac{u_{\mp}/v_{\pm}-1}{qu_{\mp}/v_{\pm}-q^{-1}}{\h}_n^{\mp}(v){\e}_{n,n+1}^{\mp}(v){\h}_i^{\pm}(u).
\eeq
Furthermore, Corollary \ref{cor:commu} also gives that
\beq
\frac{u_{\pm}/v_{\mp}-1}{qu_{\pm}/v_{\mp}-q^{-1}}{\h}_i^{\pm}(u){\h}_n^{\mp}(v)
=\frac{u_{\mp}/v_{\pm}-1}{qu_{\mp}/v_{\pm}-q^{-1}}{\h}_n^{\mp}(v){\h}_i^{\pm}(u).
\eeq
Thus, we get
\beql{eq:hienpm}
{\h}_i^{\pm}(u){\e}_{n,n+1}^{\mp}(v)
={\e}_{n,n+1}^{\mp}(v){\h}_i^{\pm}(u).
\eeq
Similarly, we have
${\h}_i^{\pm}(u){\e}_{n,n+1}^{\pm}(v)
={\e}_{n,n+1}^{\pm}(v){\h}_i^{\pm}(u)$.
Then by the definition of $\Xc^{+}(u)$ we can prove \eqref{eq:hixnp}.

Next by Corollary \ref{cor:commu} it follows that for $i<n-1$
\beql{lii+1lnn+1}
\frac{u_{\pm}/v_{\mp}-1}{qu_{\pm}/v_{\mp}-q^{-1}}{\h}_i^{\pm}(u)
{\e}_{i,i+1}^{\pm}(u){\h}_{n}^{\mp}(v){\e}_{n,n+1}^{\mp}(v)
=\frac{u_{\mp}/v_{\pm}-1}{qu_{\mp}/v_{\pm}-q^{-1}}{\h}_{n}^{\mp}v()
{\e}_{n,n+1}^{\mp}(v){\h}_i^{\pm}(u){\e}_{i,i+1}^{\pm}(u),
\eeq
\beq
\frac{u_{\pm}/v_{\mp}-1}{qu_{\pm}/v_{\mp}-q^{-1}}{\h}_i^{\pm}(u){\e}_{i,i+1}^{\pm}(u){\h}_{n}^{\mp}(v)
=\frac{u_{\mp}/v_{\pm}-1}{qu_{\mp}/v_{\pm}-q^{-1}}{\h}_{n}^{\mp}(v){\h}_i^{\pm}(u){\e}_{i,i+1}^{\pm}(u).
\eeq
Then \eqref{lii+1lnn+1} can be written as
\ben
{\h}_{n}^{\mp}(v){\h}_i^{\pm}(u){\e}_{i,i+1}^{\pm}(u){\e}_{n,n+1}^{\mp}(v)
={\h}_{n}^{\mp}(v){\e}_{n,n+1}^{\mp}(v){\h}_i^{\pm}(u){\e}_{i,i+1}^{\pm}(u).
\een
Using the above equation and the fact $\left[{\h}_i^{\pm}(u),{\e}_{n,n+1}^{\mp}(v)\right]=0$ (see \cite{df:it}), we have
\ben
{\e}_{i,i+1}^{\pm}(u){\e}_{n,n+1}^{\mp}(v)
={\e}_{n,n+1}^{\mp}(v){\e}_{i,i+1}^{\pm}(u).
\een
Similarly, we have
\ben
{\e}_{i,i+1}^{\pm}(u){\e}_{n,n+1}^{\pm}(v)
={\e}_{n,n+1}^{\pm}(v){\e}_{i,i+1}^{\pm}(u).
\een
Using these, we get $\Xc_i^{+}(u)\Xc_n^{+}(v)=\Xc_n^{+}(v)\Xc_i^{+}(u)$ easily. Similarly, we
can get $\Xc_i^{-}(u)\Xc_n^{-}(v)=\Xc_n^{-}(v)\Xc_i^{-}(u)$.
\epf

\bpr\label{prop:xn-1xn}
In $U(\overline{R}^{[n]})$, we have the following relation hold:
\beql{eq:Xn-1Xnp}
(u^2q^{2}-v^2){\Xc}_{n-1}^{+}(uq^{n-1}){\Xc}^{+}_n(vq^n)=(u^2-q^2v^2){\Xc}^{+}_n(vq^n){\Xc}_{n-1}^{+}(uq^{n-1}).
\eeq
\beql{eq:Xn-1Xnm}
(u^2q^{2}-v^2)^{-1}{\Xc}_{n-1}^{-}(uq^{n-1}){\Xc}^{-}_n(vq^n)=(u^2-q^2v^2)^{-1}{\Xc}^{-}_n(vq^n){\Xc}_{n-1}^{-}(uq^{n-1}).
\eeq
\epr

\bpf
Here we only prove \eqref{eq:Xn-1Xnp}, as \eqref{eq:Xn-1Xnm} can be treated similarly.

It follows from Corollary \ref{cor:guass-embed} that
\beql{l12l23C}
\begin{aligned}
&\frac{u_{\pm}/v_{\mp}-1}{qu_{\pm}/v_{\mp}-q^{-1}}{\ell}_{n-1,n}^{\pm\ts [2]}(u){\ell}_{n,n+1}^{\mp\ts [2]}(v)
+\frac{(q-q^{-1})u_{\pm}/v_{\mp}}{qu_{\pm}/v_{\mp}-q^{-1}}{\ell}_{n,n}^{\pm\ts [2]}(u){\ell}_{n-1,n+1}^{\mp\ts [2]}(v)=\\
&\frac{(q-q^{-1})q^{-3}(u_{\mp}/v_{\pm}-1)}{(qu_{\mp}/v_{\pm}-q^{-1})(u_{\mp}/v_{\pm}+q^{-4})}
{\ell}_{n,n-1}^{\mp\ts [2]}(v){\ell}_{n-1,n+2}^{\pm\ts [2]}(u)\\
&+\frac{(q-q^{-1})(q^{-4}+1)u_{\mp}/v_{\pm}}{(qu_{\mp}/v_{\pm}-q^{-1})(u_{\mp}/v_{\pm}+q^{-4})}
{\ell}_{n,n}^{\mp\ts [2]}(v){\ell}_{n-1,n+1}^{\pm\ts [2]}(u)\\
&+\frac{q^{-2}(u_{\mp}/v_{\pm}-1)(qu_{\mp}/v_{\pm}+q^{-1})}{(qu_{\mp}/v_{\pm}-q^{-1})(u_{\mp}/v_{\pm}+q^{-4})}
{\ell}_{n,n+1}^{\mp\ts [2]}(v){\ell}_{n-1,n}^{\pm\ts [2]}(u)\\
&-\frac{q^{-1}(q-q^{-1})(u_{\mp}/v_{\pm}-1)u_{\mp}/v_{\pm}}{(qu_{\mp}/v_{\pm}-q^{-1})(u_{\mp}/v_{\pm}+q^{-4})}
{\ell}_{n,n+2}^{\mp\ts [2]}(v){\ell}_{n-1,n-1}^{\pm\ts [2]}(u)
\end{aligned}
\eeq
Plug the Gauss decomposition of $\ell_{n,n+1}^{\mp\ts [2]}(v)$ into the left-hand side of \eqref{l12l23C}: 
\beq
\bal
&\frac{u_{\pm}/v_{\mp}-1}{qu_{\pm}/v_{\mp}-q^{-1}}
{\ell}_{n-1,n}^{\pm\ts [2]}(u){\h}_{n}^{\mp}(v){\e}_{n,n+1}^{\mp}(v)+
\frac{u_{\pm}/v_{\mp}-1}{qu_{\pm}/v_{\mp}-q^{-1}}
{\ell}_{n-1,n}^{\pm \ts [2]}(u){\ell}_{n,n-1}^{\mp \ts [2]}(v)
{\e}_{n-1,n+1}^{\mp}(v)\\
&+\frac{(q-q^{-1})u_{\pm}/v_{\mp}}{qu_{\pm}/v_{\mp}-q^{-1}}
{\ell}_{n,n}^{\pm\ts [2]}(u){\ell}_{n-1,n-1}^{\mp \ts[2]}(v){\e}_{n-1,n+1}^{\mp}(v).
\eal
\eeq
 Recalling the the commuting relations between ${\ell}_{n-1,n}^{\pm [2]}(u)$ and ${\ell}_{n,n-1}^{\mp \ts [2]}(v)$:
\ben
\begin{aligned}
&\frac{u_{\pm}/v_{\mp}-1}{qu_{\pm}/v_{\mp}-q^{-1}}{\ell}_{n-1,n}^{\pm\ts [2]}(u){\ell}_{n,n-1}^{\mp\ts[2]}(v)
+\frac{(q-q^{-1})u_{\pm}/v_{\mp}}{qu_{\pm}/v_{\mp}-q^{-1}}{\ell}_{n,n}^{\pm\ts[2]}(u){\ell}_{n-1,n-1}^{\mp\ts [2]}(v)\\
&=\frac{u_{\mp}/v_{\pm}-1}{qu_{\mp}/v_{\pm}-q^{-1}}{\ell}_{n,n-1}^{\mp\ts [2]}(v){\ell}_{n-1,n}^{\pm\ts [2]}(u)
+\frac{(q-q^{-1})u_{\mp}/v_{\pm}}{qu_{\mp}/v_{\pm}-q^{-1}}{\ell}_{n,n}^{\mp\ts [2]}(v){\ell}_{n-1,n-1}^{\pm\ts [2]}(u),
\end{aligned}
\een
we see the left hand side of \eqref{l12l23C} is equal to
\ben
\bal
\frac{u_{\pm}/v_{\mp}-1}{qu_{\pm}/v_{\mp}-q^{-1}}
{\ell}_{n-1,n}^{\pm\ts [2]}(u){\h}_{n}^{\mp}(v){\e}_{n,n+1}^{\mp}(v)&+
\frac{u_{\mp}/v_{\pm}-1}{qu_{\mp}/v_{\pm}-q^{-1}}{\ell}_{n,n-1}^{\mp\ts [2]}(v){\ell}_{n-1,n}^{\pm\ts [2]}(u)
{\e}_{n-1,n+1}^{\mp}(v)\\
&+\frac{(q-q^{-1})u_{\mp}/v_{\pm}}{qu_{\mp}/v_{\pm}-q^{-1}}{\ell}_{n,n}^{\mp\ts [2]}(v){\ell}_{n-1,n-1}^{\pm\ts [2]}(u){\e}_{n-1,n+1}^{\mp}(v),
\eal
\een
which is equivalent to the following with the help of the Gauss decomposition of ${\Lc}^{\pm\ts [2]}(u)$:
\ben
\bal
&\frac{u_{\pm}/v_{\mp}-1}{qu_{\pm}/v_{\mp}-q^{-1}}
{\ell}_{n-1,n}^{\pm\ts [2]}(u){\h}_{n}^{\mp}(v){\e}_{n,n+1}^{\mp}(v)+
\frac{u_{\mp}/v_{\pm}-1}{qu_{\mp}/v_{\pm}-q^{-1}}{\f}_{n,n-1}^{\mp}(v){\ell}_{n-1,n-1}^{\mp\ts [2]}(v){\ell}_{n-1,n}^{\pm\ts [2]}(u){\e}_{n-1,n+1}^{\mp}(v)\\
&+\frac{(q-q^{-1})u_{\mp}/v_{\pm}}{qu_{\mp}/v_{\pm}-q^{-1}}{\h}_{n}^{\mp }(v){\h}_{n-1}^{\pm}(u){\e}_{n-1,n+1}^{\mp}(v)
+\frac{(q-q^{-1})u_{\mp}/v_{\pm}}{qu_{\mp}/v_{\pm}-q^{-1}}{\f}_{n,n-1}^{\mp}(v)
{\ell}_{n-1,n}^{\pm\ts[2]}(v){\ell}_{n-1,n-1}^{\pm\ts[2]}(u){\e}_{n-1,n+1}^{\mp}(v)
\eal
\een
Then using
the relations between ${\ell}_{n-1,n}^{\pm\ts [2]}(u)$ and
${\ell}_{n-11,n-1}^{\mp\ts[2]}(v)$:
\ben
{\ell}_{n-1,n}^{\pm\ts[2]}(u){\ell}_{n-1,n-1}^{\mp\ts [2]}(v)=
\frac{u_{\mp}-v_{\pm}}{qu_{\mp}-q^{-1}v_{\pm}}{\ell}_{n-1,n-1}^{\mp\ts [2]}(v){\ell}_{n-1,n}^{\pm\ts [2]}(u)
+\frac{(q-q^{-1})u_{\mp}/v_{\pm}}{qu_{\mp}/v_{\pm}-q^{-1}}{\ell}_{n-1,n}^{\mp\ts [2]}(v){\ell}_{n-1,n-1}^{\pm\ts[2]}(u),
\een
we can get the left hand side of \eqref{l12l23C} equals to
\ben
\bal
&\frac{u_{\pm}/v_{\mp}-1}{qu_{\pm}/v_{\mp}-q^{-1}}
{\ell}_{n-1,n}^{\pm\ts [2]}(u){\h}_{n}^{\mp}(v){\e}_{n,n+1}^{\mp}(v)+
\frac{(q-q^{-1})u_{\mp}/v_{\pm}}{qu_{\mp}/v_{\pm}-q^{-1}}{\h}_{n}^{\mp }(v){\h}_{n-1}^{\pm}(u){\e}_{n-1,n+1}^{\mp}(v)\\
&+{\f}_{n,n-1}^{\mp}(v){\ell}_{n-1,n}^{\pm\ts [2]}(u){\ell}_{n-1,n+1}^{\mp\ts [2]}(v).
\eal
\een

Next by the relations between ${\ell}_{n-1,n}^{\pm\ts [2]}(u)$ and ${\ell}_{n-1,n+1}^{\mp\ts [2]}(v)$:
\beq
\begin{aligned}
{\ell}_{n-1,n}^{\pm\ts [2]}(u){\ell}_{n-1,n+1}^{\mp\ts [2]}(v)&=
\frac{(q-q^{-1})q^{-3}(u_{\mp}/v_{\pm}-1)}{(qu_{\mp}/v_{\pm}-q^{-1})(u_{\mp}/v_{\pm}+q^{-4})}
{\ell}_{n-1,n-1}^{\mp\ts [2]}(v){\ell}_{n-1,n+2}^{\pm\ts [2]}(u)\\
&+\frac{(q-q^{-1})(q^{-4}+1)u_{\mp}/v_{\pm}}{(qu_{\mp}/v_{\pm}-q^{-1})(u_{\mp}/v_{\pm}+q^{-4})}
{\ell}_{n-1,n}^{\mp\ts [2]}(v){\ell}_{n-1,n+1}^{\pm\ts [2]}(u)\\
&+\frac{q^{-2}(u_{\mp}/v_{\pm}-1)(qu_{\mp}/v_{\pm}+q^{-1})}{(qu_{\mp}/v_{\pm}-q^{-1})(u_{\mp}/v_{\pm}+q^{-4})}
{\ell}_{n-1,n+1}^{\mp\ts [2]}(v){\ell}_{n-1,n}^{\pm\ts [2]}(u)\\
&-\frac{q^{-1}(q-q^{-1})(u_{\mp}/v_{\pm}-1)u_{\mp}/v_{\pm}}{(qu_{\mp}/v_{\pm}-q^{-1})(u_{\mp}/v_{\pm}+q^{-4})}
{\ell}_{n-1,n+2}^{\mp\ts [2]}(v){\ell}_{n-1,n-1}^{\pm\ts [2]}(u),
\end{aligned}
\eeq
the equation \eqref{l12l23C} is seen as the following: 

\ben
\begin{aligned}
&\frac{u_{\pm}/v_{\mp}-1}{qu_{\pm}/v_{\mp}-q^{-1}}
{\ell}_{n-1,n}^{\pm\ts [2]}(u){\h}_{n}^{\mp}(v){\e}_{n,n+1}^{\mp}(v)+
\frac{(q-q^{-1})u_{\mp}/v_{\pm}}{qu_{\mp}/v_{\pm}-q^{-1}}{\h}_{n}^{\mp }(v){\h}_{n-1}^{\pm}(u){\e}_{n-1,n+1}^{\mp}(v)\\
&=
\frac{(q-q^{-1})(q^{-4}+1)u_{\mp}/v_{\pm}}{(qu_{\mp}/v_{\pm}-q^{-1})(u_{\mp}/v_{\pm}+q^{-4})}
{\h}_{n}^{\mp}(v){\ell}_{n-1,n+1}^{\pm\ts [2]}(u)\\
&+\frac{q^{-2}(u_{\mp}/v_{\pm}-1)(qu_{\mp}/v_{\pm}+q^{-1})}{(qu_{\mp}/v_{\pm}-q^{-1})(u_{\mp}/v_{\pm}+q^{-4})}
{\h}_{n}^{\mp}(v){\e}_{n,n+1}^{\mp}(v){\ell}_{n-1,n}^{\pm\ts [2]}(u)\\
&-\frac{q^{-1}(q-q^{-1})(u_{\mp}/v_{\pm}-1)u_{\mp}/v_{\pm}}{(qu_{\mp}/v_{\pm}-q^{-1})(u_{\mp}/v_{\pm}+q^{-4})}
{\h}_{n}^{\mp}(v){\e}_{n,n+2}^{\mp}(v){\ell}_{n-1,n-1}^{\pm\ts [2]}(u),
\end{aligned}
\een
By \eqref{ELMPj=l}, we derive the following relations between ${\e}_{n-1,n}^{\pm}(u)$ and ${\h}_{n}^{\mp}(v)$:
\ben
{\e}_{n-1,n}^{\pm}(u){\h}_n^{\mp}(v)
=\frac{qu_{\mp}/v_{\pm}-q^{-1}}{u_{\mp}/v_{\pm}-1}{\h}_n^{\mp}(v){\e}_{n-1,n}^{\pm}(u)-
\frac{(q-q^{-1})u_{\mp}/v_{\pm}}{u_{\mp}/v_{\pm}-1}
{\h}_{n}^{\mp}(v){\e}_{n-1,n}^{\mp}(v)
\een
Furthermore by the relation between ${\h}_{n-1}^{\pm}(u)$ and ${\h}_n^{\mp}(v)$:
\ben
\frac{u_{\pm}/v_{\mp}-1}{qu_{\pm}/v_{\mp}-q^{-1}}\h^{\pm}_{n-1}(u)\h^{\mp}_n(v)=
\frac{u_{\mp}/v_{\pm}-1}{qu_{\mp}/v_{\pm}-q^{-1}}\h^{\mp}_n(v)\h^{\pm}_{n-1}(u),
\een
we have
\ben
\frac{u_{\pm}/v_{\mp}-1}{qu_{\pm}/v_{\mp}-q^{-1}}
{\ell}_{n-1,n}^{\pm\ts [2]}(u){\h}_{n}^{\mp}(v)={\h}_{n}^{\mp}(v)\h_{n-1}^{\pm}(u)\Big({\e}_{n-1,n}^{\pm}(u)
-\frac{(q-q^{-1})u_{\mp}/v_{\pm}}{qu_{\mp}/v_{\pm}-q^{-1}}\e_{n-1,n}^{\mp}(v)\Big).
\een

Cancelling ${\h}_{n-1}^{\pm}(u)$ and ${\h}_{n}^{\mp}(v)$ and noticing $[{\h}_{n-1}^{\pm}(u),{\e}_{n,n+1}^{\mp}(v)]=0$
(see \eqref{eq:hienpm}), we have that
\eqref{l12l23C} is equivalent to
\beql{e12e23C}
\begin{aligned}
&{\e}_{n-1,n}^{\pm}(u){\e}_{n,n+1}^{\mp}(v)
+\frac{(q-q^{-1})u_{\mp}/v_{\pm}}{qu_{\mp}/v_{\pm}-q^{-1}}\Big({\e}_{n-1,n+1}^{\mp}(v)
-{\e}_{n-1,n}^{\mp}(v){\e}_{n,n+1}^{\mp}(v)\Big)=\\
&\frac{(q-q^{-1})(q^{-4}+1)u_{\mp}/v_{\pm}}{(qu_{\mp}/v_{\pm}-q^{-1})(u_{\mp}/v_{\pm}+q^{-4})}
{\e}_{n-1,n+1}^{\pm}(u)
+\frac{q^{-2}(u_{\mp}/v_{\pm}-1)(qu_{\mp}/v_{\pm}+q^{-1})}{(qu_{\mp}/v_{\pm}-q^{-1})(u_{\mp}/v_{\pm}+q^{-4})}
{\e}_{n,n+1}^{\mp}(v){\e}_{n-1,n}^{\pm}(u)\\
&-\frac{q^{-1}(q-q^{-1})(u_{\mp}/v_{\pm}-1)u_{\mp}/v_{\pm}}{(qu_{\mp}/v_{\pm}-q^{-1})(u_{\mp}/v_{\pm}+q^{-4})}
{\h}_{n-1}^{\pm}(u)^{-1}{\e}_{n,n+2}^{\mp}(v){\h}_{n-1}^{\pm}(u),
\end{aligned}
\eeq

By the similar process, we have the relation between ${\h}_{n-1}^{\pm}(u)$ and ${\e}_{n,n+2}^{\mp}(v)$:
\ben
\begin{aligned}
\frac{z_{\pm}-w_{\mp}}{qz_{\pm}-q^{-1}w_{\mp}}\bar{h}_{1}^{\pm}(z)\bar{h}_{2}^{\mp}(w)\bar{e}_{24}^{\mp}(w)&=
\frac{(z_{\mp}-w_{\pm})q^{-1}(z_{\mp}-w_{\pm}q^{-4})}{(qz_{\mp}-q^{-1}w_{\pm})(z_{\mp}-q^{-6}w_{\pm})}
\bar{h}_{2}^{\mp}(w)\bar{e}_{24}^{\mp}(w)\bar{h}_{1}^{\pm}(z)\\
&+\frac{(q-q^{-1})(z_{\mp}-w_{\pm})w_{\pm}q^{-3}}{(qz_{\mp}-q^{-1}w_{\pm})(z_{\mp}-q^{-6}w_{\pm})}
\bar{h}_{2}^{\mp}(w)\bar{h}_{1}^{\pm}(z)\bar{e}_{13}^{\pm}(z)\\
&-\frac{(q-q^{-1})(z_{\mp}-w_{\pm})w_{\pm}q^{-5}}{(qz_{\mp}-q^{-1}w_{\pm})(z_{\mp}-q^{-6}w_{\pm})}
\bar{h}_{2}^{\mp}(w)\bar{e}_{23}(w)\bar{h}_{1}^{\pm}(z)\bar{e}_{12}^{\pm}(z)\\
\end{aligned}
\een
Since we also have the following relation
\ben
\frac{z_{\pm}-w_{\mp}}{qz_{\pm}-q^{-1}w_{\mp}}\bar{h}_{1}^{\pm}(z)\bar{h}_{2}^{\mp}(w)
=\frac{z_{\mp}-w_{\pm}}{qz_{\mp}-q^{-1}w_{\pm}}\bar{h}_{2}^{\mp}(w)\bar{h}_{1}^{\pm}(z),
\een
we can get
\beql{h1e24C}
\begin{aligned}
{\e}_{n,n+2}^{\mp}(v)&=
\frac{q^{-2}(qu_{\mp}/v_{\pm}+q^{-1})}{u_{\mp}/v_{\pm}+q^{-4}}{\h}_{n-1}^{\pm}(u)^{-1}{\e}_{n,n+2}^{\mp}(v){\h}_{n-1}^{\pm}(u)\\
&+\frac{(q-q^{-1})q^{-3}}{u_{\mp}/v_{\pm}+q^{-4}}\Big({\e}_{n-1,n+1}^{\pm}(u)+{\e}_{n,n+1}^{\mp}(v){\e}_{n-1,n}^{\pm}(u)\Big).
\end{aligned}
\eeq
Applying \eqref{h1e24C} to \eqref{e12e23C}, we get
\begin{align}\nonumber
{\e}_{n-1,n}^{\pm}(u){\e}_{n,n+1}^{\mp}(v)&=\frac{(u_{\mp}/v_{\pm})^2-1}{q^2(u_{\mp}/v_{\pm})^2-q^{-2}}
\e_{n,n+1}^{\mp}(v)\e_{n-1,n}^{\pm}(u)
+\frac{(q^2-q^{-2})u_{\mp}/v_{\pm}}{q^2(u_{\mp}/v_{\pm})^2-q^{-2}}\e_{n-1,n+1}^{\pm}(u)\\
&-\frac{q(q-q^{-1})(u_{\mp}/v_{\pm}-1)u_{\mp}/v_{\pm}}{q^2(u_{\mp}/v_{\pm})^2-q^{-2}}\e_{n,n+2}^{\mp}(v)\\ \nonumber
&+\frac{(q-q^{-1})u_{\mp}/v_{\pm}}{qu_{\mp}/v_{\pm}-q^{-1}}\Big(\e_{n-1,n}^{\mp}(v)\e_{n,n+1}^{\mp}(v)-\e_{n-1,n+1}^{\mp}(v)\Big).
\end{align}
Similarly, we obtain that
\beq
\bal
{\e}_{n-1,n}^{\pm}(u){\e}_{n,n+1}^{\pm}(v)&=\frac{(u/v)^2-1}{q^2(u/v)^2-q^{-2}}
\e_{n,n+1}^{\pm}(v)\e_{n-1,n}^{\pm}(u)
+\frac{(q^2-q^{-2})u/v}{q^2(u/v)^2-q^{-2}}\e_{n-1,n+1}^{\pm}(u)\\
&-\frac{q(q-q^{-1})(u/v-1)u/v}{q^2(u/v)^2-q^{-2}}\e_{n,n+2}^{\pm}(v)\\
&+\frac{(q-q^{-1})u/v}{qu/v-q^{-1}}\Big(\e_{n-1,n}^{\pm}(v)\e_{n,n+1}^{\pm}(v)-\e_{n-1,n+1}^{\pm}(v)\Big).
\eal
\eeq
Then we can prove \eqref{eq:Xn-1Xnp}.
\epf

The following follows immediately from Proposition \ref{prop:xn-1xn} and \eqref{eq:xixj}. 
\bpr
In $U(\overline{R})$, we have
\beq
\bal
&\sum_{\sigma\in \mathfrak{S}_2}\sigma\Big((q^2u_1+u_2)(\Xc^{\pm}_n(v)\Xc_{n-1}^{\pm}(u_1)\Xc_{n-1}^{\pm}(u_2)\\
&-[2]_{q^2}\Xc_{n-1}^{\pm}(u_1)\Xc_{n}^{\pm}(v)\Xc_{n-1}^{\pm}(u_2)
+\Xc_{n-1}^{\pm}(u_1)\Xc_{n-1}^{\pm}(u_2)\Xc_{n}^{\pm}(v))\Big)=0
\eal
\eeq
\epr

\bpr
In $U(\overline{R}^{[n]})$, we have that
\beq
[{\Xc}_i^{+}(u),{\Xc}_n^{-}(v)]=0, \tss \text{for $i < n$}.
\eeq
\epr

\bpf
For $i\leq n-2$, it follows from Corollary \ref{cor:commu} that
\ben
\frac{u_{\pm}/v_{\mp}-1}{qu_{\pm}/v_{\mp}-q^{-1}}
{\h}_i^{\pm}(u){\e}_{i,i+1}^{\pm}(u){\f}_{n+1,n}^{\mp}(v){\h}_{n}^{\mp}(v)
=\frac{u_{\mp}/v_{\pm}-1}{qu_{\mp}/v_{\pm}-q^{-1}}
{\f}_{n+1,n}^{\mp}(v){\h}_{n}^{\mp}(v){\h}_i^{\pm}(z){\e}_{i,i+1}^{\pm}(u).
\een
Again by Corollary \ref{cor:commu} we have
\ben
\frac{u_{\pm}/v_{\mp}-1}{qu_{\pm}/v_{\mp}-q^{-1}}{\h}_i^{\pm}(u){\f}_{n+1,n}^{\mp}(v){\h}_{n}^{\mp}(v)
=\frac{u_{\mp}/v_{\pm}-1}{qu_{\mp}/v_{\pm}-q^{-1}}{\f}_{n+1,n}^{\mp}(v){\h}_{n}^{\mp}(v){\h}_i^{\pm}(u).
\een
Therefore
\ben
{\h}_i^{\pm}(u){\e}_{i,i+1}^{\pm}(u){\f}_{n+1,n}^{\mp}(v){\h}_{n}^{\mp}(v)
={\h}_i^{\pm}(u){\f}_{n+1,n}^{\mp}(v){\h}_{n}^{\mp}(v){\e}_{i,i+1}^{\pm}(u).
\een
Furthermore, since $\left[{\h}_{n}^{\mp}(v),{\e}_{i,i+1}^{\pm}(u)\right]=0$ (see result in \cite{df:it}), we can prove
\beq
\left[{\e}_{i,i+1}^{\pm}(u),{\f}_{n+1,n}^{\mp}(v)\right]=0,
\eeq
for $i=1,\dots,n-2$. Similarly, we have
\beq
\left[{\e}_{i,i+1}^{\pm}(u),{\f}_{n+1,n}^{\mp}(v)\right]=0,
\eeq
for $i=1,\dots,n-2$.
Thus, for $i=1,\dots,n-2$ we have $[{\Xc}_i^{+}(u),{\Xc}_n^{-}(v)]=0$.

For $i=n-1$, by \eqref{ELMPj=l} we have
\ben
\bal
{\e}_{n-1,n}^{\pm}(u){\f}_{n+1,n}^{\mp}(v){\h}_n^{\mp}(v)
&=\frac{qu_{\mp}/v_{\pm}-q^{-1}}{u_{\mp}/v_{\pm}-1}{\f}_{n+1,n}^{\mp}(v){\h}_{n}^{\mp}(v){\e}_{n-1,n}^{\pm}(u)\\
&-\frac{q-q^{-1}}{u_{\mp}/v_{\pm}-1}{\f}_{n+1,n}^{\mp}(v){\h}_{n}^{\mp}(v){\e}_{n-1,n}^{\mp}(v),
\eal
\een
and
\ben
{\e}_{n-1,n}^{\pm}(u){\h}_n^{\mp}(v)=
\frac{qu_{\mp}/v_{\pm}-q^{-1}}{u_{\mp}/v_{\pm}-1}{\h}_{n}^{\mp}(v){\e}_{n-1,n}^{\pm}(u)
-\frac{q-q^{-1}}{u_{\mp}/v_{\pm}-1}{\h}_{n}^{\mp}(v){\e}_{n-1,n}^{\mp}(v).
\een
Then, by the invertibility of ${\h}_{n}^{\mp}(v)$ we can prove
 $\left[{\e}_{n-1,n}^{\pm}(u),{\f}_{n+1,n}^{\mp}(v)\right]=0$.
Similarly, we have $\left[{\e}_{n-1,n}^{\pm}(u),{\f}_{n+1,n}^{\pm}(v)\right]=0$.
Thus we have $[{\Xc}_{n-1}^{+}(u),{\Xc}_n^{-}(v)]=0$.
\epf

\bth\label{thm:Ubar}
The following relations between the Gaussian generators hold in the algebra
$U(\overline{R}^{\tss[n]})$. For the relations involving $\h^{\pm}_i(u)$ we have
\begin{align} 
&\h^{\pm}_i(u)\h^{\pm}_j(v)=\h^{\pm}_j(v)\h^{\pm}_i(u), \qquad \h^{\pm}_i(u)\h^{\mp}_i(v)=\h^{\mp}_i(v)\h^{\pm}_i(u); \\ \label{hihjmp}
&\frac{u_{\pm}-v_{\mp}}{qu_{\pm}-q^{-1}v_{\mp}}\h^{\pm}_i(u)\h^{\mp}_j(v)=
\frac{u_{\mp}-v_{\pm}}{qu_{\mp}-q^{-1}v_{\pm}}\h^{\mp}_j(v)\h^{\pm}_i(u), \quad i<j, i\neq n;\\
&\frac{(u_{\pm}/v_{\mp})^{2}-1}{q^2(u_{\pm}/ v_{\mp})^2-q^{-2}}\h^{\pm}_n(u)\h_{n+1}^{\mp}(v)=
\frac{(u_{\mp}/v_{\pm})^{2}-1}{q^2(u_{\mp}/ v_{\pm})^2-q^{-2}}\h_{n+1}^{\mp}(v)\h^{\pm}_{n}(u).
\end{align}
The relations
involving $\h^{\pm}_i(u)$ and $\Xc_{j}^{\pm}(v)$ are
\ben
\bal
\h_{i}^{\pm}(u)\Xc_{j}^{+}(v)\h_{i}^{\pm}(u)^{-1}&
=\frac{u/v_{\pm}-1}{q^{(\ep_i,\alpha_j)}u/v_{\pm}-q^{-(\ep_i,\alpha_j)}}
\Xc_{j}^{+}(v),\\[0.4em]
\h_{i}^{\pm}(u)^{-1}\Xc_{j}^{-}(v)\h_{i}^{\pm}(u)&
=\frac{u/v_{\mp}-1}{q^{(\ep_i,\alpha_j)}u/v_{\mp}-q^{-(\ep_i,\alpha_j)}}
\Xc_{j}^{-}(v)
\eal
\een
for $i=1,\dots,n$, $j=1,\dots, n-1$
together with
\ben
\bal
\h_{i}^{\pm}(u)\Xc_{n}^{+}(v)\h_{i}^{\pm}(u)^{-1}&
=\frac{(u/v_{\pm})^2-1}{q^{(\ep_i,\alpha_j)}(u/v_{\pm})^2-q^{-(\ep_i,\alpha_j)}}
\Xc_{n}^{+}(v),\\[0.4em]
\h_{i}^{\pm}(u)^{-1}\Xc_{n}^{-}(v)\h_{i}^{\pm}(u)&
=\frac{(u/v_{\mp})^2-1}{q^{(\ep_i,\alpha_j)}(u/v_{\mp})^2-q^{-(\ep_i,\alpha_j)}}
\Xc_{n}^{-}(v)
\eal
\een
for $i=1,\dots n$
and
\ben
\bal
\h_{n+1}^{\pm}(u)\Xc_{n}^{+}(v)\h_{n+1}^{\pm}(u)^{-1}&
=\frac{(u/v_{\pm})^{2}-1}{q^{-2}(u/v_{\pm})^2-q^{2}}
\Xc_{n}^{+}(v),\\[0.4em]
\h_{n+1}^{\pm}(u)^{-1}\Xc_{n}^{-}(v)\h_{n+1}^{\pm}(u)&
=\frac{(u/v_{\mp})^2-1}{q^{-2}(u/v_{\mp})^2-q^{2}}
\Xc_{n}^{-}(v)
\eal
\een
and
\ben
\bal
\h_{n+1}^{\pm}(u)\Xc_{n-1}^{+}(v)\h_{n+1}^{\pm}(u)^{-1}&=
\frac{q^{-1}u/v_{\pm}+q}{u/v_{\pm}+1}\Xc^{+}_{n-1}(v),\\
\h_{n+1}^{\pm}(u)^{-1}\Xc_{n-1}^{-}(v)\h_{n+1}^{\pm}(u)&=
\frac{q^{-1}u/v_{\mp}+q}{u/v_{\mp}+1}\Xc^{-}_{n-1}(v),\\,
\eal
\een
while
\ben
\bal
\h_{n+1}^{\pm}(u)\Xc_{i}^{+}(v)
&=\Xc_{i}^{+}(v)\h_{n+1}^{\pm}(u),\\[0.4em]
\h_{n+1}^{\pm}(u)\Xc_{i}^{-}(v)
&=\Xc_{i}^{-}(v)\h_{n+1}^{\pm}(u),
\eal
\een
for $1\leqslant i\leqslant n-2$. For the relations involving $\Xc^{\pm}_i(u)$ we have
\ben
(u-q^{\pm (\alpha_i,\alpha_j)}v)\Xc_{i}^{\pm}(uq^i)\Xc_{j}^{\pm}(vq^j)
=(q^{\pm (\alpha_i,\alpha_j)}u-v) \Xc_{j}^{\pm}(vq^j)\Xc_{i}^{\pm}(uq^i)
\een
for $i,j=1,\dots,n-1$;
\ben
(u^2-q^{\pm (\alpha_i,\alpha_n)}v^2)\Xc_{i}^{\pm}(uq^i)\Xc_{n}^{\pm}(vq^n)
=(q^{\pm (\alpha_i,\alpha_n)}u^2-v^2) \Xc_{n}^{\pm}(vq^n)\Xc_{i}^{\pm}(uq^i)
\een
for $i=1,\dots,n$
and
\ben
[\Xc_i^{+}(u),\Xc_j^{-}(v)]=
\delta_{ij}(q_i-q_i^{-1})\Big(\delta\big(u\ts q^{-c}/v\big)\h_i^{-}(v_+)^{-1}\h_{i+1}^{-}(v_+)
-\delta\big(u\ts q^{c}/v\big)\h_i^{+}(u_+)^{-1}\h_{i+1}^{+}(u_+)\Big)
\een
for $j\neq n$
and
\ben
[\Xc_i^{+}(u),\Xc_n^{-}(v)]=
\delta_{in}(q_n-q_n^{-1})\Big(\delta\big((u\ts q^{-c}/v)^2\big)q^{c}\h_i^{-}(v_+)^{-1}\h_{i+1}^{-}(v_+)
-\delta\big((u\ts q^{c}/v)^2\big)q^{-c}\h_i^{+}(u_+)^{-1}\h_{i+1}^{+}(u_+)\Big)
\een
together with the
Serre relations
\beql{serrex}
\sum_{\pi\in \Sym_{r}}\sum_{l=0}^{r}(-1)^l{{r}\brack{l}}_{q_i}
  \Xc^{\pm}_{i}(u_{\pi(1)})\dots \Xc^{\pm}_{i}(u_{\pi(l)})
  \Xc^{\pm}_{j}(v)\tss \Xc^{\pm}_{i}(u_{\pi(l+1)})\dots \Xc^{\pm}_{i}(u_{\pi(r)})=0,
\eeq
which hold for all $i\neq j$ and we set $r=1-A_{ij}$ for $A_{ij}=-1$,
and
\beq
\bal
&\sum_{\sigma\in \mathfrak{S}_2}\sigma\Big((q^2u_1+u_2)(\Xc^{\pm}_n(v)\Xc_{n-1}^{\pm}(u_1)\Xc_{n-1}^{\pm}(u_2)\\
&-[2]_{q^2}\Xc_{n-1}^{\pm}(u_1)\Xc_{n}^{\pm}(v)\Xc_{n-1}^{\pm}(u_2)
+\Xc_{n-1}^{\pm}(u_1)\Xc_{n-1}^{\pm}(u_2)\Xc_{n}^{\pm}(v))\Big)=0.
\eal
\eeq
\eth

Finally we arrive at our main result. 
\bth
The mapping $DR:U^{ext}_q(A^{(2)}_{2n-1})\to U_q(R)$: $q^{c/2}\mapsto q^{c/2}$ and
\ben
\bal
X^{\pm}_{j}(u)\mapsto X^{\pm}_j(u),\qquad\quad h^{\pm}_{i}(u)\mapsto  h^{\pm}_{i}(u),
\eal
\een
for $j=1,\dots,n$, $i=1,\dots n+1$,
defines an isomorphism $U^{ext}_q(A^{(2)}_{2n-1})\to U_q(R)$.
\eth

\bpf
The relations between $U(R)$ and $U(\overline{R})$ and Theorem \ref{thm:Ubar}
imply the mapping $DR$ is a surjective homomorphism. On the other hand Proposition
\ref{prop:L-op}
  implies a homomorphism
\beq
\bal
RD: & U(R)\to U^{ext}_q(A^{(2)}_{2n-1})\\
&L^{\pm}(u)\mapsto L^{\pm}(u).
\eal
\eeq
Furthermore, we have $DR \circ RD =ID$ by using Theorem \ref{thm:tri-dec}.
Thus the homomorphism $DR$ is also injective.

\epf

\centerline{\bf Acknowledgments}
The work is partially supported by the Simons Foundation (grant no. 523868) and the National Natural Science Foundation of China (grant no. 11531004).

\end{document}